\documentclass[a4paper,12pt,fleqn]{smfart}
\usepackage[latin1]{inputenc}
\usepackage[francais,english]{babel} \AutoSpaceBeforeFDP
\usepackage[T1]{fontenc}
\usepackage{layout}
\usepackage{amstext,amssymb,amscd}
\usepackage{smfthm}
\usepackage{graphicx}
\usepackage{color}
\usepackage{epsfig}
\usepackage[all]{xy}
\usepackage{enumerate}
\usepackage{version}
\usepackage{import}
\usepackage{calc}

\newcommand{\sss}{\mathrm{SL}}
\newcommand{\e}{\textrm{ }}

\newcommand{\R}{\mathbb{R}}

\newcommand{\Z}{\mathbb{Z}}

\newcommand{\G}{\mathcal{G}}

\newcommand{\F}{\mathcal{F}}
\newcommand{\A}{\mathcal{A}}

\renewcommand{\O}{\Omega}

\renewcommand{\G}{\Gamma}

\newcommand{\PP}{\mathbb{P}}

\newcommand{\Quo}{\Omega/\!\raisebox{-.65ex}{\ensuremath{\Gamma}}}

\newcommand{\GG}{\mathcal{G}}

\theoremstyle{definition}

\theoremstyle{remark}
\newtheorem*{thm}{Théorème}

\title{Espace des modules de certains polyèdres projectifs miroirs}
\author{Ludovic Marquis}
\date{}
\email{ludovic.marquis@math.u-psud.fr}
\urladdr{www.math.u-psud.fr/~marquis}

\begin{altabstract}
A projective mirror polyhedron is a projective polyhedron endowed
with reflections across its faces. We construct an explicit
diffeomorphism between the moduli space of a mirror projective
polyhedron with fixed dihedral angles in $(0,\frac{\pi}{2}]$, and the union of $n$ copies
of $\R^d$, when the polyhedron has the combinatorics of an \emph{ecimahedron}, an
infinite class of combinatorial polyhedra we introduce here. Moreover, the integers $n$ and $d$ can be computed explicitly in terms of the combinatorics and the fixed dihedral angles.
\end{altabstract}

\begin{document}

\newenvironment{nota}{\begin{enonce}[remark]{Notations}}{\end{enonce}}

\maketitle
\section{Introduction}

\par{
Un ouvert convexe $\Omega$ de l'espace projectif réel
$\mathbb{P}^n(\R)$ est dit $\emph{divisible}$ lorsqu'il existe un
sous-groupe discret $\Gamma$ du groupe des transformations
projectives $\mathrm{PGL}_{n+1}(\R)$ qui préserve $\Omega$ et agit
proprement et cocompactement sur $\Omega$. On dit aussi que
$\Gamma$ $\emph{divise}$ $\Omega$. Vinberg a trouvé une méthode
géométrique pour construire des convexes divisibles à l'aide des
groupes de Coxeter. Cette méthode a été initiée par Poincaré pour
construire des réseaux du groupe des isométries de l'espace
hyperbolique réel.
}
\\
\par{
Dans la méthode de Poincaré, on commence par se donner un polyèdre
$P$ de l'espace hyperbolique réel $\mathbb{H}^n$ et une famille
d'isométries identifient 2 à 2 les $(n-1)$-faces de $P$. Si ces
isométries vérifient des relations de compatibilité, qui disent
que l'on pave bien autour de chaque $(n-2)$-face, alors le groupe
engendré par ces isométries agit proprement sur $\mathbb{H}^n$ et
$P$ est un domaine fondamental pour cette action. En particulier,
si $P$ est un polyèdre de $\mathbb{H}^n$ dont les angles entre les
$(n-1)$-faces adjacentes sont des sous-multiples de $\pi$, alors
le groupe engendré par les réflexions qui fixent les faces de $P$
est un groupe de Coxeter; il agit proprement sur $\mathbb{H}^n$ et
$P$ est un domaine fondamental.
}
\\
\par{
Vinberg a généralisé cette méthode au cadre projectif de la façon
suivante. On se donne un polyèdre compact $P$ de
$\mathbb{P}^n(\R)$ et des réflexions qui fixent les faces de $P$.
Cette fois-ci, on a des degrés de liberté pour choisir ces
réflexions. L'hypothèse  \og les angles dièdres de $P$ sont des sous-multiples de $\pi$\fg $\,$  devient \og le produit de deux réflexions par rapport à des faces adjacentes est conjugué à une rotation d'ordre fini \fg. Alors, le groupe $\Gamma$ engendré par ces réflexions est un groupe de Coxeter, qui agit
proprement sur un certain ouvert convexe $\Omega$ de $\mathbb{P}^n(\R)$
qui contient l'intérieur de $P$, et $P \cap \Omega$ est un domaine fondamental pour l'action
de $\Gamma$ sur $\Omega$.
}
\\
\par{
Cette construction motive l'introduction de la notion de polyèdre
projectif miroir. Il s'agit de la donnée d'un polyèdre projectif
$P$ et de réflexions par rapport aux faces de $P$. On a alors une
notion naturelle d'angle dièdre (on donnera les définitions
précises à la partie \ref{section}). Dans ce texte, on
s'intéresse au problème suivant: étant donné le graphe d'un
polyèdre $\GG$ dont les arêtes sont étiquetées par des réels $\theta
\in ]0,\frac{\pi}{2}]$, on cherche à comprendre l'espace des modules $X_{\GG}$ des
polyèdres projectifs miroirs qui ont la combinatoire de $\GG$, et dont les angles dièdres sont les angles fixés par les étiquettes de $\GG$.
}
\\
\par{
Choi s'est intéressé à ce problème dans \cite{Choi2}, sous un angle
un peu différent: il fixe un polyèdre projectif $P$ dont il
étiquette les arêtes par des réels $\theta \in ]0,\frac{\pi}{2}]$,
et il cherche à comprendre l'espace des réflexions par rapport aux
faces de $P$ qui font de $P$ un polyèdre projectif miroir dont les
angles dièdres sont les étiquettes des arêtes de $P$. Choi montre
sous une hypothèse dite \og d'ordonnabilité \fg $,$ que cet espace des
modules est une variété lisse dont il sait calculer la dimension. Il faut noter que l'hypothèse d'ordonnabilité de Choi porte non
seulement sur la combinatoire du polyèdre $P$ mais aussi sur les
étiquettes des arêtes de $P$.
}
\\
\par{
Pour pouvoir décrire la topologie de $X_{\GG}$, on aura besoin de faire une hypothèse sur la combinatoire du graphe de polyèdre $\GG$. On va définir une sous-classe de graphes de polyèdres: les écimaèdres. De façon imagée, il s'agit des polyèdres obtenus à partir du tétraèdre et par une suite finie de coupes "près d'un sommet" (voir définition \ref{defecim}). On montrera le théorème suivant:
}
\begin{thm}
[\ref{theo}]
Soit $\GG$ un écimaèdre étiqueté. L'espace $X_{\GG}$ est vide ou difféomorphe à la réunion de $n$ copies d'un certain $\R^d$, où les entiers $n$ et $d$ s'expriment à l'aide de la combinatoire de $\GG$ et de ses étiquettes.
\end{thm}
\par{
On donnera un critère précis et simple pour savoir si $X_{\GG}$ est vide ou non. On donnera aussi un système de coordonnées explicite sur l'espace $X_{\GG}$ (théorème \ref{final}).
}
\\
\par{
Rappelons que l'objectif est de construire des convexes divisibles. Il est donc essentiel de savoir qu'il existe une infinité d'écimaèdres étiquetés permettant grâce à la méthode de Vinberg de construire des convexes divisibles. On montrera aussi le théorème suivant qui sera un corollaire évident du théorème \ref{theo}.
}
\begin{thm}[Corollaire \ref{corohyp}]
Soient $P$ un écimaèdre hyperbolique de $\mathbb{H}^3$ compact dont tous les angles dièdres sont aigus ou droits, et $\GG$ le graphe de $P$ étiqueté par les angles dièdres du polyèdre $P$, alors, l'espace des modules $X_{\GG}$ des polyèdres projectifs miroirs qui réalisent $\GG$ est difféomorphe à $\R^{e_+ -3}$, où $e_+$ est le nombre d'arêtes de $P$ dont l'angle dièdre est différent de $\frac{\pi}{2}$.
\end{thm}
\par{
Terminons cette introduction en donnant le plan de cet article. La deuxième partie est constituée de rappels et de définitions. On rappelle les définitions de polyèdre projectif et de groupe de Coxeter. Ensuite, on définit la notion de polyèdre projectif miroir qui est l'objet central de cet article. Pour motiver cette définition, on rappelle le théorème de Vinberg. Pour calculer les entiers $n$ et $d$ du théorème \ref{theo}, on introduira la notion de 3-circuits de $\GG$, il s'agit des triplets $(r,s,t)$ de faces fermées de $\GG$ mutuellement adjacentes. Le théorème d'Andreev, qui est un analogue hyperbolique du théorème principal et dont on rappellera un énoncé, pousse à définir les notions de 3-circuits \emph{prismatiques}(i.e $r\cap s \cap t = \varnothing$), \emph{non prismatiques} , mais aussi les 3-circuits \emph{sphériques}, \emph{affines} et \emph{hyperboliques} selon la somme des angles dièdres $\theta_{r \cap s}$, $\theta_{s \cap t}$ et $\theta_{t \cap r}$; et enfin de 3-circuit \emph{avec angle droit} ou \emph{sans angle droit} selon que l'un des angles $\theta_{r \cap s}$, $\theta_{s \cap t}$ ou $\theta_{t \cap r}$ est droit ou non. Toutes ces notions sont cruciales dans l'étude de la topologie de $X_{\GG}$.
}
\\
\par{
Dans la troisième partie, on commence par présenter la notion d'écimaèdre qui sera la seule hypothèse du théorème \ref{theo}. On montre que les \emph{écimaèdres} peuvent être obtenus en recollant des \emph{blocs fondamentaux} le long de faces triangulaires. Un bloc fondamental est un tétraèdre tronqué en un ou plusieurs sommets distincts. Ensuite, on consacre un paragraphe à la présentation du théorème \ref{theo}.
}
\\
\par{
Les paragraphes à partir du paragraphe \ref{debutdemo} jusqu'au paragraphe \ref{3} constituent la démonstration du théorème \ref{theo}. On commencera par montrer les points les plus simples du théorème. Au paragraphe \ref{1}, on montre que l'espace des modules du triangle est homéomorphe à $\R$ (ou à un point s'il y a un angle droit), pour cela on introduit l'invariant $R$ d'un triangle miroir. Cet invariant sera la clé de voûte de la paramétrisation de l'espace $X_{\GG}$.
}
\\
\par{
Dans le paragraphe \ref{2}, on construit pour tout 3-circuit prismatique $(r,s,t)$ un plan canonique coupant $P$ à l'intérieur des arêtes communes à $r,s,t$. Ce plan permet de découper géométriquement le polyèdre $P$, en blocs fondamentaux que l'on comprendra par la suite.
}
\\
\par{
Lorsque le 3-circuit est prismatique, sans angle droit, et affine ou sphérique, l'invariant $R$ ne peut pas à appartenir à un certain intervalle compact. C'est cette obstruction qui est la cause de la non connexité de l'espace $X_{\GG}$. Dans les paragraphes \ref{5}, \ref{6} et \ref{orien}, on cherche à ramener le comptage des composantes connexes de $X_{\GG}$ à un problème combinatoire sur des forêts.
}
\\
\par{
Enfin, dans  la partie \ref{4}, on calcule l'espace des modules pour le tétraèdre et pour les blocs fondamentaux. Dans la partie \ref{3} on décrit la manière de recoller deux blocs fondamentaux.
}
\\
\par{
On termine la quatrième partie en explicitant au paragraphe \ref{last} un système de coordonnées sur l'espace des modules des polyèdres projectifs miroirs.
}
\\
\par{
Enfin, la cinquième partie décrit des exemples concrets d'espaces des modules de polyèdres projectifs miroirs.
}
\\
\par{
Je tiens à remercier Yves Benoist pour m'avoir guidé pendant ce
travail, mais aussi pour m'avoir fait découvrir les convexes
divisibles. Je remercie aussi Mickaël Crampon pour ses conseils de
rédaction. Enfin, je remercie \og le rapporteur anonyme \fg $\,$  dont les conseils m'ont permis d'améliorer grandement la qualité de ce texte.
}
\section{Notion de polyèdre projectif miroir}\label{section}

Commençons par donner des définitions précises et par rappeler le
théorème de Vinberg. On en profitera pour rappeler le
théorème d'Andreev dont le théorème \ref{theo} est un analogue naturel.

\subsection{Les convexes de $\mathbb{P}^+(V)$}

Soit $V$ un espace vectoriel réel de dimension finie $n+1\geqslant
3$. Notons $\mathbb{P}^+(V) = \{ \textrm{demi-droites vectorielles de } V\} = (V-\{0\})/_{\R_+^*}$ la sphère
projective; c'est une variété projective difféomorphe à la
$n$-sphère euclidienne usuelle. Le groupe des transformations
projectives de $\mathbb{P}^+(V)$ est $\mathrm{SL}^{\pm}(V) = \{ u
\in \mathrm{GL}(V) \textrm{ tel que } \det(u)=\pm 1\} \simeq
\mathrm{GL}(V)/_{\R_+^*}$.

\begin{defi}
Une partie $\Omega$ de $\mathbb{P}^+(V)$ est dite \emph{convexe}
lorsque son intersection avec tout grand cercle de $\mathbb{P}^+(V)$
est connexe. On dit qu'elle est $\emph{proprement convexe}$ s'il
existe en plus une hypersphère projective de $\mathbb{P}^+(V)$ qui ne
rencontre pas $\overline{\Omega}$, l'adhérence de $\Omega$, ce qui
est équivalent à l'existence d'une carte affine dans laquelle
$\Omega$ est un convexe borné. L'intérieur d'une partie $\Omega$ de $\mathbb{P}^+(V)$ sera noté $\mathring{\Omega}$. Un convexe $\Omega$ est dit \emph{strictement
convexe} si son bord $\partial \Omega = \overline{\Omega}
\backslash \mathring{\Omega}$ ne contient pas de segments non triviaux. On
dit qu'un ouvert convexe est $\emph{divisible}$ s'il existe un
sous-groupe discret $\Gamma$ de $\mathrm{SL}^{\pm}(V)$ qui
préserve $\Omega$, agit proprement sur $\Omega$ et tel que le
quotient $\Quo$ soit compact. On dit aussi que
$\Gamma$ $\emph{divise}$ $\Omega$.
\end{defi}

Pour en savoir plus sur la théorie des convexes divisibles, on pourra consulter la série d'articles de Benoist: \cite{Beno9, Beno8, Beno6, Beno2}, pour une étude complète de la dimension 2 on pourra consulter l'article de Goldman: \cite{Gold1}.

Donnons tout de suite des exemples de convexes divisibles:
\begin{itemize}
\item La sphère projective est un convexe divisible.
\newline
\item N'importe quel espace affine de dimension $n$ inclus dans
$\mathbb{P}^+(\R^{n+1})$ est divisé par le groupe $\Z^n$; $\R^n$ est
donc un exemple de convexe divisible non proprement convexe.
\newline
\item On peut construire un autre convexe divisé par $\Z^n$ de la
façon suivante. La base canonique de $\R^{n+1}$ définit
naturellement un pavage de $\mathbb{P}^+(\R^{n+1})$ en $2^{n+1}$ simplexes de dimension $n$. Si on note $\Omega$ l'intérieur de l'un deux, alors la composante neutre du stabilisateur de $\Omega$ dans $\sss_{n+1}(\R)$ est le groupe $D$ des matrices diagonales à diagonale strictement positive de
déterminant 1. Le groupe $D$ agit proprement et simplement transitivement
sur $\Omega$, et l'image $\Gamma$ de toute représentation fidèle et discrète de $\Z^n$ dans $D$ est un réseau cocompact de $D$. Ainsi, $\Gamma$ divise $\Omega$ qui est donc un
convexe divisible proprement convexe mais non strictement convexe.
\newline
\item Terminons cette liste d'exemples par la construction d'un
convexe divisible strictement convexe. Soient $q$ une forme
quadratique sur $V$ de signature $(n,1)$, et $\Omega$ l'une des 2
composantes connexes de l'ouvert $\{ [v]\in \mathbb{P}^+(V) \, | \,
q(v) < 0 \}$. Il s'agit du modèle projectif de l'espace
hyperbolique réel $\mathbb{H}^n$: il est donc divisé par tous les
réseaux cocompacts de $\textrm{Isom}(\mathbb{H}^n)$, le groupe des
isométries de $\mathbb{H}^n$.
\end{itemize}

On va construire dans ce texte des convexes divisibles
$\Omega$, divisés par des groupes de Coxeter $W$, selon une méthode
initiée par Vinberg.

\subsection{Les polyèdres projectifs}

On notera $p:V\setminus \{ 0 \} \rightarrow \mathbb{P}^+(V)$ la projection naturelle.

\begin{defi}\label{proj}
Un \emph{polyèdre projectif} est un fermé proprement
convexe $P$ d'intérieur non vide de $\mathbb{P}^+(V)$ tel qu'il
existe un nombre fini de formes linéaires $\alpha_1,...,\alpha_f $
sur $V$ telles que $P = p(\{ x \in V\setminus \{ 0 \} \,|\e \alpha_i(x) \leqslant 0,\e i=1...f\})$.
\end{defi}

Pour tout polyèdre projectif $P$ de $\mathbb{P}^+(V)$, on définit une relation
d'équivalence $\thicksim_P$ sur $P$ de la façon suivante:
$$x\thicksim_P y \e \Leftrightarrow \textrm{ le segment } [x,y] \e \textrm{se prolonge strictement des deux côtés dans } P$$
Les adhérences des classes d'équivalences de $\thicksim_P$
s'appellent les $\emph{faces}$ de $P$. On
appelle \emph{dimension} (resp. \emph{codimension}) d'une face la
dimension (resp. codimension) de l'unique sous-espace projectif minimal la
contenant.

\subsection{Groupes de Coxeter}

Les groupes de Coxeter seront au c{\oe}ur de nos motivations; nous
allons donc rappeler quelques définitions.

\begin{defi}
Un \emph{système de Coxeter} est la donnée d'un ensemble fini $S$
et d'une matrice symétrique $M=(M_{st})_{s,t \in S}$ telle que les
coefficients diagonaux vérifient $M_{ss}=1$ et les coefficients
non diagonaux vérifient $M_{st} \in \{2,3,...,\infty \}$. Le
cardinal de $S$ s'appelle le \emph{rang} du système de Coxeter
$(S,M)$. \`{A} un système de Coxeter, on associe un \emph{groupe de
Coxeter} $W_S$. Il s'agit d'un groupe défini par générateurs et
relations. Les générateurs sont les éléments de $S$ et on impose
les relations $(st)^{M_{st}}=1$ pour $s,t \in S$ tels que $M_{st}
\neq \infty$.
\end{defi}

Avec toute partie $S'$ de $S$, on peut former le groupe de Coxeter
$W_{S'}$ associé au système de Coxeter $(S',M')$, où $M'$ est la
restriction de $M$ à $S'$. Un corollaire du théorème de Vinberg
montre que le morphisme naturel $W_{S'} \rightarrow W_S$ est
injectif. Ainsi, $W_{S'}$ peut être identifié avec le sous-groupe
de $W_S$ engendré par la partie $S'$. On utilisera donc la
notation $W_{S'}$ pour désigner ces deux groupes.

\subsection{Polyèdre projectif miroir}

Dans le cadre projectif, si on se donne un polyèdre projectif, on
peut choisir pour chaque face une réflexion qui la préserve. On va utiliser cette liberté pour construire de nombreux
convexes divisibles associés à un groupe de Coxeter.

\subsubsection{Le théorème de Vinberg}\label{vinb}

Le but de ce paragraphe est de rappeler l'énoncé du théorème de Vinberg dont on peut trouver une démonstration dans \cite{Vin3} et dans \cite{Beno7}.

Soit $P$ un polyèdre projectif; on note $S$ l'ensemble des faces
de codimension 1 de $P$, et on suppose que $P= p(\{ x \in V \setminus \{ 0 \} \,|\e \alpha_s(x)
\leqslant 0,\e s \in S \})$. En particulier, pour toute face $s \in S$ et pour tout élément $x \in s$, on a $\alpha_s(x)=0$. On se donne pour chaque face $s$ de
codimension 1, une réflexion projective $\sigma_s = Id -
\alpha_s\otimes v_s$ avec $v_s \in V$ et $\alpha_s(v_s)=2$ qui
fixe la face $s$, et on note $\Gamma$ le groupe engendré par les
réflexions $\sigma_s$ pour $s\in S$, et $v_{ts} = -\alpha_s(v_t)$.

Une étude locale autour des faces de codimension 2 montre que des
conditions nécessaires pour que les $\gamma(P)$ pavent une partie
de $\mathbb{P}^+(V)$ (c'est à dire que les $\gamma(P)$ soient
d'intérieur disjoints) sont les suivantes:

$$
\begin{array}{l}
\forall s,t \in S, \textrm{ tels que } \textrm{codim}(s \cap t) =
2,
\textrm{ on a:}\\
\left\{
\begin{array}{lcccc}
1) & v_{ts} \geqslant 0  & \textrm{ et } &  & (v_{ts}=0
\Leftrightarrow
v_{st}=0)\\
2)a) & v_{st}v_{ts} \geqslant 4 & \textrm { ou } & b) &
v_{st}v_{ts} =
4\cos^2\big(\frac{\pi}{m_{st}}\big)\\
& & & &\textrm{ avec } m_{st}\geqslant 2 \textrm{ entier}
\end{array}
\right.
\end{array}
$$

Le théorème de Vinberg affirme que ces conditions sont en fait
suffisantes.

On va relâcher légèrement la condition (2) en la condition
(2\textquoteright) qui autorise les angles à prendre toutes les valeurs dans
$]0,\frac{\pi}{2}]$.

$$
\begin{array}{ccccc}
2 \textrm{\textquoteright} ) a) & v_{st}v_{ts} \geqslant 4 & \textrm { ou } & b) &
v_{st}v_{ts} =
4\cos^2(\theta_{st})\\
& & & &\textrm{ avec } \theta_{st} \in ]0, \frac{\pi}{2}]
\end{array}
$$

\begin{defi}\label{miroir}
Un \emph{polyèdre projectif miroir} est la donnée d'un polyèdre
projectif $P$ et, pour chaque face $s$ de $P$, d'une réflexion
$\sigma_s = Id- \alpha_s \otimes v_s$ qui fixe $s$, avec les
conventions suivantes: $P= p(\{ x \in V \setminus \{ 0 \} \,|\e \alpha_s(x)
\leqslant 0,\e s \in S \})$,
$\alpha_s(v_s) = 2$ et les conditions (1) et (2\textquoteright $b$) sont vérifiées.
On dit que $\emph{l'angle dièdre}$ entre deux faces $s$ et $t$
telles que  $\textrm{codim}(s \cap t) = 2$ est $\theta_{st} \in
]0,\frac{\pi}{2}]$ si $v_{st}v_{ts} = 4\cos^2(\theta_{st})$.
\end{defi}

\begin{rema}
Il est important de noter que l'on exclut le cas (2\textquoteright $a$) de la définition de polyèdre projectif miroir. En effet, la seule définition d'angle dièdre acceptable pour le cas (2\textquoteright $a$) serait un angle nul, et ce cas ne nous intéresse pas.
\end{rema}

\begin{nota}\label{equa} Soit $P$ un polyèdre miroir, on désignera par la
lettre $S$ l'ensemble de ses faces de codimension 1. On s'est
donné pour chaque $s \in S$ une réflexion $\sigma_s = Id- \alpha_s
\otimes v_s$ telle que $\alpha_s(v_s)=2$ et $P= p(\{ x \in V \setminus \{ 0 \} \,|\e \alpha_s(x)
\leqslant 0,\e s \in S \})$. Cette convention de signe détermine
un unique couple $([\alpha_s],[v_s]) \in \mathbb{P}^+(V^*) \times
\mathbb{P}^+(V)$. On appellera le point $[v_s]$ la $\emph{polaire}$
de la face $s$. On fera bien attention au fait que le couple
$(\alpha_s,v_s)$ n'est pas unique. En effet, si $(\alpha_s,v_s)$
convient alors pour tout $\lambda \in \R^*_+$ le couple
$(\lambda^{-1} \alpha_s, \lambda v_s)$ convient aussi. Enfin, si
$s,t \in S$ alors on notera $\mu_{st}$ la quantité
$v_{ts}v_{st} = \alpha_s(v_t)\alpha_t(v_s)$ qui est bien définie indépendamment de $\lambda$. Par définition,
si $s$ et $t$ sont deux faces de $P$ telles que $ \textrm{codim}(s
\cap t) = 2$ alors $\mu_{st} = 4 \cos^2(\theta_{st})$; dans ce
cas on notera $m_{st} = \frac{\pi}{\theta_{st}}$. Ces notations
et conventions seront utilisées tout au long de ce texte.
\end{nota}

\begin{defi}
Soit $P$ un polyèdre projectif miroir dont les angles dièdres
sont des sous-multiples de $\pi$. Le \emph{système de Coxeter
associé à $P$} est le système de Coxeter $(S,M)$, où $S$ est
l'ensemble des faces de codimension 1 de $P$ et pour tous $s,t \in
S$, on a $M_{st}=m_{st}$ si les faces $s$ et $t$ vérifient
$\textrm{codim}(s \cap t) = 2$ et $M_{st} =\infty$ sinon. On note
$W_S$ le groupe de Coxeter associé au système $(S,M)$.
\end{defi}

\begin{theo}[Vinberg]
Soit $P$ un polyèdre projectif miroir dont les angles dièdres
sont des sous-multiples de $\pi$. Soient $(S,M)$ le système de
Coxeter associé à $P$, $W_S$ le groupe de Coxeter associé et
$\Gamma$ le groupe engendré par les réflexions projectives $(\sigma_s)_{s \in S}$. Alors,

a) Les polyèdres $\gamma(P)_{\gamma \in \Gamma}$ pavent un convexe
$\Omega$ de $\mathbb{P}^+(V)$.

b) Le morphisme $\sigma:W_S \rightarrow \Gamma$ défini par $\sigma(s) =
\sigma_s$ est un isomorphisme.

c) Le groupe $\Gamma$ est un sous-groupe discret de $\mathrm{SL}^{\pm}(V)$.

d) Le groupe $\Gamma$ agit proprement sur $\mathring{\Omega}$, l'intérieur
de $\Omega$.

e) L'ensemble $\Omega$ est ouvert si et seulement si pour tout sommet $v$ de
$P$, $W_{S_v}$ est fini, où $S_v=\{s \in S \, | \, v \subset s
\}$.
\end{theo}


\begin{rema}
Lorsque le graphe de Coxeter de $W_S$ est connexe, que $W_S$ est
infini et que les $v_s$ engendrent $V$, alors $\Omega$ est
proprement convexe (Cela est démontré dans \cite{Beno7}).
\end{rema}

\begin{rema}\label{valence3}
Le point $e)$ du théorème de Vinberg et la classification des
groupes de Coxeter finis (\cite{Bou}) montrent que si $\Omega$ est ouvert et le polyèdre $P$ est de dimension 3, alors pour tout sommet $v$ de $P$ le groupe $W_{S_v}$ est de type
$(2,2,n)$ avec $n \geqslant 2$, ou $(2,3,3)$, $(2,3,4)$ ou encore
$(2,3,5)$. En particulier, tout sommet de $P$ doit être de valence
3 pour obtenir un convexe divisible $\Omega$.
\end{rema}

\subsubsection{Combinatoire d'un polyèdre de $\mathbb{P}^+(\R^4)$}

On se restreint désormais à la dimension 3 et on appellera une
face de codimension 1 (resp. 2, resp. 3) une face (resp. une arête,
resp. un sommet). Certaines propriétés combinatoires des polyèdres
seront essentielles, nous allons donc donner quelques définitions.

A tout polyèdre $P$, on associe un graphe $\GG_P$ planaire et
3-connexe (i.e $\GG_P$ a plus de quatre sommets et $\GG_P$ privé de
2 sommets quelconques non adjacents est encore connexe) dont les
sommets sont les sommets de $P$ et les arêtes les arêtes de $P$.

\begin{defi}
Soient $P$ un polyèdre et $\GG$ un graphe, on dira que $P$
\emph{réalise} $\GG$ lorsque $\GG$ et $\GG_P$ sont des graphes planaires isomorphes.
\end{defi}

\begin{rema}
En fait, on peut montrer qu'un graphe est réalisé par un polyèdre
si et seulement s'il est planaire et 3-connexe (Théorème de
Steinitz, voir le livre \cite{JG}). De plus, le plongement d'un tel graphe dans le plan est unique à isotopie et réflexion près.
\end{rema}

Ce qui nous intéresse, c'est de réaliser des polyèdres avec des
angles dièdres prescrits; on en vient donc à la définition
suivante:

\begin{defi}
Un \emph{graphe étiqueté} est la donnée d'un graphe de valence 3,
planaire et 3-connexe $\GG$ et pour chaque arête $e$ de $\GG$ d'un
réel $\theta_e \in ]0,\frac{\pi}{2}]$.
\end{defi}

\begin{rema}
Un graphe étiqueté $\GG$ est en particulier un graphe planaire et
3-connexe, la notion de face de $\GG$ est donc bien définie.
\end{rema}

\begin{nota}
Soit $\GG$ un graphe étiqueté, si $e$ est une arête de $\GG$, on
désignera par $\theta_e$ l'angle associé. De plus, on désignera
par $\mu_e$ le réel $4\cos^2(\theta_e)$; et si $\theta_e =
\frac{\pi}{m}$, où $m$ est un entier supérieur ou égale à 2, on dira que l'arête est
d'ordre $m$. De plus, si $s$ et $t$ désignent deux faces de $\GG$
qui partagent une arête alors cette arête sera notée $st$, et
$\theta_{st}$ désigne alors l'angle qu'elle porte.
\end{nota}

\begin{defi}\label{marqué}
Un graphe étiqueté (resp. polyèdre miroir) est dit $\emph{marqué}$
lorsqu'une numérotation de ses faces a été choisie.
\end{defi}

\begin{defi}\label{realisation}
Soient $\GG$ un graphe étiqueté marqué, et $P$ un polyèdre miroir
marqué; on dit que $P$ réalise $\GG$ lorsque le polyèdre
sous-jacent à $P$ $\emph{réalise}$ le graphe sous-jacent à $\GG$
via une identification qui respecte leur marquage et que les
angles dièdres du polyèdre projectif miroir $P$ correspondent aux
étiquettes de $\GG$.
\end{defi}

\begin{rema}
Pour alléger la rédaction, on supposera implicitement tout au long
de ce texte que tous les graphes étiquetés et tous les polyèdres
projectifs miroirs sont marqués.
\end{rema}

\begin{defi}
Soit $\GG$ un graphe planaire et 3-connexe, un \emph{$k$-circuit orienté} (resp. \emph{$k$-circuit}) $\Gamma$
de $\GG$ est une suite $(f_1,e_1,f_2,e_2,...,f_{k},e_{k})$ définie
à permutation circulaire près (resp. à permutation circulaire près
et à sens de parcours près) telle que:

$$
\left\{
\begin{array}{l}
\textrm{les } f_i \textrm{ sont des faces distinctes de } \GG,\\
\textrm{les } e_i \textrm{ sont des arêtes de } \GG,\\
\forall i = 1...k, \, f_{i} \cap f_{i+1} = e_i, \textrm{ où }
f_{k+1}=f_1.
\end{array}
\right.
$$

De plus, si toutes les extrémités des arêtes de $\Gamma$ sont
distinctes, alors on dit que $\Gamma$ est \emph{prismatique}.
\end{defi}

\begin{rema}
On voit aisément que si $\GG$ est un graphe planaire, 3-connexe et de valence 3 alors les trois arêtes de tout 3-circuit non prismatique ont un sommet en commun.
\end{rema}

\begin{defi}
Soient $\GG$ un graphe étiqueté et $\Gamma$ un $k$-circuit de $\GG$, on
note $(\theta_i)_{i=1...k}$ les angles des arêtes de $\Gamma$;
posons $\Sigma = \displaystyle{\sum_{i=1...k}{\theta_i}}$.

$$
\textrm{ On dira que } \Gamma \textrm{ est }
\left\{
\begin{array}{lll}\label{droithyp}
\emph{sans angle droit} & \textrm{ si les } \theta_i  \textrm{
sont tous différents
de } \frac{\pi}{2} \textrm{ pour } i=1...k,\\

\emph{avec angle
droit} & \textrm{ si l'un au moins des } \theta_i  \textrm{ est égal à } \frac{\pi}{2} \textrm{ pour } i=1...k,\\

\emph{sphérique
} & \textrm{ lorsque } \Sigma > (k-2)\pi,\\

\emph{affine } &
\textrm{ lorsque } \Sigma = (k-2)\pi,\\

\emph{hyperbolique } & \textrm{ lorsque } \Sigma < (k-2)\pi.
\end{array}
\right.
$$
\end{defi}

\subsubsection{Le théorème d'Andreev}

\begin{figure}[!h]
\begin{center}
\includegraphics[width=5cm]
{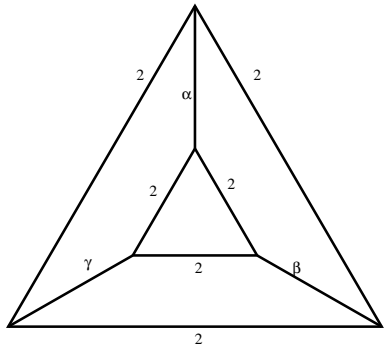} \caption{Prisme exceptionnel, avec $\alpha,\beta,\gamma \in ]0, \frac{\pi}{2}]$} \label{prisme}
\end{center}
\end{figure}

L'étude des polyèdres hyperboliques et des polyèdres projectifs miroirs fait apparaître une famille de graphes étiquetés qui nécessite un traitement à part. Il s'agit des graphes étiquetés de la figure \ref{prisme} qu'on appellera les \emph{prismes exceptionnels} et on les notera
$\GG_{\alpha,\beta,\gamma}$.

\begin{theo}[Andreev]
Soit $\mathcal{G}$ un graphe étiqueté qui n'est pas le graphe d'un
tétraèdre. Alors, il existe un polyèdre compact hyperbolique $P$
qui réalise $\mathcal{G}$ si et seulement si les quatre conditions
suivantes sont vérifiées:

\begin{itemize}
\item Tout 3-circuit non prismatique de $\GG$ est sphérique.

\item Tout 3-circuit prismatique de $\GG$ est hyperbolique.

\item Tout 4-circuit prismatique de $\GG$ est hyperbolique.

\item $\GG$ n'est pas un prisme exceptionnel.

\end{itemize}
De plus, ce polyèdre est unique à isométrie près.
\end{theo}

On peut trouver une démonstration du théorème d'Andreev dans
\cite{And} ou \cite{RHD}. On va s'intéresser à présent à un
analogue projectif de ce résultat.

\subsubsection{Espace des modules d'un polyèdre projectif miroir}

Soit $\GG$ un graphe étiqueté marqué; on introduit l'espace suivant:
$$Y_{\GG} = \{ P \textrm{ polyèdre projectif miroir marqué tel que } P
\textrm{ réalise } \GG \}.$$

Le groupe $\sss^{\pm}_4(\R)$ agit naturellement sur $Y_{\GG}$. On souhaite
comprendre l'espace quotient $X_{\GG}= Y_{\GG}/\!\raisebox{-.65ex}{\ensuremath{\sss^{\pm}_4(\R)}}$ que l'on appelle $\emph{espace des modules des polyèdres
projectifs miroirs marqués qui réalisent} \, \, \GG$.

\section{Présentation des résultats}

Dans ce texte, on s'intéresse exclusivement à une classe
très particulière de polyèdres, que nous allons définir et étudier
dans toute la partie \ref{ecima}; on appellera ces polyèdres les
écimaèdres.

\subsection{Les écimaèdres combinatoires}\label{ecima}

\begin{figure}[!h]
\begin{center}
\includegraphics[width=10cm]
{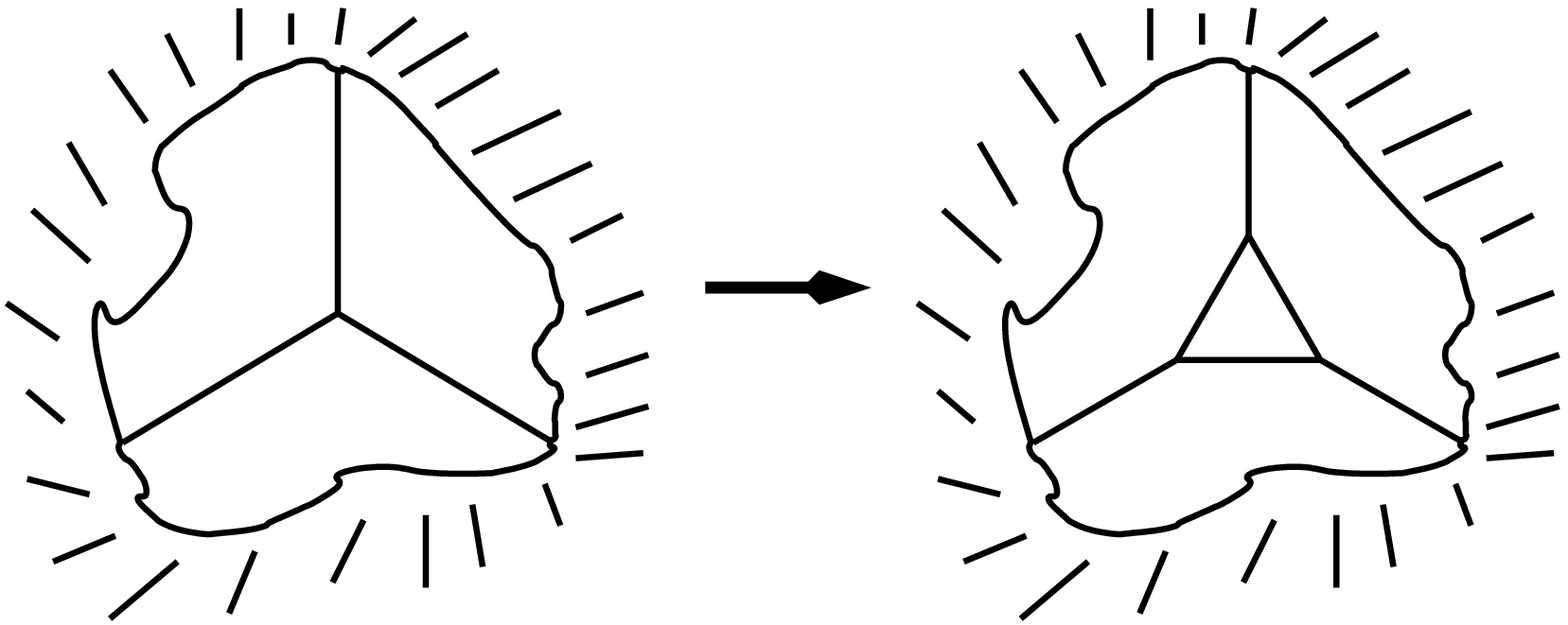} \caption{\'{E}cimage} \label{ecimage}
\end{center}
\end{figure}

Soient un graphe planaire 3-connexe $\GG$ de valence 3 et un
sommet $v$ de $\GG$; on définit un nouveau graphe $\GG'$ planaire
3-connexe de valence 3 où le sommet $v$ a été remplacé par
un triangle dont les sommets appartiennent aux arêtes issues de
$v$, comme sur la figure \ref{ecimage}. Cette opération sera
appelée l'$\emph{écimage}$ du graphe $\GG$ au sommet $v$.

\begin{defi}\label{defecim}
Un graphe $\GG$ est un $\emph{écimaèdre combinatoire}$ lorsqu'on
peut l'obtenir à partir du graphe du tétraèdre combinatoire (i.e
le graphe complet à 4 sommets) et d'un nombre fini d'écimages. Un
polyèdre $P$ est un $\emph{écimaèdre}$ si son graphe $\GG_P$ est
un écimaèdre combinatoire.
\end{defi}

\begin{defi}
Soient $\GG_1$ (resp. $\GG_2$) un polyèdre combinatoire qui possède une face triangulaire $T_1$ (resp. $T_2$) et une identification $\phi$ de $T_1$ et $T_2$ (qui renverse l'orientation). Le \emph{polyèdre combinatoire $\GG$ obtenu en recollant les polyèdres combinatoires $\GG_1$ et $\GG_2$ le long des faces triangulaires $T_1$ et $T_2$ via l'identification $\phi$} est le polyèdre combinatoire obtenu par le procédé suivant:
\begin{enumerate}
\item On identifie les arêtes et les sommets de $T_1$ (resp. $T_2$) via $\phi$, on obtient ainsi un triangle $\mathcal{T}$ inclus dans $\GG$.

\item On retire l'intérieur des arêtes du triangle $T$. Ainsi les 3 sommets $v_1,v_2,v_3$ de $T$ deviennent de valence 2.

\item Pour $i=1,..,3$, on note $e_i$ et $f_i$ les arêtes incidentes en $v_i$ et on fusionne ces deux arêtes en une seule arête en oubliant le sommet $v_i$.
\end{enumerate}
On obtient ainsi un graphe $\GG$ planaire, 3-connexe et de valence 3.
\end{defi}

\begin{figure}
\centerline{
\psfig{figure=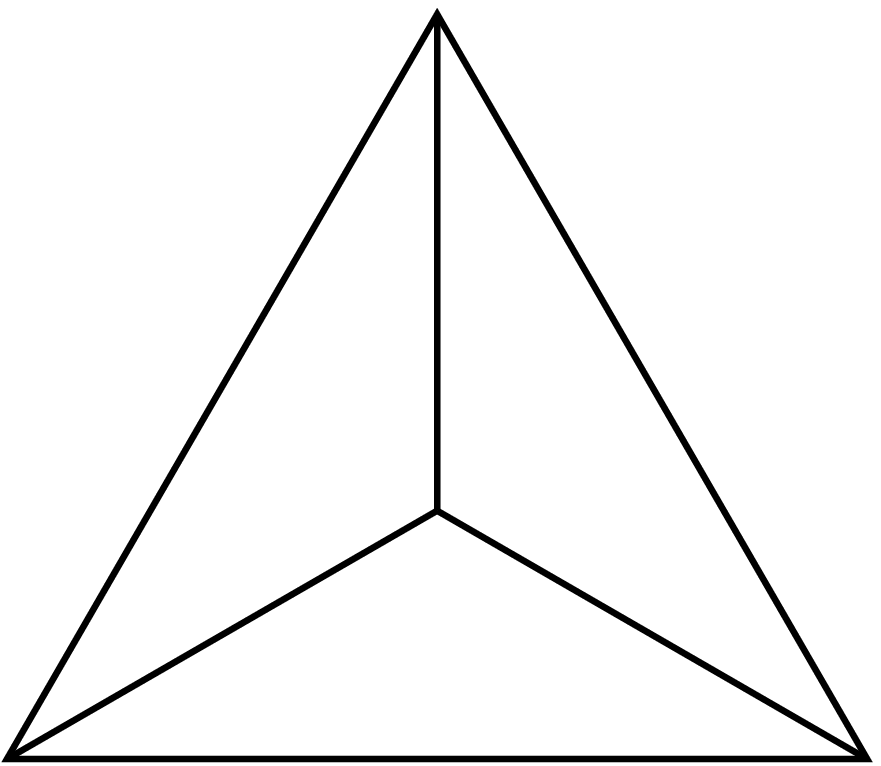, width=3cm}
\hspace{2em}
\psfig{figure=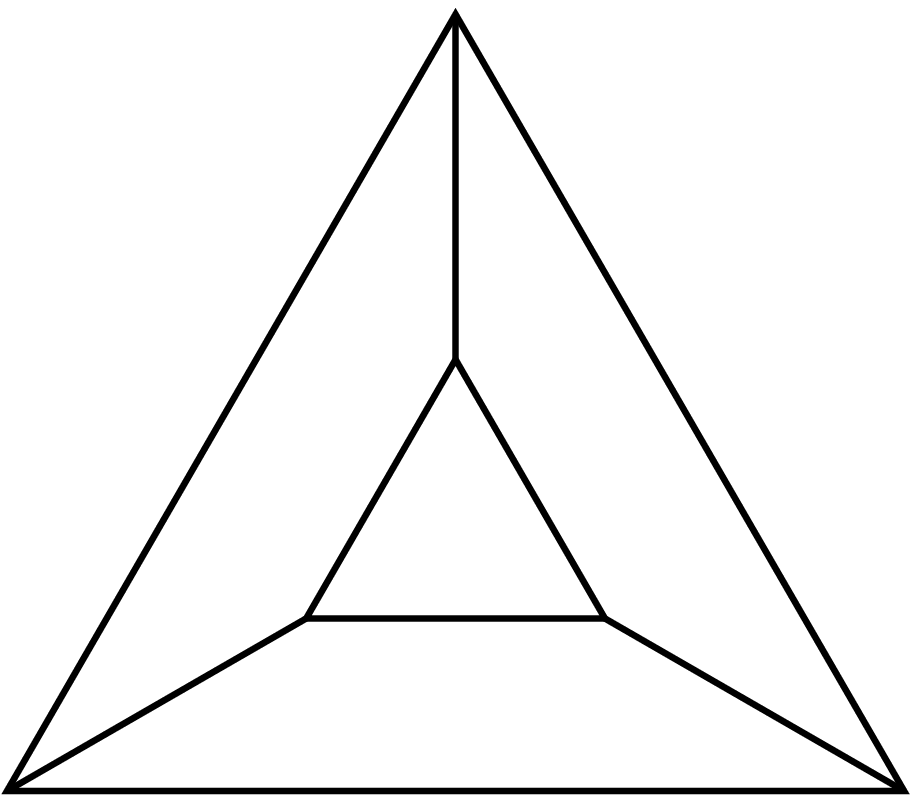, width=3cm}
\hspace{2em}
\psfig{figure=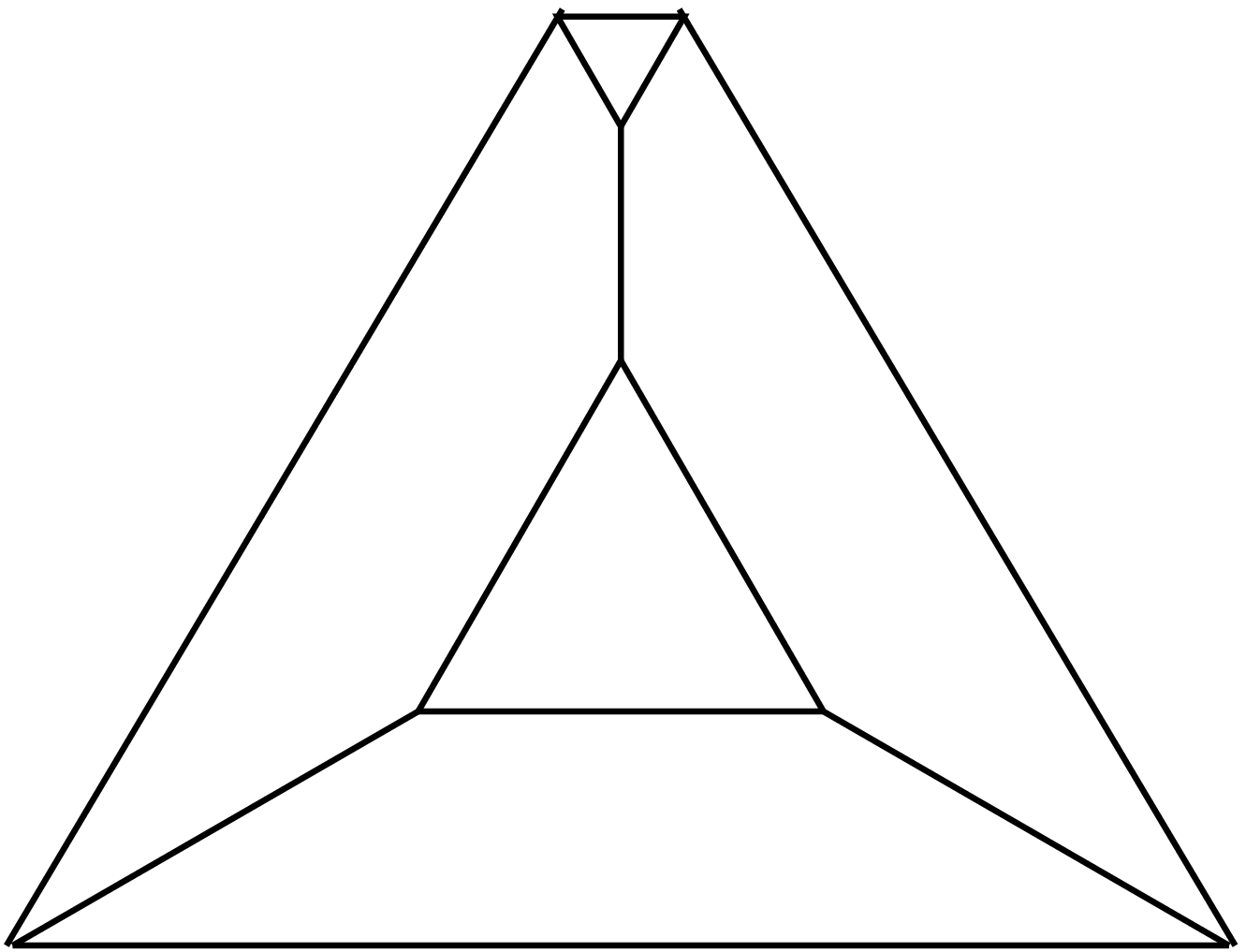, width=3cm}
}
\vspace{2em}
\centerline{
\psfig{figure=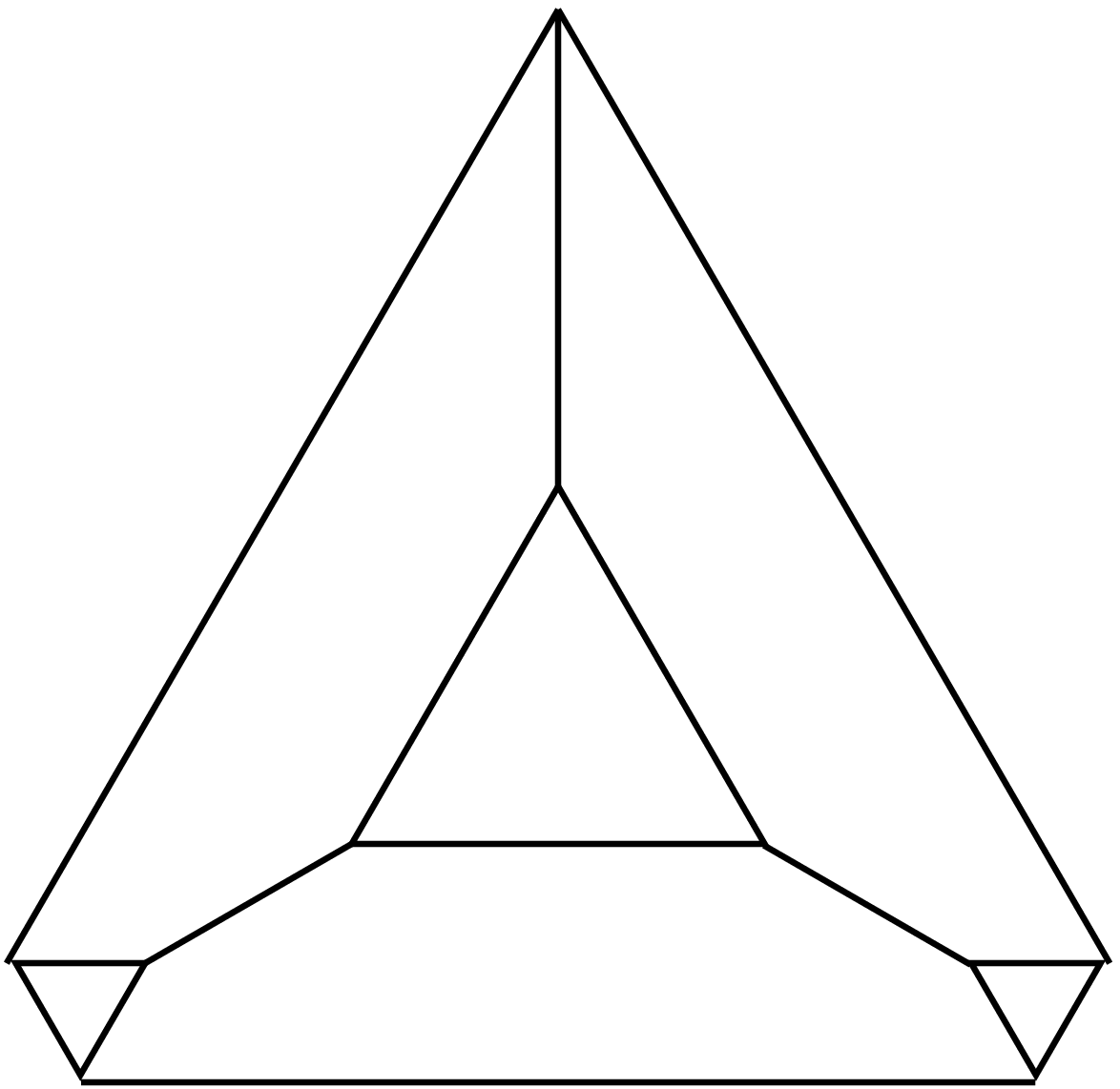, width=3cm}
\hspace{2em}
\psfig{figure=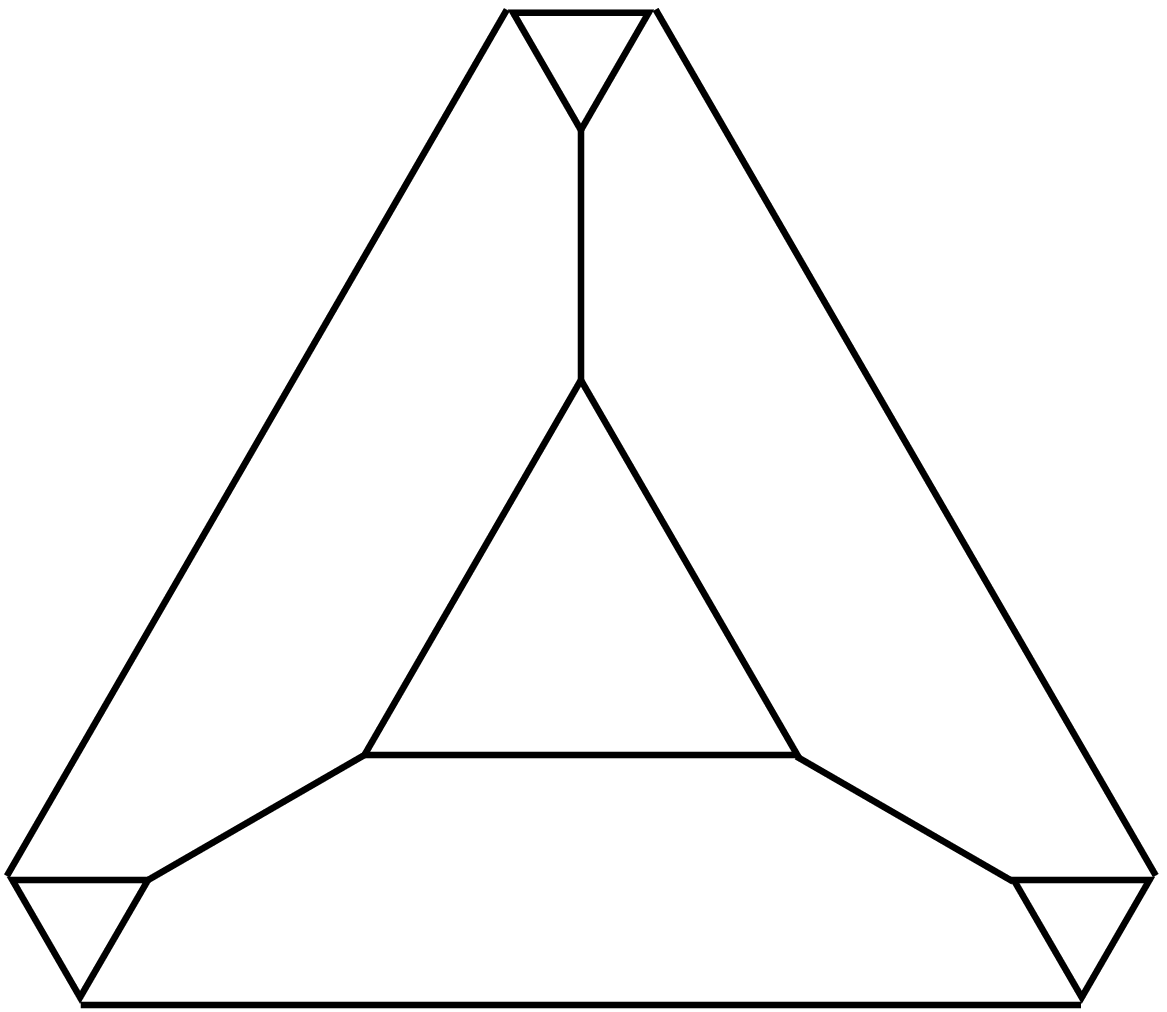, width=3cm}
}

\caption{Graphes de $\mathcal{T}_0$, $\mathcal{T}_1$, $\mathcal{T}_2$, $\mathcal{T}_3$ et $\mathcal{T}_4$} \label{tétra}
\end{figure}

Nous allons décomposer les écimaèdres combinatoires en \emph{blocs fondamentaux}. Il s'agit des écimaèdres $(\mathcal{T}_i)_{i=0,...,4}$. L'écimaèdre $\mathcal{T}_i$ est obtenu en écimant le tétraèdre combinatoire $\mathcal{T}_0$ en $i$ sommets distincts. Leurs graphes sont représentés par la figure \ref{tétra}.

\begin{prop}\label{constecima}
Tout écimaèdre combinatoire est obtenu en recollant un nombre fini de blocs fondamentaux le long de faces triangulaires. Inversement tout polyèdre combinatoire qui est obtenu comme recollement de blocs fondamentaux le long de faces triangulaires est un écimaèdre combinatoire.
\end{prop}

\begin{proof}
Les deux directions se montrent par récurrence.

Commençons par montrer que tout écimaèdre combinatoire $\GG$ est obtenu en recollant un nombre fini de blocs fondamentaux le long de faces triangulaires, par récurrence sur le nombre d'écimage subis par $\GG$.

L'initialisation de la récurrence est évidente. \`{A} présent, si $\GG$ a subi $n$ écimages alors il est obtenu en écimant un sommet $v$ d'un écimaèdre $\mathcal{H}$ qui est par hypothèse de récurrence obtenu en recollant des blocs fondamentaux $(\mathcal{B}_i)_{i \in I}$ en nombre fini. Le sommet $v$ appartient donc à un bloc fondamental $\mathcal{B}_{i_0}$. Distinguons deux cas, ou bien le sommet $v$ de $\mathcal{B}_{i_0}$ est un sommet du tétraèdre sous-jacent à $\mathcal{B}_{i_0}$ ou bien il est sur face provenant de l'écimage.

Si le sommet $v$ est un sommet du tétraèdre sous-jacent à $\mathcal{B}_{i_0}$ alors l'écimaèdre obtenu en écimant $\mathcal{B}_{i_0}$ au sommet $v$ est encore un bloc fondamental, par conséquent $\GG$ est obtenue en recollant un nombre fini de blocs fondamentaux le long de faces triangulaires.

Si le sommet $v$ est sur une face de $\mathcal{B}_{i_0}$ venant de l'écimage alors $v$ est sur une face triangulaire $s$ de $\mathcal{B}_{i_0}$. Dans ce cas, $\GG$ est obtenu en recollant l'une des deux faces triangulaires du bloc fondamental $\mathcal{T}_2$ avec la face $s$ de $\mathcal{H}$. Ainsi, $\GG$ est obtenu en recollant un nombre fini de blocs fondamentaux le long de faces triangulaires. Ce qui conclut la première partie de la démonstration.

Montrons à présent que tout polyèdre combinatoire qui est obtenu comme recollement de blocs fondamentaux le long de faces triangulaires est un écimaèdre, par récurrence sur le nombre de blocs fondamentaux du recollement.

L'initialisation de la récurrence est encore évidente. \`{A} présent, si le polyèdre combinatoire $\GG$ est obtenu en recollant $n$ blocs fondamentaux le long de faces triangulaires, alors par hypothèse de récurrence $\GG$ est obtenu en recollant un écimaèdre $\mathcal{H}$ et un bloc fondamental le long d'une face triangulaire $s$. Par conséquent, le polyèdre $\GG$ est obtenu en écimant la face $s$ en 1,2 ou 3 sommets. C'est donc un écimaèdre.
\end{proof}

On a d'abord besoin de quelques définitions purement combinatoire.

\begin{figure}[!h]
\begin{center}
\includegraphics[width=10cm]
{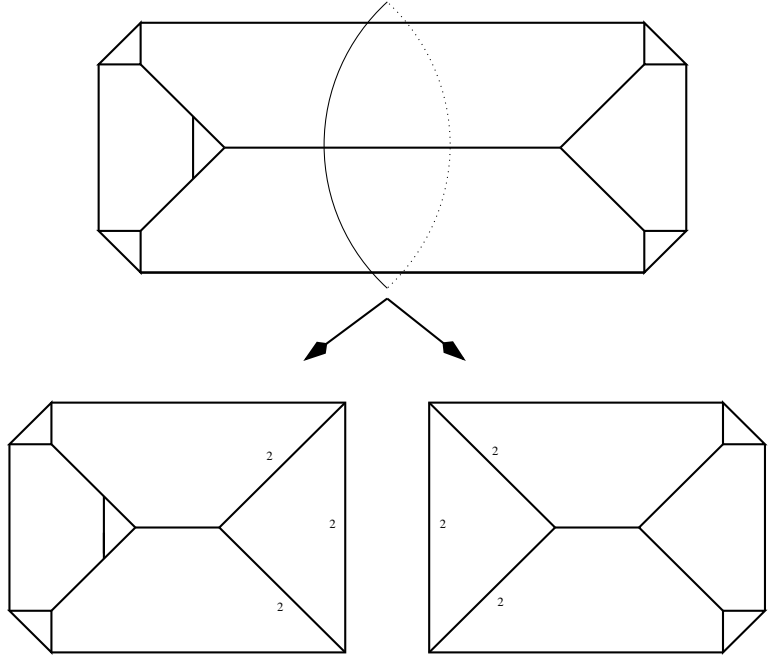}\caption{Coupe d'un graphe étiqueté} \label{coupe}
\end{center}
\end{figure}

\begin{defi}
Soient $\mathcal{G}$ un graphe planaire et 3-connexe et $\Gamma$
un 3-circuit orienté de $\GG$. La \emph{coupe droite (resp. gauche)
de $\GG$ le long de $\Gamma$} est le graphe $\GG_{\Gamma}^d$ (resp.
$\GG_{\Gamma}^g$) obtenu en retirant tous les sommets et les arêtes
de $\GG$ à gauche (resp. droite) de $\Gamma$ et en ajoutant un
triangle dont les sommets sont les extrémités des arêtes de
$\Gamma$ qui sont à gauche (resp. droite) de $\Gamma$. (Voir figure
\ref{coupe}.)
\end{defi}

\begin{defi}\label{deessens}
Soit $\GG$ un écimaèdre, on dit qu'un 3-circuit prismatique est \emph{combinatoirement essentiel} lorsqu'il ne
fait pas le tour d'une face triangulaire.
\end{defi}

\begin{prop}\label{3ciretri}
Si un écimaèdre $\GG$ est obtenu en recollant des blocs fondamentaux (qui ne sont pas des prismes triangulaires) le long de faces triangulaires, alors, les 3-circuits combinatoirement essentiels de $\GG$ sont exactement les triangles le long desquels se recollent les blocs fondamentaux.
\end{prop}

\begin{proof}
En effet, tous les triangles le long desquels se recollent les blocs fondamentaux correspondent naturellement à un 3-circuit prismatique combinatoirement essentiel. Enfin, il faut montrer par récurrence (de façon analogue à la proposition \ref{constecima}) que tout 3-circuit prismatique combinatoirement essentiel correspond à un triangle. On ne fait que la propriété d'hérédité, l'initialisation étant évidente. Soit $\G$ un 3-circuit prismatique combinatoirement essentiel d'un écimaèdre $\GG$. On suppose que $\GG$ est obtenu en recollant un écimaèdre $\mathcal{H}$ et un bloc fondamental $\mathcal{B}$ (qui n'est pas $\mathcal{T}_1$) le long d'une face triangulaire. Si $\G$ est inclus dans $\mathcal{H}$ ou $\mathcal{B}$ alors l'hypothèse de récurrence conclut cette remarque. Sinon $\G$ possède une face $f_1$ (resp. $f_3$) qui est une face de $\mathcal{H}$ (resp. $\mathcal{B}$) et qui n'est pas une face de $\mathcal{B}$ (resp. $\mathcal{H}$). Mais les faces $f_1$ et $f_3$ ne peuvent s'intersecter le long d'une arête, ce qui est absurde (le 3-circuit $\G$ possède une face $f_2$ qui est une face de $\mathcal{H}$ et $\mathcal{B}$, et qui est une face du triangle de recollement de $\mathcal{H}$ et $\mathcal{B}$).
\end{proof}

\begin{rema}\label{3circourbe}
Comme tout écimaèdre est obtenu en recollant des blocs fondamentaux le long de faces triangulaires, il est facile de voir que tout 3-circuit non prismatique de $\GG$ est inclus dans un et un seul des blocs fondamentaux qui décomposent $\GG$. La proposition \ref{3ciretri} montre que tout 3-circuit prismatique non combinatoirement essentiel de $\GG$ est inclus dans un ou deux blocs fondamentaux qui décomposent $\GG$, et s'il est inclus dans 2 alors l'un d'eux est un prisme triangulaire. Enfin, la proposition \ref{3ciretri} montre que tout 3-circuit prismatique combinatoirement essentiel de $\GG$ est inclus dans exactement deux blocs fondamentaux. Par conséquent, l'ensemble des 3-circuits de $\GG$ peut être réalisé comme un ensemble de courbes simples, continues, tranverses aux arêtes de $\GG$ et disjointes. Nous allons découper $\GG$ le long de ces courbes.
\end{rema}

On a des définitions analogues pour les graphes étiquetés.

\begin{defi}
Soient $\mathcal{G}$ un graphe étiqueté et $\Gamma$ un 3-circuit
orienté de $\GG$; on définit la \emph{coupe droite (resp. gauche) de $\GG$ le
long de $\Gamma$}, comme le graphe étiqueté $\GG_{\Gamma}^d$ (resp.
$\GG_{\Gamma}^g$) dont le graphe sous-jacent est la coupe droite
(resp. gauche) du graphe sous-jacent à $\GG$, et on étiquette les
nouvelles arêtes de $\GG_{\Gamma}^d$ (resp. $\GG_{\Gamma}^g$) par des
2 et les anciennes gardent leurs étiquettes. (Voir figure
\ref{coupe}.)
\end{defi}

\begin{rema}
Lorsqu'il n'y a pas d'ambiguïté sur le 3-circuit $\Gamma$, on
allégera les notations $\GG_{\Gamma}^d$ et $\GG_{\Gamma}^g$ en
notant $\GG_{\Gamma}^d = \GG^d$ et $\GG_{\Gamma}^g=\GG^g$.
\end{rema}

\begin{defi}\label{combiessens}
Soit $\GG$ un écimaèdre étiqueté, on dit qu'un 3-circuit prismatique est \emph{essentiel} lorsqu'il ne
fait pas le tour d'une face triangulaire dont toutes les arêtes sont d'ordre 2.
\end{defi}

\begin{defi}
Un $\emph{bloc fondamental étiqueté}$ est un graphe étiqueté et dont le
graphe sous-jacent est un bloc fondamental combinatoire, et les
faces triangulaires issues de l'écimage ne possèdent que des arêtes
d'ordre 2.
\end{defi}

\begin{rema}
Il faut bien faire attention aux blocs fondamentaux qui sont des prismes triangulaires. En effet, ces derniers possèdent deux faces triangulaires mais une seule de ces deux faces provient de l'écimage. Par conséquent, les blocs fondamentaux qui sont des prismes triangulaires ne possèdent a priori qu'une seule face triangulaire qui porte uniquement des arêtes d'ordre 2.
\end{rema}

Tout cela permet d'associer à tout écimaèdre $\GG$, un arbre $\A_{\GG}$ qui \og mémorise \fg $\,$ les écimages effectués. Cet arbre permet de coder une partie de la combinatoire de $\GG$.

Nous allons découper $\GG$ le long de ces 3-circuits prismatiques. Pour cela, on se donne une famille $(\G_j)_{j \in J}$ de 3-circuits orientés de $\GG$ où chaque 3-circuit est représenté exactement une fois. La proposition \ref{3ciretri} et la remarque \ref{3circourbe} montrent que si l'on découpe $\GG$ le long d'un 3-circuit prismatique $\G_{i_0}$ de $\GG$ alors les 3-circuits $(\G_i)_{i \neq i_0}$ sont des 3-circuits de $\GG_{\Gamma}^d$ ou bien de $\GG_{\Gamma}^g$. On découpe $\GG$ le long de tous les 3-circuits $(\G_j)_{j \in J}$ qui sont prismatiques. Les écimaèdres obtenus par ce découpage ne possèdent plus aucun 3-circuit prismatique essentiel. Par conséquent, la proposition \ref{3ciretri} et la remarque \ref{3circourbe} montrent que ce sont des blocs fondamentaux étiquetés.

\begin{defi}
Soit $\GG$ un écimaèdre combinatoire, notons $(\G_j)_{j \in J}$ la famille complète des 3-circuits orientés de $\GG$. \emph{L'arbre associé à $\GG$} sera noté $\A_{\GG}$ et il est défini de la façon suivante: les sommets de $\A_{\GG}$ sont d'une part les sommets de $\GG$, d'autre part les blocs fondamentaux $(\mathcal{B}_i)_{i=1...n+1}$ résultant du découpage de $\GG$ le long des 3-circuits prismatiques de la famille $(\G_j)_{j \in J}$. Deux sommets de $\A_{\GG}$ qui sont des blocs fondamentaux $\mathcal{B}_i$ et $\mathcal{B}_k$ de $\A_{\GG}$ sont reliés par une arête lorsque $\mathcal{B}_i$ et $\mathcal{B}_k$ partagent un 3-circuit de $(\G_j)_{j \in J}$. Un sommet $v$ de $\GG$ est relié à un bloc fondamental lorsque $v$ appartient à ce bloc.
\end{defi}

On représentera les arêtes extrémales de $\A_{\GG}$ en pointillés, car elles sont en bijection avec les sommets de $\GG$, alors que les arêtes non extrémales sont en bijection avec les 3-circuits
prismatiques de $\GG$. De plus, les arêtes extrémales du sous-arbre en trait plein sont en bijection avec les 3-circuits prismatiques non combinatoirement essentiels.

\begin{figure}
\centerline{
\psfig{figure=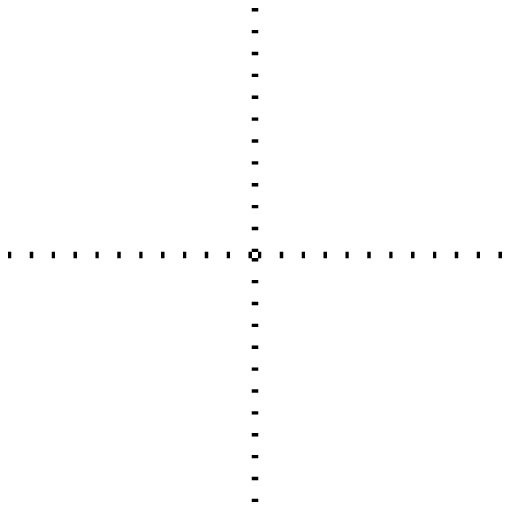, width=3cm}
\hspace{3em}
\psfig{figure=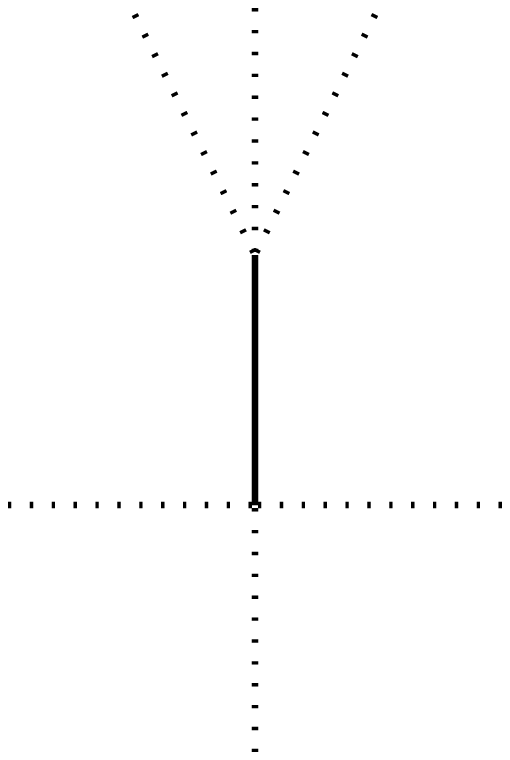, width=2.3cm}
\hspace{3em}
\psfig{figure=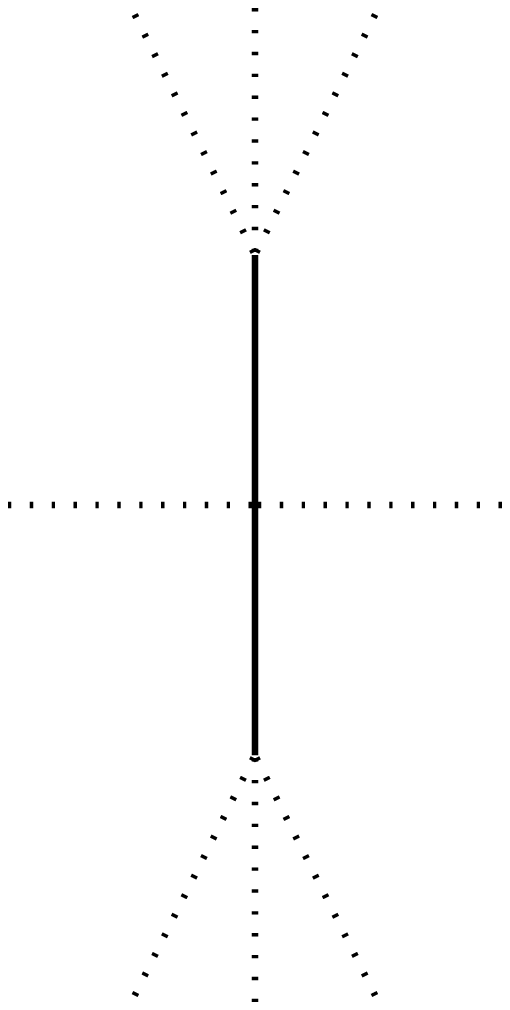, width=1.8cm}
}
\vspace{2em}
\centerline{
\psfig{figure=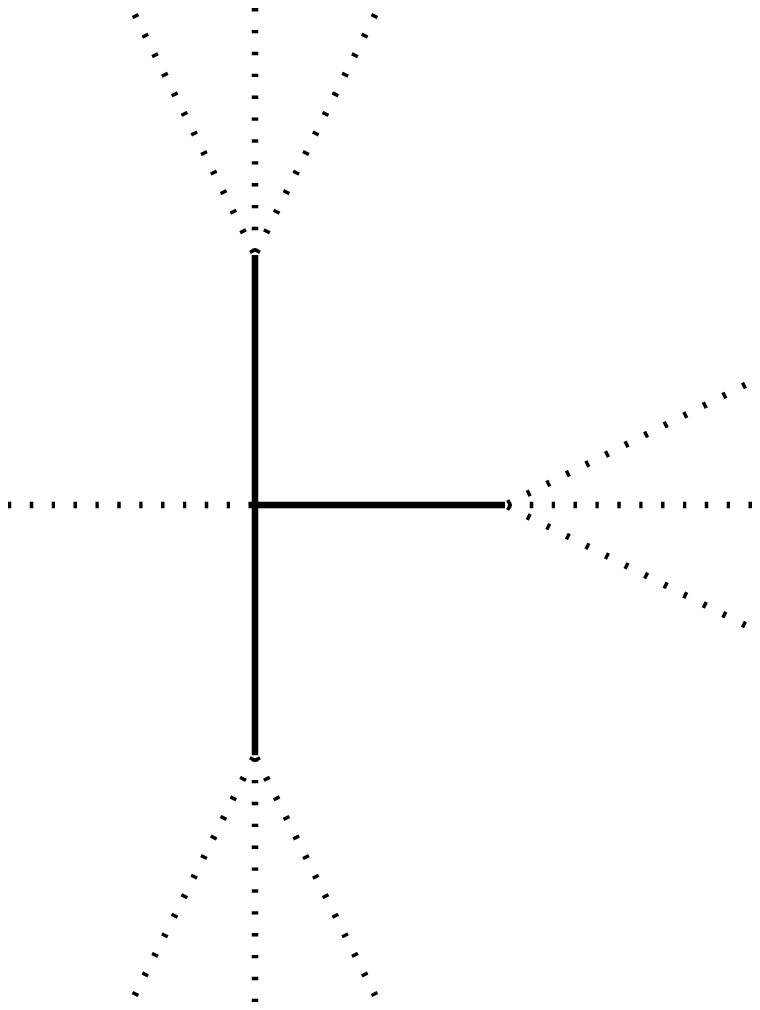, width=2.6cm}
\hspace{3em}
\psfig{figure=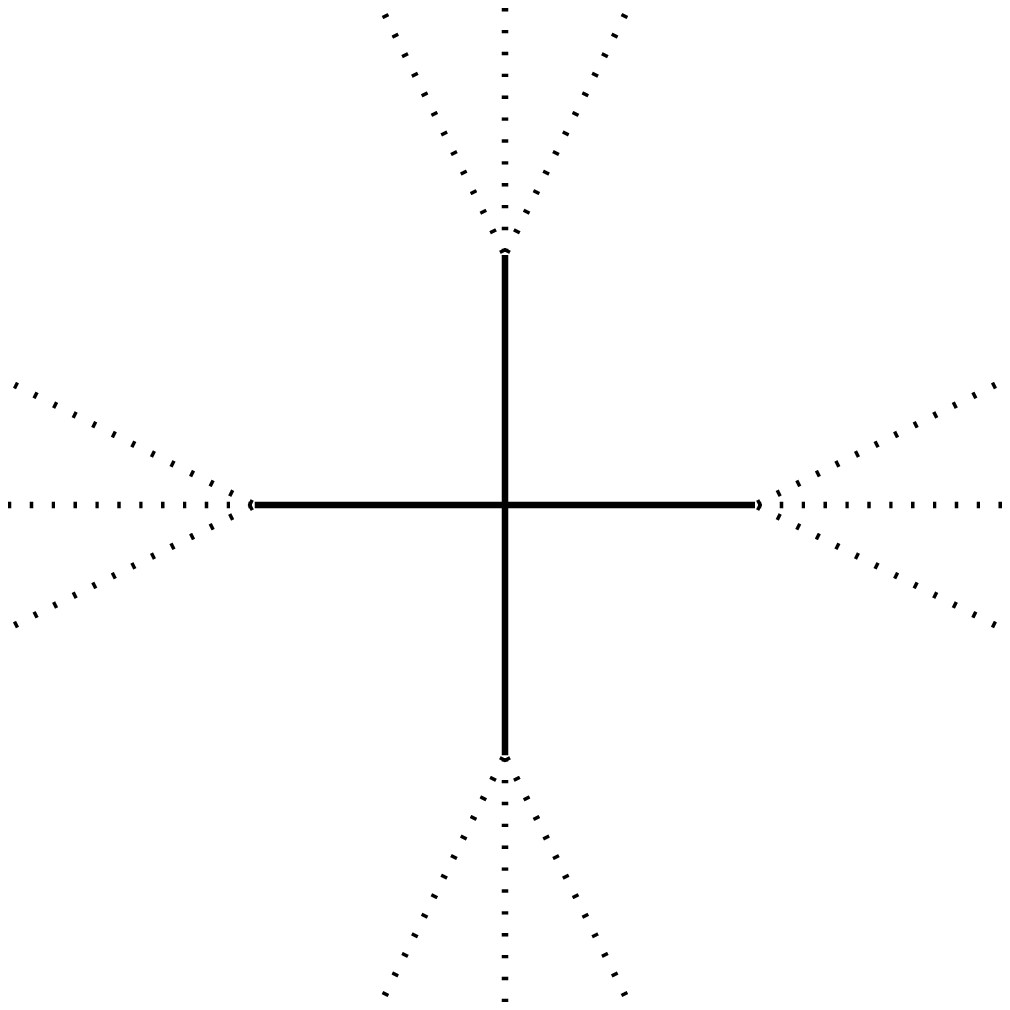, width=3.5cm}
}
\caption{Les arbres $\A_{\mathcal{T}_0}$, $\A_{\mathcal{T}_1}$, $\A_{\mathcal{T}_2}$, $\A_{\mathcal{T}_3}$ et $\A_{\mathcal{T}_4}$.}

\label{tétra1c}
\end{figure}

On a représenté les arbres $\A_{\mathcal{T}_i}$ pour $i=1,...,4$ à l'aide de la figure \ref{tétra1c}.

\begin{rema}
L'application $\GG \mapsto A_{\GG}$ n'est pas injective. En revanche, on vérifie facilement que l'image de cette application est l'ensemble des arbres dont les sommets non extrémaux sont de valence 4.
\end{rema}

Nous allons montrer que l'on peut construire tous les écimaèdres
miroirs dont le graphe sous-jacent n'est pas celui d'un tétraèdre en recollant \og des blocs
fondamentaux miroirs \fg.

\subsection{Enoncé du résultat}

Pour énoncer le résultat, on a besoin d'introduire
plusieurs quantités associées à $\GG$. Nous allons calculer
la dimension attendue de $X_{\GG}$ de façon heuristique.

Rappelons que lorsque deux faces adjacentes $s$ et $t$ partagent une arête $s \cap t$ qui porte l'étiquette $\mu_{st}$ alors (définition \ref{miroir} et notations \ref{equa}) :
\begin{itemize}
\item On a $v_{ts} = -\alpha_s(v_t) \geqslant 0$.

\item Si l'arête $s \cap t$ est d'ordre 2 alors les deux équations suivantes doivent être vérifiés : $v_{ts}=v_{st}=0$.

\item Si l'arête $s \cap t$ n'est pas d'ordre 2 alors une seule équation doit être vérifiée : $v_{st}v_{ts}= \mu_{st}$.
\end{itemize}

On peut voir $X_{\GG}$ de la façon suivante: soient $S$ l'ensemble
des faces de $\GG$ et $f$ le cardinal de $S$.
\begin{enumerate}

\item Nombre d'inconnues: $([\alpha_s],[v_s])_{s \in S} \in (\mathbb{P}^+((\R^4)^*) \times
\mathbb{P}^+(\R^4))^S$ est une variété de dimension $6f$.

\item Nombre d'équations: on a $e+e_2$ équations, où $e$ est le
nombre d'arêtes de $\GG$ et $e_2$ le nombre d'arêtes d'ordre 2.

\item Modulo $\sss^{\pm}_4(\R)$ qui agit librement et proprement sur
$Y_{\GG}$ et $\dim(\sss^{\pm}_4(\R))=15$.

\end{enumerate}

Donc la dimension attendue de $\GG$ est $d(\GG) = 6f-e-e_2-15$. Or,
$\GG$ est de valence 3, il possède donc $\frac{2}{3}e$ sommets. La relation d'Euler montre que $f = \frac{6+e}{3}$, il vient donc $d(\GG)=(e-e_2)-3=e_+-3$, où $e_+=e-e_2$ est le nombre
d'arêtes de $\GG$ d'ordre différent de 2. Il est intéressant de noter que ce calcul heuristique ne suppose pas que $\GG$ est un écimaèdre.

\begin{defi}
Soit $\GG$ un écimaèdre étiqueté, on définit les quantités suivantes:
\begin{itemize}
\item $d(\GG)=e_+-3$.

\item $n(\GG)$ est le nombre de 3-circuits affines ou sphériques, prismatiques, sans
angle droit de $\GG$.

\item $m(\GG)$ est le nombre de 3-circuits affines ou sphériques, prismatiques, avec
angle droit de $\GG$.
\end{itemize}
\end{defi}

\begin{theo}\label{theo}
Soit $\GG$ un écimaèdre étiqueté marqué. Alors,
\[
\begin{tabular}{llll}

1) & Si $m(\GG)=0$ et $d(\GG) \geqslant 0$   & alors & $X_{\GG}$ est difféomorphe à $\kappa(\GG)$ copies de $\R^{d(\GG)}$,\\
   & et $\GG$ n'est pas un prisme            &       &   où $\kappa(\GG)$ est un entier pair ou égal à 1,\\
   & exceptionnel                            &       & qui vérifie $1 \leqslant \kappa(\GG) \leqslant 2^{n(\GG)}$.\\
\\
2) & Si $m(\GG)=0$ et $d(\GG) < 0$         & alors & on a deux cas:\\

   &                                     &       & a) Le graphe
sous-jacent à $\GG$ est un tétraèdre\\
   &                                     &       & combinatoire et $X_{\GG}$ est
un singleton.\\

   &                                     &       & b) $\GG =
\GG_{\alpha,\beta,\frac{\pi}{2}}$ pour un certain couple $(\alpha,\beta)$\\

   &                                     &       & avec
$\alpha+\beta < \frac{\pi}{2}$ et $X_{\GG}= \varnothing$.\\
\\
3) & Si $m(\GG) > 0$                    & alors  & on a deux cas \\
   &                                    &        & a) si $\GG =\GG_{\alpha,\frac{\pi}{2}-\alpha,\frac{\pi}{2}}$ alors $X_{\GG}$ est un singleton\\
   &                                    &        & et $m(\GG)=1$.\\
   &                                    &        & b) sinon (en particulier si $\GG = \GG_{\alpha,\beta,\frac{\pi}{2}}$ \\
   &                                    &        & avec $\alpha + \beta > \frac{\pi}{2}$) alors $X_{\GG} = \varnothing$.\\
\\
4) & Si $\GG = \GG_{\alpha,\beta,\gamma}$ est un prisme & alors & on a 3 cas\\
   &                                                    & & a) si $\alpha+\beta+\gamma < \pi$ alors $X_{\GG}=
\varnothing$.\\
   &  exceptionnel avec  & & b) si $\alpha+\beta+\gamma = \pi$ alors $X_{\GG}$ est un
   singleton.\\
   &      $0<\alpha,\beta,\gamma < \frac{\pi}{2}$             & & c) si $\alpha+\beta+\gamma > \pi$ alors $Card (X_{\GG}) =2$.
\end{tabular}
\]
\end{theo}

\begin{rema}
Il faut bien remarquer que les points 1, 2, 3, et 4 sont disjoints. La seule difficulté est de voir que le cas 4) est disjoint des trois autres. Mais dans le cas 4, on a $m(\GG)=0$, $d(\GG)=0$ et $\GG$ est un prisme exceptionnel.
\end{rema}

\begin{coro}
Soit $\GG$ un écimaèdre étiqueté marqué: alors $X_{\GG}$ est une variété qui possède un nombre fini de composantes connexes, toutes les composantes connexes de $X_{\GG}$ sont homéomorphes à des boules et elles ont toutes la même dimension.

Enfin, l'espace $X_{\GG}$ est vide si et seulement si on est dans l'un des cas suivants:
\begin{enumerate}
\item $X_{\GG}$ est un prisme exceptionnel dont l'unique 3-circuit prismatique est hyperbolique.


\item $m(\GG) \geqslant 1$ et $\GG$ n'est pas un prisme exceptionnel dont l'unique 3-circuit prismatique est affine avec angle droit.

\end{enumerate}
\end{coro}

\begin{proof}
En effet, si $m(\GG)=0$ alors les cas 1), 2) et 4) du théorème \ref{theo} montrent que $X_{\GG} = \varnothing$ si et seulement si $\GG$ est un prisme exceptionnel dont l'unique 3-circuit prismatique est hyperbolique. Si  $m(\GG) \geqslant 1$ alors c'est le cas 3) du théorème \ref{theo} qui conclut.
\end{proof}

\begin{coro}\label{corohyp}
Si $\GG$ est un écimaèdre étiqueté qui vérifie les conditions du
théorème d'Andreev alors $X_{\GG}$ est difféomorphe à $\R^{d(\GG)}$.
\end{coro}

\begin{proof}
\par{
Il est clair que $m(\GG)=n(\GG)=0$ et que $\GG$ n'est pas un prisme exceptionnel, il reste donc à vérifier que $d(\GG) \geqslant 0$, pour appliquer le point 1). Il faut noter que $\GG$ n'est pas un tétraèdre puisqu'il vérifie les hypothèses du théorème d'Andreev. Le plus simple est de distinguer deux cas: $\GG$ est un prisme (dont le circuit prismatique est nécessairement hyperbolique) ou bien $\GG$ n'est ni un tétraèdre, ni un prisme.
}
\\
\par{
Dans le premier cas, il est facile de voir que tout prisme non exceptionnel qui vérifie que $m(\GG)=0$ possède au moins 3 arêtes qui ne sont pas d'ordre 2. Par suite, $d(\GG) \geqslant 0$.
}
\\
\par{
Si $\GG$ est un écimaèdre qui n'est ni un tétraèdre, ni un prisme alors $\GG$ possède au moins deux 3-circuits prismatiques. Comme $\GG$ vérifie les hypothèses du théorème d'Andreev, tout 3-circuit prismatique de $\GG$ est hyperbolique et possède donc au plus une seule arête d'ordre 2. L'écimaèdre étiqueté $\GG$ possède donc au moins 3 arêtes d'ordre différent de 2 (la 3-connexité de $\GG$ entraîne que deux 3-circuits distincts de $\GG$ ont au plus une arête en commun). Par suite, $d(\GG) \geqslant 0$.
}
\end{proof}

\begin{rema}
L'hypothèse d'ordonnabilité de Choi est plus forte dans le cas où le graphe $\GG$ est trivalent que l'hypothèse de cet article selon laquelle $\GG$ est un écimaèdre. Nous avons supposé que les graphes étiquetés étaient des graphes de valence 3, et dans ce cas, on peut montrer que, si $\GG$ est un graphe étiqueté ordonnable, alors le graphe sous-jacent à $\GG$ est un écimaèdre combinatoire. Mais il existe
des graphes planaires et 3-connexes de valence supérieure à 3 qui sont ordonnables (s'ils sont étiquetés correctement). Il faut noter que pour construire un convexe divisible via la méthode de
Vinberg, une condition nécessaire (voir la remarque \ref{valence3}) est que le polyèdre de départ soit de valence 3.
\end{rema}

\begin{rema}
Pour calculer $\kappa(\GG)$ nous allons retirer certaines arêtes de
l'arbre $\A_{\GG}$ et on obtient une forêt que l'on notera
$\F_{\GG}$. Puis on introduira la notion d'orientation partielle et
d'orientation partielle admissible de la forêt $\F_{\GG}$ (voir
partie \ref{5}). Enfin on démontrera que $\kappa(\GG)$ est le
nombre d'orientations partielles admissibles de la forêt $\F_{\GG}$
(voir le théorème \ref{final}).
\end{rema}

\begin{rema}
Nous énoncerons à la fin de ce texte le théorème \ref{final} qui
précise la partie 1) du théorème \ref{theo} en donnant une
paramétrisation explicite de $X_{\GG}$.
\end{rema}

Nous allons commencer par les points 2) et 3) qui sont beaucoup plus
faciles à montrer. Puis on attaquera la démonstration
des points 1) et 4).

\subsection{Démonstration des points 2 et 3}\label{debutdemo}

\begin{lemm}\label{lemme1}
Soit $P$ un polyèdre miroir, soient $s$ et $t$ deux faces distinctes
quelconques de $P$; alors $v_{ts} = -\alpha_s(v_t) \geqslant 0$ et
$v_{ts} = -\alpha_s(v_t) =0$ si et seulement si $s$ et $t$ partagent une
arête d'ordre 2.
\end{lemm}

\begin{proof}
Si $s$ et $t$ sont adjacentes alors les conditions 1)-2\textquoteright $b$) (définition \ref{miroir} et notations \ref{equa}) entraînent le lemme. \`{A} présent, soient $s$ et $t$ deux faces non adjacentes de $P$, et $r_1,...,r_k$ les faces de $P$ adjacentes à $t$. Comme $P$ est un
polyèdre proprement convexe, son dual dans $\mathbb{P}^+(\R^4)^*$ est aussi un polyèdre proprement convexe dont les sommets sont les faces de $P$. Par conséquent, le segment qui relie $[-\alpha_s]$ à $[-\alpha_t]$ traverse l'enveloppe convexe des $([-\alpha_{r_i}])_{i=1,...,k}$ dans $\mathbb{P}^+(\R^4)^*$. Par conséquent, $-\alpha_s$ peut s'écrire comme une combinaison linéaire à coefficients positifs des $(-\alpha_{r_i})_{i=1,...,k}$ et de $\alpha_t$ dont le coefficient associé à $\alpha_t$ est strictement positif. Les inégalités $-\alpha_{r_i}(v_t) \geqslant 0$ pour $i=1...k$ et l'égalité $\alpha_t(v_t)=2$ entraînent l'inégalité $-\alpha_s(v_t) > 0$.
\end{proof}

\begin{rema}
Rappelons que pour tout polyèdre projectif miroir $P$, pour toute face $s$ de $P$, la polaire $[v_s]$ de $s$ n'appartient pas à $P$ puisque $-\alpha_s(v_s) < 0$. Le lemme \ref{lemme1} montre que tout segment projectif issu de $[v_s]$ et rencontrant $P$, rencontre en premier la face $s$.
\end{rema}

Le lemme suivant démontre le point 3) et un résultat sur les
3-circuits affines ou sphériques, prismatiques  qui  sera utile
à plusieurs reprises.

\begin{lemm}\label{tropde2}
Soit $\GG$ un graphe étiqueté tel que $X_{\GG} \neq \varnothing$, alors
tout 3-circuit $\Gamma$ affine ou sphérique, prismatique de $\GG$
possède au plus un angle droit. Et s'il en possède un alors  $\GG$ est le prisme exceptionnel $\GG_{\alpha,\frac{\pi}{2}-\alpha,\frac{\pi}{2}}$ pour un certain $\alpha \in ]0,\frac{\pi}{2}[$ et $X_{\GG}$ est un
singleton.
\end{lemm}

\begin{proof}
Soient $P$ un polyèdre miroir qui réalise $\GG$ et $\Gamma$ un
3-circuit orienté prismatique avec angle droit, affine ou
sphérique.

Commençons par numéroter de 1 à 3 les faces qui forment $\Gamma$
de telle sorte que l'arête adjacente aux faces 1 et 2 soit d'ordre
2. Le 3-circuit $\Gamma$ est prismatique donc il existe une face à gauche de
$\Gamma$ et une à droite de $\Gamma$ qu'on numérote 4 et 5. On
peut supposer que l'on a la configuration suivante dans une base
$(e_j)_{j=1...4}$ de $\R^4$: $\alpha_i = -e_i^*$ pour $i=1...4$.
On souhaite avoir une forme simple pour la forme linéaire $\alpha_5$, on va montrer que l'on peut supposer que $\alpha_5 = -e_1^* - e_2^* - e_3^* + e_4^*$.
Le stabilisateur des $([\alpha_i])_{i=1...4}$ est le groupe des matrices diagonales à diagonales strictement positives dans la base $(e_j)_{j=1...4}$ de $\R^4$. La forme linéaire $\alpha_5$ se décompose dans la base duale de la base $(e_j)_{j=1...4}$. Il suffit de vérifier que les signes des coefficients de cette décomposition sont les bons. Le polyèdre dual de $P$ est un polyèdre convexe dont les sommets sont les faces de $P$, par conséquent le segment reliant $[-\alpha_4]$ à $[-\alpha_5]$ dans $\mathbb{P}^+(\R^4)$ traverse le triangle dont les sommets sont les points $[-\alpha_1]$, $[-\alpha_2]$ et $[-\alpha_3]$. Par conséquent, $\alpha_5$ est une combinaison linéaire à coefficients positifs de $-\alpha_1$, $-\alpha_2$, $-\alpha_3$ et $\alpha_4$. On peut donc supposer que $\alpha_5 = -e_1^* - e_2^* - e_3^* + e_4^*$.
\newline
Soient $v_{ij}$ les coordonnées de $(v_i)_{i=1...5}$. On a:

$$
\left(
\begin{array}{c}
v_1 \\
v_2 \\
v_3
\end{array}
\right)
=
\left(
\begin{array}{cccc}
-2 & 0 & v_{13} & v_{14}\\
 0 & -2 & v_{23} & v_{24}\\
v_{31}& v_{32} & -2 & v_{34}
\end{array}
\right)
$$

On rappelle que si les faces $i$ et $j$ sont adjacentes alors $v_{ij}v_{ji} = \mu_{ij}=  4\cos^2(\theta_{ij})$ (définition \ref{miroir} et notations \ref{equa}). Le lemme \ref{lemme1} affirme que $v_{ij} \geqslant 0$ pour $i \neq j$ et
que les inégalités suivantes doivent être vérifiées:

$$
\left\{
\begin{array}{lll}
-\alpha_5(v_1) \geqslant 0 & \textrm{c'est-à-dire} & v_{13} \geqslant 2 + v_{14},\\
-\alpha_5(v_2) \geqslant 0 & \textrm{c'est-à-dire} & v_{23} \geqslant 2 + v_{24},\\
-\alpha_5(v_3) \geqslant 0 & \textrm{c'est-à-dire} &
\frac{\mu_{13}}{v_{13}}+\frac{\mu_{23}}{v_{23}} \geqslant 2 +
v_{34}.
\end{array} \right.
$$

Donc $v_{13} \geqslant 2$ et $v_{23}\geqslant 2$ et ainsi,
$\mu_{13}+\mu_{23} \geqslant 4$ et donc $\theta_{13}+\theta_{23}
\leqslant \frac{\pi}{2}$. Or, $\theta_{13}+\theta_{23} \geqslant
\frac{\pi}{2}$ (car $\Gamma$ est affine ou sphérique) donc
$\theta_{13}+\theta_{23} = \frac{\pi}{2}$, et $\mu_{13}+\mu_{23} =
4$. Il vient finalement que $v_{13}=v_{23}=2$, et donc $v_{31} =\frac{\mu_{13}}{2}$ et $v_{32} =\frac{\mu_{23}}{2}$, puis
$v_{14}=v_{24}=v_{34}=0$, et enfin $\alpha_5(v_1)= \alpha_5(v_2) =
\alpha_5(v_3) = 0$. Les faces 4 et 5 sont donc adjacentes aux
faces 1, 2 et 3, d'après le lemme \ref{lemme1}.

Ainsi, le polyèdre miroir $Q$ obtenu à l'aide des faces 1,2,3,4,5
de $P$ réalise un $\GG_{\alpha,\frac{\pi}{2}-\alpha,\frac{\pi}{2}}$ pour $\alpha=\arccos(\sqrt{\frac{\mu_{13}}{4}}) \in ]0,\frac{\pi}{2}[$ et le polyèdre $Q$ est clairement
unique puisque les coordonnées des vecteurs $(v_i)_{i=1...5}$ ont
été calculées explicitement.

Cela montre aussi que si $X_{\GG}$ est non vide alors tout 3-circuit prismatique affine ou sphérique avec angle droit est affine. Par conséquent, le polyèdre $Q$ est le polyèdre $P$ lui-même. En effet, si $P$ possédait d'autres faces, alors comme $P$ est un écimaèdre, $P$ posséderait un 3-circuit prismatique sphérique avec angle droit provenant d'un des six 3-circuits non prismatiques sphériques de $Q$.
\end{proof}

Maintenant que le point 3) est acquis, on peut montrer le point 2).

\begin{proof}[Démonstration du point 2]
On a par hypothèse $m(\GG) = 0$ (i.e. tout 3-circuit prismatique avec angle droit de $\GG$ est hyperbolique). Par conséquent, aucun 3-circuit prismatique de $\GG$ ne possède deux angles droits. Mais on a aussi supposé que $d(\GG) < 0$, par conséquent $e_+ \leqslant 2$. Comme dans un graphe 3-connexe deux 3-circuits distincts ont au plus une arête en commun, l'écimaèdre étiqueté $\GG$ contient au plus un 3-circuit prismatique.

Les seuls écimaèdres qui possèdent au plus un 3-circuit
prismatique sont le tétraèdre combinatoire et le prisme
triangulaire combinatoire.

Commençons par le cas où le graphe sous-jacent à $\GG$ est un prisme triangulaire combinatoire. On a supposé que $d(\GG) < 0$ donc $\GG$ possède au moins 7 arêtes d'ordre 2, et comme $\GG$ n'a pas de 3-circuit affine ou sphérique, prismatique, avec angle droit ($m(\GG) = 0$), $\GG$ possède exactement 7 arêtes d'ordre 2, et $\GG = \GG_{\alpha,\beta,\frac{\pi}{2}}$ pour un certain couple $(\alpha,\beta) \in ]0,\frac{\pi}{2}[^2$ tels que $\alpha +\beta < \frac{\pi}{2}$. On reprend le calcul de la démonstration du lemme \ref{tropde2}.

On rappelle que dans la démonstration  du lemme \ref{tropde2}, on a en particulier calculé l'espace des modules d'un prisme exceptionnel dont l'unique 3-circuit prismatique est à angle droit et affine ou sphérique. \`{A} présent, on doit calculer l'espace des modules d'un prisme exceptionnel dont l'unique 3-circuit prismatique est à angle droit et hyperbolique.

On reprend les mêmes notations pour les faces et on suppose toujours que c'est l'arête partagée par les faces 1 et 2 qui est d'ordre 2. On obtient les coordonnées suivantes pour les formes linéaires $(\alpha_i)_{i=1,...,5}$ et les vecteurs $(v_i)_{i=1,...,5}$.
$$
\begin{array}{cccc}
\alpha_i = -e^*_i, \, i=1,...,4\\
\alpha_5 = -e_1^*-e_2^*-e_3^*+e_4^*\\
\end{array}
$$

$$
\left(
\begin{array}{c}
v_1 \\
v_2 \\
v_3
\end{array}
\right)
=
\left(
\begin{array}{cccc}
-2      & 0      & v_{13} & 0\\
0       & -2      & v_{23} & 0\\
v_{31} & v_{32}  & -2     & 0
\end{array}
\right)
$$

La face numérotée 5 de $\GG$ ne contient que des arêtes d'ordre 2 par conséquent on a:

$$
\left\{
\begin{array}{lll}
-\alpha_5(v_1) = 0 & \textrm{c'est-à-dire} & v_{13} = 2 ,\\
-\alpha_5(v_2) = 0 & \textrm{c'est-à-dire} & v_{23} = 2 ,\\
-\alpha_5(v_3) = 0 & \textrm{c'est-à-dire} &
\frac{\mu_{13}}{v_{13}}+\frac{\mu_{23}}{v_{23}} = 2.
\end{array} \right.
$$

Donc $v_{13} = v_{23} = \frac{\mu_{13}}{v_{13}}+\frac{\mu_{23}}{v_{23}} = 2$, il vient que $\mu_{13}+\mu_{23} = 4$. Ce qui est absurde puisque par hypothèse $\theta_{23} + \theta_{31} < \frac{\pi}{2}$ et donc $\mu_{13}+\mu_{23} > 4$. Par conséquent $X_{\GG}= \varnothing$.

Il  reste le cas où le graphe sous-jacent à $\GG$ est un
tétraèdre combinatoire. La démonstration de ce point sera faite au
paragraphe \ref{paratetra}, proposition \ref{tetra}, où l'on
calcule l'espace des modules d'un tétraèdre étiqueté marqué
quelconque.
\end{proof}

\section{Démonstration des résultats}
\subsection{Plan de la démonstration des points 1) et 4)}\label{hyp}

\`A partir de maintenant tous les graphes étiquetés que l'on considère n'ont aucun 3-circuit affine ou sphérique, prismatique, avec angle droit, autrement dit ils vérifient $m(\GG) = 0$.

La démonstration se déroule en 6 étapes:

\begin{enumerate}
\item Comme les 3-circuits vont jouer un rôle essentiel dans la
compréhension des polyèdres miroirs, on va commencer par
s'intéresser aux triangles miroirs. Dans la partie \ref{1} on introduit l'invariant $R$ qui paramètre l'espace des modules d'un triangle combinatoire étiqueté.

\item L'idée pour comprendre l'espace $X_{\GG}$ est de découper les
polyèdres miroirs qui réalisent $\GG$ le long de leurs 3-circuits
prismatiques essentiels. On obtient ainsi des polyèdres miroirs qui réalisent
des blocs fondamentaux. Les démonstrations des lemmes nécessaires
à la construction des blocs fondamentaux et des plans pour
découper $P$ font l'objet de la partie \ref{2}.

\item Une bonne méthode pour comprendre les composantes connexes
de $X_{\GG}$ est de parler d'orientation partielle admissible de la
forêt $\F_{\GG}$ qui est une sous-forêt de l'arbre $\A_{\GG}$. Nous donnerons les définitions de tout cela dans la partie \ref{5}.

\item Pour que la bijection entre les composantes connexes de
$X_{\GG}$ et les orientations partielles admissibles de $\F_{\GG}$
soit la plus simple possible, il faut introduire un outil
technique: les systèmes puits-source de 3-circuits de $\GG$. Cela
sera fait dans la partie \ref{6}.

\item Il faut ensuite comprendre $X_{\GG}$ pour les blocs
fondamentaux et les prismes exceptionnels sans angle droit, ce qui
est fait dans la partie \ref{4}.

\item Enfin, une fois que l'on a découpé $P$ en blocs
fondamentaux, il faut comprendre comment on peut les recoller,
c'est la partie \ref{3}.
\end{enumerate}

\subsection{Les triangles miroirs}\label{1}

\begin{defi}
On appellera $\emph{triangle combinatoire}$ le graphe complet à 3
sommets. Un triangle $\emph{étiqueté}$ est la donnée d'un triangle
combinatoire $\mathcal{T}$ et pour chaque sommet $s$ de
$\mathcal{T}$ d'un réel $\theta \in ]0,\frac{\pi}{2}]$. Un
triangle combinatoire est dit $\emph{marqué}$ lorsque l'on a
numéroté ses arêtes de 1 à 3. On a une définition analogue à celle
des paragraphes précédents 
de triangle avec angle droit, sans angle droit, affine, sphérique,
hyperbolique, projectif, miroir, marqué et aussi d'espace des
modules de triangles miroirs marqués.
\end{defi}

\begin{rema}
Le marquage de $\mathcal{T}$  fournit une orientation
naturelle de $\mathcal{T}$.
\end{rema}

\begin{prop}\label{modtri}
Soient $\mathcal{T}$ un triangle étiqueté marqué, et
$X_{\mathcal{T}}$ l'espace des modules associé.
\begin{itemize}
\item Si $\mathcal{T}$ ne possède aucun sommet d'ordre 2 alors
$X_{\mathcal{T}}$ est difféomorphe à $\R$. De plus, si l'on
numérote les arêtes du triangle projectif miroir $T$ via le marquage de $\mathcal{T}$, alors
l'application:

$$
T \mapsto R_T =
\log\bigg(\frac{\alpha_1(v_2)\alpha_2(v_3)\alpha_3(v_1)}{\alpha_1(v_3)\alpha_3(v_2)\alpha_2(v_1)}\bigg)
$$

est un difféomorphisme de $X_{\mathcal{T}}$ sur $\R$.

\item Si $\mathcal{T}$ possède un sommet d'ordre 2 alors
$X_{\mathcal{T}}$ est un singleton.
\end{itemize}
\end{prop}

\begin{proof}
On ne démontre que le cas où tous les sommets sont d'ordre
différent de 2, l'autre cas se démontre par un calcul analogue.
Numérotons les faces du triangle étiqueté $\mathcal{T}$ de 1 à 3, et notons
$(e_i)_{i=1...3}$ la base canonique de $\R^3$. On peut supposer
que $\alpha_i = -e_i^*$. Le stabilisateur des classes projectives
$[\alpha_1],[\alpha_2],[\alpha_3]$ dans $\sss_3^{\pm}(\R)$ est l'ensemble
des matrices diagonales, on peut donc supposer que dans la base
$(e_i)_{i=1...3}$, le vecteur $v_1=(-2,\sqrt{\mu_{12}},\sqrt{\mu_{13}})$. Le
stabilisateur de $[\alpha_1],[\alpha_2],[\alpha_3]$ et $[v_1]$ est réduit
à l'identité.

On note encore pour tout $i,j = 1,...,3$, $v_{ij} = -\alpha_j(v_i)$, on rappelle que pour tout $i,j = 1,...,3$, les équations suivantes doit être vérifiées  $v_{ij}v_{ji} = \mu_{ij}$ (Définition \ref{miroir} et notations \ref{equa})

On obtient donc pour un certain $x \in \R^*_+$ :

$$
\left(
\begin{array}{c}
v_2 \\
v_3
\end{array}
\right)
=
\left(
\begin{array}{cccc}
\sqrt{\mu_{12}} & -2  & \sqrt{\mu_{23}}x \\
\sqrt{\mu_{13}} & \frac{\sqrt{\mu_{23}}}{x}&  -2
\end{array}
\right)
$$

Cela montre que l'application $T \mapsto
\log\Big(\frac{\alpha_1(v_2)\alpha_2(v_3)\alpha_3(v_1)}{\alpha_1(v_3)\alpha_3(v_2)\alpha_2(v_1)}\Big)=-2\log(x)$
est un difféomorphisme de $X_{\mathcal{T}}$ sur $\R$.
\end{proof}

Soit $\mathcal{T}$ un triangle étiqueté marqué, les quantités
suivantes vont se révéler cruciales dans la suite. On note $\Sigma
= \mu_{12} + \mu_{23} + \mu_{31}$, $p =
\sqrt{\mu_{12}\mu_{23}\mu_{31}}$. Lorsque $\mathcal{T}$ est sphérique (et sans angle droit, i.e $p>0$) on
définit aussi la quantité $r_{\mathcal{T}}=
2 \log\Big(\frac{4-\Sigma+((4-\Sigma)^2-p^2)^{\frac{1}{2}}}{p}\Big)
> 0$. Il n'est pas clair que la quantité $r_{\mathcal{T}}$ est bien définie mais on le montrera dans la démonstration du lemme \ref{triangle}.

\begin{lemm}\label{triangle}
Soient un triangle étiqueté marqué $\mathcal{T}$ avec au plus un angle droit et un triangle miroir marqué $T$ qui réalise $\mathcal{T}$. On note $\mathcal{B}$ une base de $\R^3$ dont l'orientation est opposée à celle de la base antéduale du triplet ($\alpha_1$, $\alpha_2$, $\alpha_3$).

\begin{tabular}{lll}
1)Si $\mathcal{T}$ est & hyperbolique & alors  $\, \det_{\mathcal{B}}(v_1,v_2,v_3) > 0$\\
& affine sans angle droit et $R_T \neq 0$ &\\
 & sphérique sans angle droit et $ |R_T| > r_{\mathcal{T}}$ & \\
&&\\
2) Si $\mathcal{T}$ est & affine avec angle droit & alors $\, \det_{\mathcal{B}}(v_1,v_2,v_3) = 0$.\\
& affine sans angle droit et $R_T=0$\\
& sphérique sans angle droit et $|R_T|=r_{\mathcal{T}}$& \\
&&\\
3) Si $\mathcal{T}$ est & sphérique sans angle droit et $|R_T| <
r_{\mathcal{T}}$ & alors $\, \det_{\mathcal{B}}(v_1,v_2,v_3) < 0$.\\
& sphérique avec angle droit & \\
\end{tabular}
\end{lemm}

\begin{rema}
On rappelle que les couples $(\alpha_i,v_i)_{i=1,...,3}$ sont bien définis à une constante multiplicative strictement positive près (notations \ref{equa}). Par conséquent, la valeur de la quantité $\text{det}_{\mathcal{B}}(v_1,v_2,v_3)$ n'a pas de sens, mais son signe et sa nullité ont un sens.
\end{rema}

\begin{proof}
Notons $(e_i)_{i=1...3}$ la base canonique de $\R^3$. On peut supposer que $\alpha_i = -e_i^*$. La base canonique de $\R^3$ possède bien une orientation opposée à celle de la base antéduale du triplet ($\alpha_1$, $\alpha_2$, $\alpha_3$). Commençons par le cas où $\Gamma$ est avec angle droit. On suppose que c'est l'angle entre l'arête 1 et l'arête 2 qui est égal à $\frac{\pi}{2}$. On obtient facilement (après multiplication par une matrice diagonale positive qui fixe donc les classes projectives de $\alpha_1,\alpha_2, \alpha_3$) que:

$$
\left(
\begin{array}{c}
v_1 \\
v_2 \\
v_3
\end{array}
\right)
=
\left(
\begin{array}{cccc}
-2             & 0               &\sqrt{\mu_{31}} \\
0               & -2                & \sqrt{\mu_{23}}\\
\sqrt{\mu_{31}} & \sqrt{\mu_{23}}  & -2
\end{array}
\right)
$$

Par conséquent, $D=\text{det}(v_1,v_2,v_3) = 2(\mu_{23} + \mu_{13} - 4) = 8(\cos^2(\theta_{23}) + \cos^2(\theta_{31})-1)$. Par conséquent, $D>0$ si $\theta_{23} + \theta_{31} < \frac{\pi}{2}$, nul si $\theta_{23} + \theta_{31} = \frac{\pi}{2}$ et $D<0$ si $\theta_{23} + \theta_{31} > \frac{\pi}{2}$.

On suppose à présent $\mathcal{T}$ sans angle droit et on pose $D
= \text{det}(v_1,v_2,v_3)$. Le lemme \ref{modtri} montre (après multiplication par une matrice diagonale positive qui fixe donc les classes projectives de $\alpha_1,\alpha_2, \alpha_3$) que
$$
\left(
\begin{array}{c}
v_1 \\
v_2 \\
v_3
\end{array}
\right)
=
\left(
\begin{array}{ccc}
-2              & \sqrt{\mu_{12}}                   & \sqrt{\mu_{13}}\\

\sqrt{\mu_{12}} & -2                                & \sqrt{\mu_{23}}e^{-\frac{R_T}{2}}\\

\sqrt{\mu_{13}} & \sqrt{\mu_{23}}e^{\frac{R_T}{2}} & -2 \\
\end{array}
\right)
$$

On a $D=\text{det}(v_1,v_2,v_3)=pe^{\frac{R_T}{2}}+pe^{\frac{-R_T}{2}}+2(\Sigma-4)$. On pose
$P(x) = p x^2 -2(4-\Sigma)x + p$; ainsi $D
=e^{-\frac{R_{T}}{2}}P(e^{\frac{R_{T}}{2}})$. Si $\Delta$ est le
discriminant de $P$, on a $\Delta = 4(4-\Sigma-p)(4-\Sigma+p)$.
Commençons par remarquer que comme $p>0$, si les racines de $P$
sont réelles alors elles sont toutes les deux du signe de
$4-\Sigma$. Ensuite, un peu de trigonométrie  montre que

$$
\textstyle{4-\Sigma-p =
-16\cos\Big(\frac{\theta_{12}+\theta_{23}+\theta_{31}}{2}\Big)\cos\Big(\frac{\theta_{12}-\theta_{23}+\theta_{31}}{2}\Big)\cos\Big(\frac{\theta_{12}+\theta_{23}-\theta_{31}}{2}\Big)\cos\Big(\frac{-\theta_{12}+\theta_{23}+\theta_{31}}{2}\Big)}.
$$

On remarque que comme les angles $\theta_{ij}$ sont aigus le produit des 3
derniers membres est strictement positif.

On distingue donc 3 cas:
\begin{enumerate}
\item Si $\mathcal{T}$ est sphérique, alors $4-\Sigma-p
> 0 $, par suite $\Delta > 0$ et $4-\Sigma >0$. Les racines $u_{\mathcal{T}},u_{\mathcal{T}}^{-1}$ de $P$ sont donc
réelles, strictement positives et inverses l'une de l'autre. On
peut supposer que $u_{\mathcal{T}} > 1$ et en remarquant que
$r_{\mathcal{T}}=2\log(u_{\mathcal{T}})$, on obtient donc que la quantité $r_{\mathcal{T}}$ est bien définie et strictement positive. La quantité $D$ est du signe de $|R_T| - r_{\mathcal{T}}$.

\item Si $\mathcal{T}$ est affine, alors $4-\Sigma-p = 0 $; par suite
$\Delta = 0$ et $4-\Sigma=p>0$, et ainsi 1 est racine double. La quantité $D$ est strictement positive si et seulement si $R_T \neq 0$, et nulle sinon.

\item Si $\mathcal{T}$ est hyperbolique, alors $4-\Sigma-p < 0 $. Alors
soit $4-\Sigma+p \leqslant 0$ et les racines de $P$ sont
strictement négatives, soit $4-\Sigma+p > 0$ et les racines de $P$
sont complexes. La quantité $D$ est toujours strictement positive.
\end{enumerate}
\end{proof}

\begin{nota}\label{defirgamma}
Soient $\GG$ un écimaèdre étiqueté marqué, $P$ un polyèdre miroir qui réalise $\GG$ et $\Gamma$ un 3-circuit orienté de
$\GG$. Alors $\Gamma$ définit un triangle étiqueté marqué $\mathcal{T}$ dont les
sommets sont les arêtes traversées par $\Gamma$, les arêtes de $\mathcal{T}$ sont
les faces traversées par $\Gamma$, les étiquettes portées par les
sommets sont les étiquettes portées par les arêtes de $\Gamma$, et on
marque les arêtes à l'aide de l'orientation de $\Gamma$.

Lorsque $\G$ est sans angle droit, on notera $R_{\Gamma}(P)$ le réel $R_T$ introduit dans l'énoncé du lemme
\ref{modtri}; si $\mathcal{T}$ est sphérique, on notera $r_{\Gamma}$ le réel
$r_{\mathcal{T}}$ introduit avant le lemme \ref{triangle} et si $\mathcal{T}$ est
affine, on pose $r_{\Gamma}=0$. Cette convention permet de ne pas distinguer le cas affine du cas sphérique.

On utilisera une autre convention, si un 3-circuit $\G$ est avec angle droit alors on pose $R_{\G}(P) = 0$. Cette convention permet de distinguer le cas sans angle droit du cas avec angle droit uniquement lorsque c'est nécessaire.

Enfin, plus généralement, si $(f_1,...,f_k)$ est une suite de faces de
$P$, on peut définir la quantité suivante:
$$
R_{(f_1,...,f_k)}(P)= \log \bigg( \frac{\alpha_{f_1} (v_{f_2})
\alpha_{f_2} (v_{f_3}) \cdot\cdot\cdot \alpha_{f_k}
(v_{f_1})}{\alpha_{f_1} (v_{f_k}) \alpha_{f_2} (v_{f_1})
\cdot\cdot\cdot \alpha_{f_k} (v_{f_{k-1}})} \bigg).
$$
Ainsi, lorsque $\Gamma$ est un 3-circuit orienté, il définit
naturellement une suite de faces et l'on retrouve la définition
précédente.

Lorsqu'il n'y a pas d'ambiguïté, on allégera les notations en
notant $R_{(f_1,...,f_k)}(P)=R_{(f_1,...,f_k)}$ et
$R_{\Gamma}(P)=R_{\Gamma}$.
\end{nota}

\subsection{Les lemmes de coupe}\label{2}

On cherche tout d'abord à construire les blocs fondamentaux
miroirs en coupant un tétraèdre miroir le long de ses 3-circuits.
C'est l'objet des lemmes \ref{ecimage2}, \ref{ecimage3} et \ref{nonchevau}, des deux prochains paragraphes. Ensuite,
on se donne un polyèdre miroir $P$ et on souhaite découper $P$ le
long de ses 3-circuits prismatiques essentiels pour obtenir des blocs
fondamentaux. Nous allons avoir besoin des définitions suivantes:

\begin{defi}
Soient $P$ un polyèdre miroir et $\Gamma$ un 3-circuit de $P$ qui
traverse les faces $r,s,t$ de $P$. On désignera par $\Pi_{\Gamma}$
le sous-espace projectif engendré par les points polaires
$[v_r],[v_s],[v_t]$. On dira que \emph{$\Pi_{\Gamma}$ coupe $P$ le
long de $\Gamma$} lorsque $\Pi_{\Gamma}$ est un plan et
l'intersection de $\Pi_{\Gamma}$ avec les arêtes de $P$ est
incluse dans les arêtes ouvertes le long desquelles se rencontrent les faces $r,s,t$.
\end{defi}


\begin{rema}\label{uniqueref}
Soient $P$ un polyèdre miroir de $X_{\GG}$ et $\Gamma$ un 3-circuit prismatique non essentiel de $P$ (définition \ref{combiessens}). On numérote de 1 à 3 les faces traversées par $\G$, et on numérote 4 la face triangulaire dont $\G$ fait le tour. Nous allons montrer que si $\Pi_{\Gamma}$ est un plan alors la réflexion $\sigma_4$ est entièrement déterminée par les réflexions $\sigma_1,\sigma_2,\sigma_3$.

En effet, comme $\G$ n'est pas essentiel, les arêtes de la face 4 sont toutes d'ordre 2 (en particulier $\sigma_4$ commute avec $\sigma_1$, $\sigma_2$ et $\sigma_3$), par conséquent les points $v_1$, $v_2$ et $v_3$ vérifient $\alpha_4(v_1)=\alpha_4(v_2)=\alpha_4(v_3)=0$ (définition \ref{miroir} et notations \ref{equa}). Ainsi, le plan projectif définissant la face $4$ est le plan projectif $\Pi_{\G}$ qui est engendré par les points polaires $[v_1]$, $[v_2]$ et $[v_3]$. La polaire $[v_4]$ de la face 4 vérifie les équations $\alpha_1(v_4)=\alpha_2(v_4)=\alpha_3(v_4)=0$, cette intersection est réduite à deux points car $P$ est un polyèdre convexe. La polaire $[v_4]$ doit aussi vérifier que $\alpha_4(v_4)=2$, il vient qu'un seul des deux points de l'intersection des plans projectifs associés aux faces 1, 2 et 3 vérifie ces propriétés.
\end{rema}

Nous allons montrer que, pour tout polyèdre miroir $P$ de $X_{\GG}$ et
tout 3-circuit prismatique essentiel $\Gamma$, $\Pi_{\Gamma}$ est un plan et que ce dernier coupe $P$ le long de $\Gamma$. On peut donc définir
la coupe $\emph{droite}$ $P^d_{\Gamma}$ (resp. $\emph{gauche}$
$P^g_{\Gamma}$) de $P$ le long de $\Gamma$ comme le polyèdre
miroir obtenu en retirant les faces et les réflexions par rapport
aux faces à gauche (resp. droite) de $\Gamma$ et en ajoutant
l'unique (remarque \ref{uniqueref}) réflexion par rapport au plan $\Pi_{\Gamma}$ qui commute
avec les réflexions associées aux faces traversées par $\Gamma$.

\begin{rema}
Lorsqu'il n'y a pas d'ambiguïté sur le 3-circuit $\Gamma$, on
allégera les notations $P_{\Gamma}^d$ et $P_{\Gamma}^g$ en notant
$P_{\Gamma}^d = P^d$ et $P_{\Gamma}^g=P^g$.
\end{rema}

\subsubsection{Lemme d'écimage}

\begin{lemm}\label{ecimage2}
Soient $P$ un polyèdre miroir et $\Gamma$ un 3-circuit orienté
de $P$. Le sous-espace $\Pi_{\Gamma}$ est un plan projectif sauf si les
conditions 1)2)3) ou 1)2)3)$'$ ou 1)2)3)$''$ sont réunies.

\begin{tabular}{cl}
1) & Le polyèdre $P$ est un tétraèdre ou un prisme triangulaire.\\
2) & Les arêtes des faces qui ne sont pas traversées par $\Gamma$
sont d'ordre 2. \\
3) & Le 3-circuit $\Gamma$ est affine sans angle droit et $R_{\Gamma}=0$.\\
3)$'$ & Le 3-circuit  $\Gamma$ est sphérique sans angle droit et
$|R_{\Gamma}|=r_{\Gamma}$.\\
3)$''$ & Le 3-circuit $\Gamma$ est affine avec angle droit.
\end{tabular}
\end{lemm}

\begin{proof}
Soient $s_1,s_2,s_3$ les faces de $P$ traversées par $\Gamma$.
Soit $s_4$ une face quelconque de $P$ différente de $s_1,s_2,s_3$.
On peut supposer que les faces $(s_i)_{i=1,...,4}$ sont définies
par les formes linéaires $(-e_i^*)_{i=1,...,4}$, où
$(e_i)_{i=1,...,4}$ est une base de $\R^4$, car $P$ est un
polyèdre convexe. Le sous-espace $\Pi_{\Gamma}$ n'est pas un plan projectif si et seulement
si la matrice $M=(-\alpha_j(v_i))_{i=1...3, \, j=1...4}$ est de
rang 2, donc si et seulement si tous les mineurs de taille 3
extraits de la matrice $M$ sont nuls. On a:

$$
M =\left( \begin{array}{cccc}
-2     & v_{12} & v_{13} & v_{14}\\
v_{21} & -2     & v_{23} & v_{24}\\
v_{31} & v_{32} & -2     & v_{34}\\
\end{array}\right)
$$

où les $v_{ij}$ désignent les coordonnées de $(v_i)_{i=1...3}$
dans la base $(e_j)_{j=1...4}$.

Le mineur obtenu en rayant la dernière colonne  donne le
déterminant du lemme \ref{triangle}. Donc, si $\Pi_{\Gamma}$ n'est
pas un plan, alors l'une des conditions 3), 3)$'$ ou 3)$''$ est
réalisée.

Calculons à présent le mineur $D_1$ obtenu en rayant la 1ère
colonne:
$$
D_1=v_{12}(v_{23}v_{34}+2v_{24})+v_{13}(2v_{34}+v_{32}v_{24})+v_{14}(4-\mu_{23}).
$$
Le lemme \ref{lemme1} montre que $v_{ij} \geqslant 0$ pour $i \neq
j$; de plus $\mu_{23}< 4$ et donc $v_{14}=0$ quand $D_1=0$.

De la même façon, le calcul des mineurs obtenus en rayant la 2ème
ou la 3ème colonne montre que $v_{24}=v_{34}=0$. Le lemme
\ref{lemme1} montre donc que la face $s_4$ est adjacente à $s_1$,
$s_2$ et $s_3$ . Les conditions 1) et 2) sont alors réalisées car
il y a au plus 2 faces adjacentes simultanément à $s_1$, $s_2$ et
$s_3$.

Les mêmes calculs montrent que si les conditions 1)2)3) ou 1)2)3)$'$ ou 1)2)3)$''$ sont réunies alors $\Pi_{\G}$ est une droite projective.
\end{proof}

\begin{lemm}\label{ecimage3}
Soient $P$ un polyèdre miroir et $\Gamma$ un 3-circuit qui entoure
un sommet $v$. On suppose que $P$ n'est pas un tétraèdre dont la
face opposée à $v$ ne possède que des arêtes d'ordre 2. On suppose
aussi que $\Pi_{\Gamma}$ est un plan projectif et que $\G$ possède au plus un seul angle droit. On note $P_0$ le polyèdre projectif sous-jacent à $P$.

\begin{tabular}{lll}
1)Si $\Gamma$ est & hyperbolique & alors $\Pi_{\Gamma}$ coupe $P_0$ le long de $\Gamma$.\\
& affine sans angle droit et $R_{\Gamma} \neq 0$ &\\
 & sphérique sans angle droit et $ |R_{\Gamma}| > r_{\Gamma}$ & \\
&&\\
2) Si $\Gamma$ est & affine avec angle droit & alors $\Pi_{\Gamma} \cap P_0 = \{v\}$.\\
& affine sans angle droit et $R_{\Gamma}=0$\\
& sphérique sans angle droit et $|R_{\Gamma}|=r_{\Gamma}$& \\
&&\\
3) Si $\Gamma$ est & sphérique sans angle droit et $|R_{\Gamma}| <
r_{\Gamma}$ & alors $\Pi_{\Gamma} \cap P_0 = \varnothing$.\\
& sphérique avec angle droit & \\
\end{tabular}
\end{lemm}

\begin{rema}\label{apasrate}
Pour alléger l'énoncé du lemme \ref{ecimage3} nous apportons une précision au cas 3) sous la forme d'une remarque. On numérote les faces de $\G$ de 1 à 3. On note $e_i$ la $i$-ème arête fermée traversée par $\G$ et $C$ le convexe polyédral obtenu à partir de $P_0$ en oubliant les inégalités associées aux faces traversées par $\Gamma$ (ces faces partagent un sommet). On note $L_i$ la droite engendrée par $e_i$ et $D_i$ l'intersection $C \cap L_i$. On note $M_i$ le point d'intersection du plan $\Pi_{\G}$ et $D_i$, pour $i=1,...,3$. Il n'est pas évident que les points $(M_i)_{i=1,...,3}$ sont bien définis mais nous montrerons que c'est le cas. Nous allons montrer que dans le cas 3), les points $M_i$ sont dans l'intérieur de $D_i \setminus e_i$.
\end{rema}

\begin{rema}
Avant de démontrer ce lemme, on pourra remarquer que les conditions 1)2)3) sont celles du lemme \ref{triangle}.
\end{rema}

\begin{proof}
Soient $(s_i)_{i=1,...,3}$ les faces de $P$ traversées par
$\Gamma$ et soit $s_4$ une face quelconque de $P$. Soit $\Pi_{\Gamma}$ le sous-espace
engendré par les points polaires $[v_1],[v_2]$ et $[v_3]$. On a
supposé que $\Pi_{\Gamma}$ est un plan projectif.

Comme on travaille dans le revêtement à deux feuillets de $\PP(\R^4)$ l'intersections de trois plans en positions génériques est la réunion de deux points opposés. On définit les six points $M_1^{\pm},M_2^{\pm}$ et $M_3^{\pm}$, de la façon suivante:
$\{ M_1^+,M_1^- \} = \Pi_{\Gamma} \cap \{\alpha_{2} =0\} \cap \{\alpha_{3}=0\}$,
$\{ M_2^+,M_2^- \} = \Pi_{\Gamma} \cap \{\alpha_{3} =0\} \cap \{\alpha_{1}=0\}$,
$\{ M_3^+,M_3^- \} = \Pi_{\Gamma} \cap \{\alpha_{1} =0\} \cap \{\alpha_{2}=0\}$.

Commençons par vérifier que les points $M_1^{\pm},M_2^{\pm}$ et $M_3^{\pm}$ sont bien définis. On ne traite que les points $M_1^{\pm}$. Il s'agit de montrer que l'ensemble des solutions des équations suivantes en les inconnues $x,y,z$ est un espace vectoriel de dimension 1.

$$
\left\{
\begin{array}{cccc}
\alpha_2(x v_1 + y v_2 + z v_3) = 0\\
\alpha_3(x v_1 + y v_2 + z v_3) = 0\\
\end{array}
\right.
$$

Le calcul explicite de ce système donne:

$$
\left\{
\begin{array}{ccccccc}
x \alpha_2(v_1) & + & 2y              & + &  z \alpha_2(v_3) & = 0\\
x \alpha_3(v_1) & + & y \alpha_3(v_2) & + &  2z              & = 0\\
\end{array}
\right.
$$

Ce système est de rang 2 puisque $4 - \alpha_2(v_3)\alpha_3(v_2) = 4-\mu_{23} > 0$. Les points $M_1^{\pm},M_2^{\pm}$ et $M_3^{\pm}$ sont donc bien définis.

On peut supposer que les faces $(s_i)_{i=1,...,4}$ sont définies par
$\alpha_{s_i}=-e_i^*$ pour $i=1...4$ dans une base
$(e_i)_{i=1...4}$ de $\R^4$. Ainsi, le plan projectif $\Pi_{\G}$ est défini par
une équation de la forme $\{ ax_1+bx_2+cx_3+dx_4=0 \}$, où
$a,b,c,d \in \R$. Il faut prendre garde au fait que les quantités $a,b,c,d$ sont bien définies seulement à une constante multiplicative près.

Introduisons la quantité:

$$D=\left| \begin{array}{c c c} -2 & v_{21} & v_{31} \\
v_{12} & -2 & v_{32} \\
v_{13} & v_{23} & -2
\end{array} \right|,$$

\noindent où les $v_{ij}$ désignent les coordonnées de $(v_i)_{i=1...3}$
dans la base $(e_j)_{j=1...4}$. La quantité $D$ est le déterminant du lemme \ref{triangle}.

Le plan projectif $\Pi_{\G}$ d'équation $ax_1+bx_2+cx_3+dx_4=0$ est le plan engendré par les vecteurs $v_1$, $v_2$ et $v_3$. Par conséquent, il existe un $\lambda \neq 0$ tel que pour tout $x,y,z,t \in \R$ on a:

$$
\left|
\begin{array}{c c c c}
 -2    & v_{21} & v_{31} & x\\
v_{12} & -2     & v_{32} & y\\
v_{13} & v_{23} & -2     & z\\
v_{14} & v_{24} & v_{34} & t
\end{array}
\right|
=\lambda(ax+by+cz+dt)
$$

Nous allons distinguer les cas $D > 0$, $D=0$ et $D < 0$.

Commençons par supposer que $D > 0$. D'après le lemme \ref{triangle}, cette hypothèse revient à se placer dans le cas $1)$, on a alors $D = \lambda d \neq 0$ par la formule de Cramer. Pour fixer les idées, on suppose à présent que $d=-1$.

Les formules de Cramer  permettent de calculer
$a$. On a
$-a\lambda = aD = v_{14}(4-\mu_{23})+v_{24}(2v_{12}+v_{13}v_{32})+v_{34}(v_{12}v_{23}+2v_{13})$.

Par hypothèse $\G$ possède au plus un angle droit donc au plus une seule des quantités $v_{12},
\, v_{13}, \, v_{23}$ est nulle. De plus, on a supposé que $P$
n'était pas un tétraèdre où la face opposée à $v$ ne possédait que
des arêtes d'ordre 2. Par conséquent, l'une des arêtes de la face
4 n'est pas d'ordre 2. Il vient donc que $v_{14}$ ou $v_{24}$ ou
$v_{34}$ est strictement positif. Il est donc clair à présent que
$aD > 0$. De la même façon, on obtient
que $b$ et $c$ sont strictement positifs.

On peut donc à présent écrire explicitement les coordonnées des points $M_1^{\pm},M_2^{\pm}$ et $M_3^{\pm}$ dans la base $(e_1, ..., e_4)$. On a: $M_1^+ =[-d:0:0:a]$, $M_2^+ = [0:-d:0:b]$, $M_3^+=[0:0:-d:c]$ et $M_1^- =[d:0:0:-a]$, $M_2^- = [0:d:0:-b]$ et $M_3^-=[0:0:d:-c]$.

En résumé $\alpha_i(M_i^+) = d < 0$ et $(\alpha_4(M_1^+),\alpha_4(M_2^+),\alpha_4(M_3^+))=(-a,-b,-c)<0$. Ainsi, on obtient que le point $M_i^+$ est à l'intérieur de l'arête adjacente aux faces $i+1$ et $i+2$ [mod 3] du tétraèdre formé par les faces 1,2,3 et 4. Mais la face $s_4$ est une face quelconque de $P$. Par conséquent, le point $M_i^+$ est à l'intérieur de l'arête adjacente aux faces $i+1$ et $i+2$ [mod 3] du polyèdre $P_{0}$. Ce qui entraîne que $\Pi_{\Gamma}$ coupe $P_0$ le long de $\Gamma$.

\`{A} présent, si $D< 0$ alors le lemme \ref{triangle} montre que l'on est dans le cas $3)$. On suppose encore que $d=-1$. On a encore l'égalité $aD=v_{14}(4-\mu_{23})+v_{24}(2v_{12}+v_{13}v_{32})+v_{34}(v_{12}v_{23}+2v_{13}) > 0$. On obtient donc que $a,b$ et $c$ sont strictement négatifs.

On écrit explicitement les coordonnées des points $M_1^{\pm},M_2^{\pm}$ et $M_3^{\pm}$ dans la base $(e_1, ..., e_4)$. On a: $M_1^+ =[-d:0:0:a]$, $M_2^+ = [0:-d:0:b]$, $M_3^+=[0:0:-d:c]$ et $M_1^- =[d:0:0:-a]$, $M_2^- = [0:d:0:-b]$ et $M_3^-=[0:0:d:-c]$.

Cette fois-ci on obtient que $\alpha_i(M_i^-) = -d >0$ et $(\alpha_4(M_1^-),\alpha_4(M_2^-),\alpha_4(M_3^-))=(a,b,c)<0$. Ainsi, le point $M_i^-$ n'appartient pas à l'arête adjacente aux faces $i+1$ et $i+2$ [mod 3] du tétraèdre formé par les faces 1,2,3 et 4 mais le point $M_i^-$ vérifie  $\alpha_4(M_i^-)<0$. La face $s_4$ est une face quelconque de $P$. Par conséquent, le point $M_i^-$ n'appartient pas à l'arête adjacente aux faces $i+1$ et $i+2$ [mod 3] du polyèdre $P_0$ mais il appartient au convexe $C=p(\{ x \in \R^4 \, | \,  \alpha_s(x) \leqslant 0, \,  s\neq 1,2,3 \})$ (le point $M_i^-$ est le point $M_i$ de la remarque \ref{apasrate}). Comme la partie $C$ est incluse dans une carte affine, il vient que le polyèdre $P$ et les points $M_i^-$ sont eux aussi inclus dans une carte affine. Ce qui entraîne que $\Pi_{\Gamma} \cap P_0 = \varnothing$ et la remarque \ref{apasrate}.

Enfin, si $D= 0$ alors le lemme \ref{triangle} montre que l'on est dans le cas $2)$, mais dans ce cas, on a $d=0$. Rappelons que $v$ est le point dont $\G$ fait le tour. Comme $d=0$ l'intersection $\Pi_{\Gamma} \cap \{\alpha_{1} =0\} \cap \{\alpha_{2} =0\} \cap \{\alpha_{3}=0\} = \{ v,-v \}$ et par suite les ensembles $\{M_i^+,M_i^-\}$ (pour $i=1,...,3$) sont confondus et égaux à $\{\pm v \}$. Ce qui conclut la démonstration.
\end{proof}

\begin{rema}
Il est important de noter que cette étude est exhaustive si les hypothèses du lemme \ref{ecimage2} sont vérifiées.
\end{rema}

\subsubsection{Lemme de coupe le long d'un 3-circuit prismatique essentiel}

\begin{lemm}\label{coupris}
Soient $P$ un polyèdre miroir et $\Gamma$ un 3-circuit
prismatique essentiel. Alors $\Pi_{\Gamma}$ coupe $P$ le long de $\Gamma$.
\end{lemm}

\begin{proof}
Comme $P$ possède un 3-circuit essentiel, $P$ n'est ni un tétraèdre, ni un prisme exceptionnel et par conséquent, le lemme \ref{ecimage2} permet d'affirmer que $\Pi_{\Gamma}$ est un plan. On suppose que les faces traversées par $\Gamma$ sont numérotées de 1 à 3. On note $S^d$ (resp. $S^g$) l'ensemble des faces du polyèdre $P$ strictement à droite (resp. strictement à gauche) de $\Gamma$ (on considère que les faces 1,2 et 3 ne sont ni strictement à droite ni strictement à gauche de $\Gamma$).


On note $e_i$ l'arête adjacente aux faces $i+1$ et $i+2$ [mod 3]. On note $D_i$ l'intersection entre la droite engendrée par l'arête $e_i$ et le convexe $C^d = p(\{ x \in \R^4 \,|\, \forall s \in S^d,\, \alpha_s(x) < 0 \})$. On note $M_i$ le point d'intersection du plan $\Pi_{\G}$ avec $D_i$ [mod 3]. Le lemme \ref{ecimage3} est exhaustif, ainsi en utilisant la remarque \ref{apasrate}, seul l'un de ces trois cas est possible:

\begin{itemize}
\item Le point $M_i$ appartient à l'intérieur de $e_i$, pour $i=1,...,3$.

\item Les points $M_i$ sont tous égaux au sommet $e_1 \cap e_2 \cap e_3$.

\item Les points $M_i$ appartiennent à l'intérieur de $D_i \smallsetminus e_i$, pour $i=1,...,3$.
\end{itemize}

Il est important de remarquer que $D_i$ prolonge l'arête $e_i$ uniquement à un bout, à savoir le bout à gauche de $\G$.

Comme $\Gamma$ est prismatique essentiel, on peut faire la même chose avec le convexe $C^g = p(\{ x \in \R^4 \,|\, \forall s \in S^g,\, \alpha_s(x) < 0 \})$. Par conséquent, seul le premier cas peut être réalisé, ainsi $\Pi_{\Gamma}$ coupe $P$ le long de $\Gamma$.
\end{proof}

\subsubsection{Lemme de non-chevauchement}

\begin{lemm}\label{nonchevau}
Soient $P$ un polyèdre miroir, $u,v$ deux sommets de $P$ reliés
par une arête $e$. Soit $\Gamma_u$ (resp. $\Gamma_v$) un 3-circuit
qui entoure $u$ (resp. $v$). On suppose que $\Pi_{\Gamma_u}$ (resp.
$\Pi_{\Gamma_v}$) coupe $P$ le long de $\Gamma_u$ (resp.
$\Gamma_v$). Si $P_0$ désigne le polyèdre projectif sous-jacent à
$P$, alors $P_0 \cap \Pi_{\Gamma_u} \cap \Pi_{\Gamma_v} =
\varnothing$.
\end{lemm}

\begin{proof}
On peut supposer que le sommet $u$ est à droite de $\G_u$ et on note $P^g_{\G_u}$ la coupe gauche de $P$ le long de $\Pi_{\G_u}$. Pour montrer ce lemme il suffit de montrer que $\Pi_{\G_v}$ coupe $P^g_{\G_u}$ le long de $\G_v$.

D'après le lemme \ref{ecimage3}, il suffit de vérifier que le 3-circuit $\G_v$ vérifie la condition 1) du lemme \ref{ecimage3}. Or le 3-circuit $\Pi_{\G_v}$ coupe $P$ le long de $\G_v$, il vérifie donc la condition 1) du lemme \ref{ecimage3}.

\end{proof}

\subsection{La forêt $\F_{\GG}$ et son orientation}\label{5}

\subsubsection{Classification des arêtes de $\A_{\GG}$}\label{classi}

Nous sommes à présent en mesure de définir l'ensemble $I_{\Gamma}^{\GG}$ des valeurs possibles a priori
pour la quantité $R_{\Gamma}$.

\[
\begin{tabular}{lll}

1. $\left.
\begin{tabular}{l}
$\Gamma$ est affine  ou sphérique\\
\end{tabular}
\right.$
&
\begin{tabular}{ll}
avec angle droit  & prismatique\\
\end{tabular}
&
\begin{tabular}{ll}
alors&   $I_{\Gamma}^{\GG}:=\varnothing$\\
\end{tabular}

\\
\\
2. $\left\{
\begin{tabular}{l}
$\Gamma$ est hyperbolique\\
$\Gamma$ est affine ou sphérique
\end{tabular}
\right.$
&
\begin{tabular}{ll}
avec angle droit  & quelconque\\
avec angle droit  &  non prismatique
\end{tabular}
&
\begin{tabular}{ll}
alors & $I_{\Gamma}^{\GG}:=\{0\}$\\
\end{tabular}

\\
\\
3. $\left\{
\begin{tabular}{l}
$\Gamma$ est hyperbolique\\
$\Gamma$ est affine ou sphérique
\end{tabular}
\right.$
&
\begin{tabular}{ll}
sans angle droit & quelconque\\
sans angle droit  &  non prismatique
\end{tabular}
&
\begin{tabular}{ll}
alors &    $I_{\Gamma}^{\GG}:=\R$\\
\end{tabular}

\\
\\
4. $\left.
\begin{tabular}{l}
$\Gamma$ est affine ou sphérique
\end{tabular}
\right.$
&
\begin{tabular}{llll}
sans angle droit  & prismatique
\end{tabular}
&
\begin{tabular}{ll}
alors &  $I_{\Gamma}^{\GG}:=\R \smallsetminus [-r_{\Gamma},r_{\Gamma}]$\\
\end{tabular}
\end{tabular}
\]

\begin{rema}
On rappelle que lorsque le 3-circuit $\Gamma$ est affine sans angle droit on utilise la convention $r_{\Gamma}=0$ (voir notations \ref{defirgamma}). On a aussi convenu que $R_{\G}=0$ lorsque $\G$ est à angle droit.
\end{rema}

Ces notations sont justifiées par les propositions suivantes:

\begin{prop}\label{propoublie2}
Soient $\GG$ un écimaèdre étiqueté qui vérifie $m(\GG)=0$ et $P$ un polyèdre miroir qui réalise $\GG$. Soit $v$ un sommet de $\GG$, on note $\G$ le 3-circuit orienté qui entoure le sommet $v$ et tel que la coupe gauche $\GG^g_{\G}$ ne contienne pas le sommet $v$. Alors, $\Pi_{\G}$ coupe $P$ le long de $\G$ si et seulement si $R_{\G}(P) \in I_{\G}^{\GG^g_{\G}}$.
\end{prop}

\begin{proof}
Cette proposition est essentiellement une relecture du lemme \ref{ecimage3}. Comme $\G$ devient prismatique dans $\GG^g_{\G}$, l'intervalle $I_{\G}^{\GG^g_{\G}}$ est égale à l'intervalle $I_{\G}^{\GG}$ moins un segment $S$. Il suffit de maintenant, par une étude exhaustive, de voir que $S$ correspond aux valeurs de $R_{\G}$ pour lesquels le plan $\Pi_{\G}$ ne coupe pas $P$ le long de $\G$. Remarquons que $\G$ (vu dans $\GG^g$) peut être de 3 types de la classification énoncée au paragraphe \ref{classi}: 2.(première ligne) ; 3.(première ligne) ; et 4.

Le lemme \ref{ecimage3} montre que si $\Pi_{\G}$ coupe $P$ le long de $\G$ alors $R_{\G}(P) \in I_{\G}^{\GG^g_{\G}}$.

Inversement, le lemme \ref{ecimage2} montre que si $R_{\G}(P) \in I_{\G}^{\GG^g_{\G}}$ alors $\Pi_{\G}$ est un plan. Le lemme \ref{ecimage3} montre que si $\Pi_{\G}$ est un plan et que $R_{\G}(P) \in I_{\G}^{\GG^g_{\G}}$ alors $\Pi_{\G}$ coupe $P$ le long de $\G$.

%
%
\end{proof}

\begin{prop}\label{propoublie}
Soient $\GG$ un écimaèdre étiqueté qui vérifie $m(\GG) = 0$ et qui n'est pas un prisme exceptionnel, $P$ un polyèdre miroir qui réalise $\GG$ et $\G$ un 3-circuit orienté de $\GG$. Alors $R_{\G}(P) \in I_{\G}^{\GG}$.
\end{prop}

\begin{proof}
Il faut distinguer les 4 cas de la définition de $I_{\G}^{\GG}$. Le cas 1) ne peut pas apparaître puisqu'on a supposé $m(\GG) = 0$. Si $\G$ est du type 2 de la classification alors on a la convention $R_{\G}(P) = 0$ (notations \ref{defirgamma}). Si $\G$ est de type 3 alors il n'y a rien à montrer. Enfin, si $\G$ est de type 4, aucun des deux polyèdres combinatoires étiquetés $\GG^g_{\G}$ ou $\GG^d_{\G}$ n'est un tétraèdre car $\G$ est prismatique. Notons $\mathcal{Q}$ le polyèdre combinatoire étiqueté obtenu en collant un tétraèdre sur la face triangulaire de $\GG^g_{\G}$ dont $\G$ fait le tour (autrement dit $\mathcal{Q}$ est le polyèdre combinatoire obtenu à partir de $\GG^g_{\G}$ en oubliant l'inégalité associée à la face dont $\G$ fait le tour dans $\GG^g_{\G}$). On note $Q$ le polyèdre obtenu à partir de $P$ par le même procédé.

Si le 3-circuit $\G$ de $\GG$ est essentiel alors le lemme \ref{coupris} montre que $\Pi_{\G}$ coupe $P$ le long de $\G$. En particulier, $\Pi_{\G}$ coupe $Q$ le long de $\G$.

Si le 3-circuit $\G$ de $\GG$ n'est pas essentiel, cela signifie que l'un des deux polyèdres combinatoires $\GG^g_{\G}$ ou $\GG^d_{\G}$ est un prisme exceptionnel. On peut supposer que $\GG^g_{\G}$ n'est pas un prisme exceptionnel. On va montrer que $\Pi_{\G}$ coupe $Q$ le long de $\G$. Le lemme \ref{ecimage2} montre que $\Pi_{\G}$ est un plan car $\G$ est prismatique et $\GG$ n'est pas un prisme exceptionnel. Donc le plan $\Pi_{\G}$ et le plan engendré par la face de $P$ dont $\G$ fait le tour sont égaux, il vient que $\Pi_{\G}$ coupe $Q$ le long de $\G$.

Ainsi, le lemme \ref{ecimage3} appliqué à $Q$ montre que $R_{\G}(P) \in I_{\G}^{\GG}$.
\end{proof}

Si $\GG$ est un écimaèdre étiqueté, alors les arêtes de $\A_{\GG}$
sont par définition en bijection avec les 3-circuits de $\GG$. On a
donc une définition naturelle d'arête affine, sphérique,
hyperbolique, sans angle droit, avec angle droit, prismatique et
enfin, non prismatique.

\begin{rema}\label{pasvide}
On rappelle que l'on suppose toujours, depuis le paragraphe
\ref{hyp}, que tout écimaèdre étiqueté vérifie $m(\GG)=0$; il n'y a
donc pas d'arête du type 1 de la classification que l'on vient de
donner dans $\A_{\GG}$.
\end{rema}

\begin{defi}
Soit $\GG$ un écimaèdre étiqueté. \emph{La forêt associée à $\GG$}
est la forêt $\F_{\GG}$ obtenue en supprimant les arêtes avec angle droit de
l'arbre $\A_{\GG}$ (i.e. les arêtes de type 2 de la classification). Les arêtes (resp. les 3-circuits de $\GG$) de
$\F_{\GG}$ se séparent en deux familles via la classification
énoncée plus haut:

\begin{itemize}
\item les arêtes (resp. les 3-circuits) hyperboliques sans angle
droit quelconques et les arêtes (resp. les 3-circuits) affines ou sphériques sans angle droit non-prismatiques;

\item les arêtes (resp. les 3-circuits) affines ou sphériques sans angle
droit prismatiques qu'on désignera dorénavant par le
terme \emph{spéciales} (resp. \emph{spéciaux}).
\end{itemize}
\end{defi}

\subsubsection{Orientation de la forêt $\F_{\GG}$}

\begin{defi}
Soit $\mathcal{F}$ une forêt orientée. On dira qu'un sommet $s$ de
$\mathcal{F}$ est un \emph{puits} (resp. une \emph{source}) lorsque
toutes les arêtes incidentes à $s$ possèdent le même but (resp. la
même source).
\end{defi}

\begin{defi}
Soient $\GG$ un écimaèdre étiqueté, et $\F_{\GG}$ la forêt associée.
\begin{itemize}
\item Une $\emph{orientation globale admissible}$ de $\F_{\GG}$ est
une orientation de $\F_{\GG}$ qui ne contient aucun sommet de
valence 4 qui soit un puits ou une source.

\item Une $\emph{orientation partielle}$ de $\F_{\GG}$ est une
orientation de toutes les arêtes spéciales de $\F_{\GG}$.

\item Une $\emph{orientation partielle admissible}$ de $\F_{\GG}$
est une orientation partielle de $\F_{\GG}$ telle qu'il existe une
orientation globale admissible de $\F_{\GG}$ qui la prolonge.
\end{itemize}
\end{defi}

Nous allons voir que si $\GG$ est un écimaèdre étiqueté alors les
composantes connexes de $X_{\GG}$ sont en bijection avec les
orientations partielles admissibles de la forêt $\F_{\GG}$. Cette
bijection est construite à l'aide du signe des quantités
$R_{\Gamma}(P)$, où $\Gamma$ est un 3-circuit orienté spécial de
$\GG$ et $P$ un polyèdre projectif miroir qui réalise $\GG$. Pour
uniformiser la construction de cette bijection, il faut choisir
correctement un système de 3-circuits orientés. Nous appellerons
ces systèmes les systèmes puits-source de 3-circuits et nous
allons les définir dès à présent.

\subsection{Système puits-source de $\A_{\GG}$}\label{6}

\subsubsection{Orientation d'une arête de $\A_{\GG}$ induite par l'orientation d'un 3-circuit de $\GG$}

Chaque arête $e$ de $\A_{\GG}$ définit une tripartition des arêtes
de $\GG$: $\{e\}$ et les deux composantes connexes de
$\A_{\GG}-\{e\}$. Le choix d'une orientation de l'arête $e$ permet
de distinguer les deux composantes connexes de $\A_{\GG}-\{e\}$. On
a celle donnée par le but de l'arête orientée $e$ et celle donnée
par sa source.

De même, chaque 3-circuit $\Gamma$ de $\GG$ définit une
tripartition des 3-circuits de $\GG$. Le choix d'une orientation du
3-circuit $\Gamma$ permet de les distinguer. On a le 3-circuit
$\Gamma$, les 3-circuits de $\GG$ à droite de $\Gamma$ et les
3-circuits de $\GG$ à gauche de $\Gamma$.

Il est évident que la bijection entre les 3-circuits de $\GG$ et
les arêtes de $\A_{\GG}$ respecte cette tripartition. On a donc une
définition naturelle d'orientation d'une arête de $\A_{\GG}$
induite par l'orientation d'un 3-circuit de $\GG$:

\begin{defi}
Soient $\GG$ un graphe étiqueté, $\Gamma$ un 3-circuit orienté de
$\GG$ et $e$ l'arête correspondante de l'arbre $\A_{\GG}$. On dira
que $e$ est orienté dans le \emph{sens} (resp. \emph{sens
contraire}) de $\Gamma$ lorsque la composante connexe de $\A_{\GG}-\{e\}$ donnée par le but de l'arête orientée $e$ correspond aux 3-circuits de $\GG$ à droite (resp. à gauche) de $\Gamma$.
\end{defi}

\subsubsection{Système puits-source de $\A_{\GG}$}

\begin{defi}
Soit $\GG$ un écimaèdre étiqueté, un $\emph{système puits-source de
3-circuits de}$ $\GG$ est le choix d'une orientation de chaque
3-circuit de $\GG$, de telle sorte que si toutes les arêtes de
l'arbre $\A_{\GG}$ sont orientées dans le sens des 3-circuits
orientés choisis, alors tous les sommets de $\A_{\GG}$ sont des
puits ou des sources.
\end{defi}

\begin{rema}
Puisque le graphe $\A_{\GG}$ est un arbre, il est clair que tout écimaèdre étiqueté $\GG$ possède exactement deux systèmes puits-sources inverses l'un de l'autre.
\end{rema}

\subsection{Orientation partielle et globale de $\F_{\GG}$ induite par $P \in X_{\GG}$ via un système
puits-source}\label{orien}

\subsubsection{Orientation partielle induite}

\begin{defi}
Soient $\GG$ un écimaèdre étiqueté et $P$ un polyèdre miroir qui
réalise $\GG$. On suppose que l'on s'est donné un système
puits-source de 3-circuits de $\GG$. \emph{L'orientation partielle de $\F_{\GG}$ induite par $P$} est l'orientation obtenue en
orientant toute arête spéciale $e$ de $\F_{\GG}$ dans le sens de
$\Gamma$ (l'unique 3-circuit orienté correspondant à $e$ via
le système puits-source que l'on s'est donné), lorsque $R_{\Gamma}(P) > 0$ et dans
le sens contraire lorsque $R_{\Gamma}(P) < 0$.
\end{defi}

\begin{rema}
Il n'y a pas d'ambiguïté dans la définition précédente. En effet, si
$\Gamma$ est un 3-circuit spécial de $\GG$, alors $R_{\Gamma}(P)
\neq 0$.
\end{rema}

\subsubsection{Orientation globale induite}\label{orieninduite}

Soit $\GG$ un écimaèdre étiqueté. On se donne un système
puits-source de 3-circuits orientés de $\GG$. Supposons à présent que pour
tout 3-circuit sans angle droit $\Gamma$ de $\GG$, on ait
$R_{\Gamma}(P) \neq 0$. Alors on peut définir de façon analogue
une \emph{orientation globale de} $\F_{\GG}$ \emph{induite par} $P$
via le système puits-source que l'on s'est donné.

\subsubsection{Obstruction}

Il faut bien faire attention au fait que toutes les orientations
globales ou partielles de $\F_{\GG}$ ne peuvent pas être induites
par des polyèdres miroirs.

En effet, si on note, pour chaque sommet $s$ de $\F_{\GG}$ de
valence 4, $(\Gamma^s_i)_{i=1...4}$ les quatre 3-circuits orientés
via notre système puits-source, qui correspondent aux arêtes de
$\F_{\GG}$ incidentes à $s$, et si $P$ est un polyèdre miroir qui
réalise $\GG$, alors la définition même des $(R_{\Gamma^s_i}(P))_{i=1...4}$ (voir notations \ref{defirgamma}) montre qu'ils vérifient
la relation suivante :
$$
\begin{array}{cc}
R_{\Gamma^s_1}+R_{\Gamma^s_2}+R_{\Gamma^s_3}+R_{\Gamma^s_4}=0. & (*)\\
\end{array}
$$
Par conséquent, aucune orientation induite partielle ou globale ne
peut posséder de sommet de valence 4 qui soit un puits ou une
source, ce qui explique la définition d'orientation partielle ou
globale \emph{admissible}.

Nous allons montrer que les composantes connexes de $X_{\GG}$ sont
en bijection avec les orientations partielles admissibles de
$\F_{\GG}$.

\subsection{Le tétraèdre miroir et les blocs fondamentaux}\label{4}

\subsubsection{Le tétraèdre miroir}\label{paratetra}

La proposition suivante donne une paramétrisation de l'espace $X_{\mathcal{G}}$ lorsque $\GG$ est un tétraèdre combinatoire étiqueté. Cette proposition conclut la démonstration du point 2)a) du théorème \ref{theo} et est le point de départ pour l'étude des blocs fondamentaux.

\begin{prop}\label{tetra}
Soit $\mathcal{G}$ un tétraèdre combinatoire étiqueté marqué. Alors
$X_{\mathcal{G}}$ est difféomorphe à $\R^{d(\GG)}$ si $d(\GG)
\geqslant 0$ et est un singleton sinon.
\newline
\begin{itemize}
\item Cas 1: Si toutes les arêtes de $\GG$ sont d'ordre différent
de 2, on choisit un système puits-source de 3-circuits de $\GG$,
$\{\Gamma_1, \Gamma_2, \Gamma_3, \Gamma_4 \}$. L'application $P
\in X_{\GG} \mapsto (R_{\Gamma_1}, \, R_{\Gamma_2}, \,
R_{\Gamma_3}, \, R_{\Gamma_4}) \in \{ (r_1,r_2,r_3,r_4) \in \R^4
\textrm{ tels que } r_1 + r_2 + r_3 + r_4=0\}$ est un
difféomorphisme.
\newline
\item Cas 2: Si une seule arête de $\GG$ est d'ordre 2, alors $\GG$
possède exactement deux 3-circuits orientés (à orientation près)
sans angle droit $\Gamma_1$ et $\Gamma_2$. L'application $P \in
X_{\GG} \mapsto (R_{\Gamma_1},R_{\Gamma_2}) \in \R^2$ est un
difféomorphisme.
\newline
\item Cas 3: Si exactement 2 arêtes de $\GG$ sont d'ordre 2, et si
elles sont sur la même face, alors $\GG$ possède exactement un
3-circuit orienté (à orientation près) $\Gamma$ sans angle droit.
L'application $P \in X_{\GG} \mapsto R_{\Gamma} \in \R$ est un
difféomorphisme.
\newline
\item Cas 4: Si exactement 2 arêtes de $\GG$ sont d'ordre 2, et si
elles ne sont pas sur la même face, alors il existe un unique
4-circuit sans angle droit $\Gamma$ de $\GG$ (à orientation près).
L'application $P \in X_{\GG} \mapsto R_{\Gamma} \in \R$ est un
difféomorphisme.
\newline
\item Cas 5: Si au moins 3 arêtes de $\GG$ sont d'ordre 2, alors
$X_{\GG}$ est un singleton.
\end{itemize}
\end{prop}

\begin{proof}
On ne démontre la proposition que lorsque $\GG$ ne possède pas
d'arête d'ordre 2, les autres cas se démontrent de façon analogue.
Soit $P$ un tétraèdre qui réalise $X_{\GG}$. On procède comme pour
le triangle étiqueté (proposition \ref{modtri}). On note $(e_i)_{i=1...4}$ la base canonique de
$\R^4$. On peut supposer que $\alpha_i=-e_i^*$ pour $i=1...4$ et
$v_1=(-2,\sqrt{\mu_{12}},\sqrt{\mu_{13}},\sqrt{\mu_{14}})$. Le
stabilisateur de $\alpha_1,...,\alpha_4, v_1$ est réduit à
l'identité. On rappelle que les notations \ref{equa} fournissent les équations suivantes sur $v_2,v_3,v_4$:

$$
\begin{array}{cc}
\alpha_i(v_i) = 2 & \textrm{pour } i=1,...,4\\
\alpha_i(v_j)\alpha_j(v_i) = \mu_{ij} & \textrm{pour } i,j=1,...,4 \textrm{ et } i \neq j
\end{array}
$$

Ainsi, on obtient les coordonnées des vecteurs $v_2,v_3$ et $v_4$.

$$
\begin{array}{cccccc}
v_2= ( & \sqrt{\mu_{12}} & -2 & \sqrt{\mu_{23}}x & \frac{\sqrt{\mu_{24}}}{y} & )\\
v_3= ( & \sqrt{\mu_{13}} & \frac{\sqrt{\mu_{23}}}{x} & -2 & \sqrt{\mu_{34}}z & )\\
v_4= ( & \sqrt{\mu_{14}} & \sqrt{\mu_{24}}y & \frac{\sqrt{\mu_{34}}}{z} & -2 & )\\
\end{array}
$$

avec $x,y,z > 0$. On a alors comme dans la preuve de la proposition \ref{modtri}: $R_{\Gamma_1} = -2\ln(xyz)$,
$R_{\Gamma_2}=2\ln(z)$, $R_{\Gamma_3} = 2\ln(y)$ et $R_{\Gamma_4}
= 2\ln(x)$. Ce qui montre le résultat.
\end{proof}

\subsubsection{Les blocs fondamentaux}

\begin{prop}\label{blocfix}
Soit $\GG$ un bloc fondamental étiqueté qui n'est pas un prisme
exceptionnel et tel que $m(\GG)=0$. On se donne un système
puits-source de 3-circuits de $\GG$ et on note
$(\Gamma_i)_{i=1...4}$ les quatre 3-circuits de $\GG$ de notre
système puits-source qui viennent du tétraèdre $\mathcal{T}$
sous-jacent à $\GG$.

L'application $\varphi:P \in X_{\GG} \mapsto Q \in X_{\mathcal{T}}$, où $Q$ est le
tétraèdre miroir sous-jacent à $P$, est un difféomorphisme sur son
image $\mathcal{E} = \{ Q \in X_{\mathcal{T}} \, | \, \forall i=1,...,4 \,\, ,
\, R_{\Gamma_i}(Q) \in I_{\Gamma_i}^{\GG} \}$. En particulier,
$X_{\GG}$ est difféomorphe à $\kappa(\GG)$ copies de $\R^{d(\GG)}$ où
$\kappa(\GG)$ est le nombre d'orientations partielles admissibles de
$\F_{\GG}$.
\end{prop}

\begin{proof}
Soit $P$ un polyèdre miroir de $X_{\GG}$. Comme $\GG$ n'est pas un
prisme exceptionnel, le lemme \ref{ecimage2} montre que pour
tout 3-circuit prismatique $\Gamma$ de $\GG$, $\Pi_{\Gamma}$ est un
plan.

Commençons par montrer que l'application $\varphi:P \in X_{\GG} \mapsto Q
\in X_{\mathcal{T}}$ (le tétraèdre sous-jacent) est injective. On appellera faces tronquées les faces de $P$ qui ne sont pas des faces de $Q$. Les faces tronquées sont des faces triangulaires où toutes les arêtes sont d'ordre 2. La remarque \ref{uniqueref} montre que la réflexion par rapport à une face tronquée $s$ est entièrement déterminée par les réflexions par rapport aux faces adjacentes à $s$.

On rappelle que $X_{\mathcal{T}}$ est l'espace des modules des polyèdres
projectifs miroirs marqués qui réalisent le tetraèdre étiqueté $\mathcal{T}$.

Montrons ensuite que l'image de $\varphi$ est incluse dans $\mathcal{E} =\{ Q \in X_{\mathcal{T}} \, | \, \forall i=1,...,4 \,\, , \, R_{\Gamma_i}(Q) \in I_{\Gamma_i}^{\GG} \}$. C'est exactement le contenu de la proposition \ref{propoublie}.

Ensuite, il faut montrer que l'image de $\varphi$ est exactement $\mathcal{E}$. On se donne un tétraèdre $Q \in X_{\mathcal{T}}$ qui vérifie que $\forall i=1,...,4 \,\, , \, R_{\Gamma_i}(Q) \in I_{\Gamma_i}^{\GG}$, alors la proposition \ref{propoublie2} montre si $\G_i$ est prismatique pour $\GG$ et $R_{\Gamma_i}(Q) \in I_{\Gamma_i}^{\GG}$ alors le plan $\Pi_{\G_i}$ coupe $Q$ le long de $\G_i$. Enfin, le lemme \ref{nonchevau} montre que les nouvelles faces créées par ces écimages ne se chevauchent pas, par conséquent il existe un polyèdre miroir $P \in X_{\GG}$ tel que $\varphi(P) = Q$.

Enfin, il ne reste plus qu'à montrer que l'inverse de $\varphi$ est continue. Mais il s'agit simplement de voir que la réflexion par rapport à une face tronquée $s$ dépend continûment des réflexions par rapport aux faces adjacentes à $s$. Et cela est clair à la vue de la remarque \ref{uniqueref}.

Comme $m(\GG)=0$ tous les ensembles $I_{\Gamma_i}^{\GG}$ sont non vides (remarque \ref{pasvide}) et
ainsi la proposition \ref{tetra} montre que $\mathcal{E}$ est une
réunion d'ouverts convexes disjoints, naturellement indexée par
les orientations partielles admissibles de $\F_{\GG}$.
\end{proof}

Pour illustrer cette proposition on donne un tableau dont les entrées en colonnes indiquent le nombre de 3-circuits affines ou sphériques, sans angle droit et prismatiques de $\GG$. Les entrées en lignes indiquent de quel écimaèdre $\GG$ il s'agit (on rappelle que $\mathcal{T}_i$ est le tétraèdre écimé en $i$ sommets distincts). Les sorties sont les nombres $\kappa(\GG)$ correspondants.

$$
\begin{tabular}{|c|ccccc|}
\hline
 & $\mathcal{T}_0$ & $\mathcal{T}_1$ & $\mathcal{T}_2$ & $\mathcal{T}_3$ & $\mathcal{T}_4$\\
\hline
0 & 1 & 1 & 1 & 1 & 1\\
1 &  & 2 & 2 & 2 & 2\\
2 &  &  & 4 & 4 & 4\\
3 &  &  &  & 8 & 8\\
4 &  &  &  &  & 14\\
\hline
\end{tabular}
$$

\subsubsection{Les prismes exceptionnels}

Nous avons déjà traité le cas des prismes exceptionnels avec angles droits à l'aide du lemme \ref{tropde2} et du second point du théorème \ref{theo} (proposition \ref{tetra}).

\begin{prop}
Soit $\GG$ un prisme exceptionnel sans angle droit. On note
$\Gamma$ l'unique 3-circuit prismatique de $\GG$.
\begin{itemize}
\item Si $\Gamma$ est sphérique, alors $Card(X_{\GG})=2$.

\item Si $\Gamma$ est affine, alors $X_{\GG}$ est un singleton.

\item Si $\Gamma$ est hyperbolique, alors $X_{\GG} = \varnothing$.
\end{itemize}
\end{prop}

\begin{proof}
On numérote les faces de l'unique 3-circuit prismatique $\Gamma$
de $\GG$ de 1 à 3, et on numérote les deux autres faces 4 et 5. Les mêmes arguments que ceux de la démonstration du lemme \ref{tropde2} montrent que l'on peut supposer que l'on a la configuration suivante dans une base
$(e_i)_{i=1...4}$ de $\R^4$: $\alpha_i=-e_i^*$ pour $i=1...4$ et
$\alpha_5=-e_1^*-e_2^*-e_3^*+e_4^*$. Enfin, pour simplifier la
discussion, on suppose que $\mu_{31} \geqslant \mu_{12}$. On utilise les notations \ref{equa}. Les équations $\alpha_5(v_1) = \alpha_5(v_2) = \alpha_5(v_3) = 0$ donnent les équations suivantes:

$$
\left\{
\begin{array}{cccccc}
-2     &+& v_{12} &+&v_{13}  &= 0\\
v_{23} &-& 2      &+&v_{21}  &= 0\\
v_{31} &+& v_{32} &-& 2      &=0\\
\end{array}
\right.
$$

On obtient ainsi que $X_{\GG}$ est homéomorphe à:
$$
\left\{ (v_{12},v_{23},v_{31}) \in \R^3 \textrm{ tels que:
}
\begin{array}{c}
v_{12},v_{23},v_{31} > 0\\
v_{12}+\frac{\mu_{13}}{v_{31}} = 2\\
v_{23}+\frac{\mu_{12}}{v_{12}} = 2\\
v_{31}+\frac{\mu_{23}}{v_{23}} = 2
\end{array}
\right\}.
$$

Ainsi, si on pose $x=v_{12}$ et que l'on substitue les quantités $v_{31}$ et $v_{23}$ dans la dernière ligne du système précédent, alors on obtient que $X_{\GG}$ est en bijection avec les racines de
$f(x)=\frac{\mu_{31}}{2-x}+\frac{\mu_{23}}{2-\frac{\mu_{12}}{x}}-2$
qui appartiennent à l'intervalle $]\frac{\mu_{12}}{2},2[$.

On pose $\sigma =
-\frac{\mu_{23}+\mu_{31}-\mu_{12}-4}{4-\mu_{23}}$ et $p =
\frac{\mu_{12}(4-\mu_{31})}{4-\mu_{23}}$. On est ainsi ramené à
l'étude du polynôme $Q=x^2-2\sigma x +p$ sur
$]\frac{\mu_{12}}{2},2[$. On pose alors $Y=x-\sigma$,
$\delta=\sigma^2-p$, $u=\frac{\mu_{12}}{2}-\sigma$ et
$v=2-\sigma$ (noter $u < v$). On est ainsi ramené à l'étude du polynôme
$R=Y^2-\delta$ sur $]u,v[$.

Le discriminant $\delta$ de $R$ est égal à
$$\delta =
\frac{(4-\mu_{31}-\mu_{12}-\mu_{23}-\sqrt{\mu_{31}\mu_{12}\mu_{23}})(4-\mu_{31}-\mu_{12}-\mu_{23}+\sqrt{\mu_{31}\mu_{12}\mu_{23}})}{(4-\mu_{23})^2}.$$
Comme dans la preuve du lemme \ref{triangle}, on a:
$$
\begin{array}{ll}
4-\mu_{31}-\mu_{12}-\mu_{23} & -\sqrt{\mu_{31}\mu_{12}\mu_{23}}
=\\
 &
 \textstyle{-16\cos\Big(\frac{\theta_{12}+\theta_{23}+\theta_{31}}{2}\Big)\cos\Big(\frac{\theta_{12}-\theta_{23}+\theta_{31}}{2}\Big)\cos\Big(\frac{\theta_{12}+\theta_{23}-\theta_{31}}{2}\Big)\cos\Big(\frac{-\theta_{12}+\theta_{23}+\theta_{31}}{2}\Big)}.
\end{array}
$$
On remarque que, comme les angles sont aigus, le produit des 3
derniers membres est strictement positif.

On rappelle que par définition $\mu_{ij} = 4\cos^2(\theta_{ij})$, où $0 < \theta_{ij} < \frac{\pi}{2}$ car on a supposé que $m(\GG) = 0$ (i.e. $\Gamma$ sans angle droit). On a les inégalités suivantes:

$$
\left\{
\begin{array}{cc cc cc ccc}
v      & = & \frac{4-\mu_{23}+\mu_{31}-\mu_{12}}{4-\mu_{23}}                             & \geqslant 1 & \\
u^2-\delta & = &
\frac{\mu_{12}\mu_{23}(4-\mu_{12})}{4(4-\mu_{23})}
& > 0  & &&(S_1)\\
v^2-\delta & = & \frac{\mu_{31}(4-\mu_{12})}{4-\mu_{23}}& > 0 & \\
\frac{\mu_{12}\mu_{23}}{2} & <  & \sqrt{\mu_{31}\mu_{12}\mu_{23}} &\\
\end{array} \right.
$$

Ces quatres inégalités viennent de l'hypothèse $\mu_{12} \leqslant \mu_{31}$ et des inégalités $0 < \mu_{12}, \mu_{23},\mu_{31} < 4$.

Le signe de la quantité $D=4-\mu_{31}-\mu_{12}-\mu_{23} -\sqrt{\mu_{31}\mu_{12}\mu_{23}}$ est crucial dans cette preuve.

Pour conclure, il faut distinguer 2 cas:

1) $\Gamma$ est sphérique (i.e. $D > 0$) ou affine (i.e.
$D = 0$). On calcule

$$
\left\{
\begin{array}{cc cc cc}
\delta  & \geqslant & 0
                                                                                          &  \\
-u      & = &  \frac{4-\mu_{31}-\mu_{12}-\mu_{23}+\frac{\mu_{12}\mu_{23}}{2}}{4-\mu_{23}} & \geqslant 0\\
\end{array}
\right.
$$

Ces deux inégalités sont évidentes puisqu'on a supposé que $D \geqslant 0$.

Donc si $\Gamma$ est sphérique, alors $u < -\sqrt{\delta} < \sqrt{\delta} < v$ donc $R$ possède 2 racines sur
l'intervalle $]u,v[$, et si $\Gamma$ est affine, alors $R$ possède une racine double sur l'intervalle $]u,v[$.

2) Si $\Gamma$ est hyperbolique (i.e. $D < 0$), alors il faut de nouveau distinguer 2 cas, posons $\overline{D} = 4-\mu_{31}-\mu_{12}-\mu_{23}+\sqrt{\mu_{31}\mu_{12}\mu_{23}}$:

a) Si $\overline{D} > 0$ alors $\delta <0$ et les racines de $R$ sont complexes.

b) Si $\overline{D} \leqslant 0$ alors $-u(4-\mu_{23}) < \overline{D} \leqslant 0$, puisque la quatrième inégalité du système $S_1$ entraîne l'inégalité suivante: $4-\mu_{31}-\mu_{12}-\mu_{23}+\frac{\mu_{12}\mu_{23}}{2} < \overline{D} \leqslant 0$.

On a $\delta>0$ mais $\sqrt{\delta} < u < v$: les racines de $R$ sont réelles mais à l'extérieur de $]u,v[$.
\end{proof}

\subsubsection{Blocs fondamentaux à structure projective fixée}

Dans ce paragraphe, on cherche à comprendre l'espace des modules
d'un bloc fondamental $\GG$ dont la structure projective d'un
triangle est fixée.

\begin{nota}
Soit $\GG$ un bloc fondamental étiqueté marqué, on suppose que l'on
s'est donné un système puits-source de 3-circuits de $\GG$. Soit
$\Gamma$ un 3-circuit prismatique orienté de $\GG$. On se donne $r
\in I_{\Gamma}^{\GG}$ et on pose $X_{\GG}^r= \{ P \in X_{\GG} |
R_{\Gamma} = r \}$.
\end{nota}

On dira qu'une orientation de $\F_{\GG}$ est $\emph{compatible avec
r}$ lorsque $r$ appartient à la composante connexe de
$I_{\Gamma}^{\GG}$ donnée par l'orientation de $\F_{\GG}$.

\begin{prop}\label{coro3}
Soit $\GG$ un bloc fondamental étiqueté qui n'est pas un prisme
exceptionnel et tel que $m(\GG)=0$, muni d'un système puits-source.
Soit $\Gamma$ un 3-circuit prismatique de $\GG$ orienté via le
système puits-source, et $r \in I_{\Gamma}^{\GG}$.

Alors, $X_{\GG}^r$ est difféomorphe à $\kappa(\GG,\Gamma,r)$ copies
de $\R^{d(\GG)-1}$ où $\kappa(\GG,\Gamma,r)$ est le nombre
d'orientations partielles admissibles de $\F_{\GG}$ compatibles avec $r$.
\end{prop}

\begin{proof}
L'idée de la démonstration est que les propositions \ref{tetra} et \ref{blocfix} fournissent une paramétrisation de l'espace $X_{\GG}$ uniquement à l'aide des quantités $(R_{\G_i})_{i=1,...,4}$ où les $(\G_i)_{i=1,...,4}$ sont les 3-circuits du tétraèdre étiqueté $\mathcal{T}$ sous-jacent à $\GG$ de notre système puits-source.

Le moyen le plus simple de démontrer cette proposition est de faire une étude exhaustive. Pour restreindre le nombre de cas, on peut remarquer que comme $m(\GG) = 0$ et que $\GG$ n'est pas un prisme exceptionnel, le tétraèdre étiqueté $\mathcal{T}$ possède au plus 3 arêtes d'ordre 2.

On obtient facilement le résultat, en distinguant les cas: $\mathcal{T}$ ne possède aucune arête d'ordre 2, une seule arête d'ordre 2, deux arêtes opposées d'ordre 2, deux arêtes non opposées d'ordre 2 ou exactement 3 arêtes d'ordre 2.
\end{proof}

\begin{rema}\label{nombrefix}
Pour illustrer cette proposition on donne un tableau dont les entrées en colonnes indiquent le nombre de 3-circuits affines ou sphériques, sans angle droit et prismatiques de $\GG$ différents de $\G$. Les entrées en lignes indiquent de quel écimaèdre $\GG$ il s'agit. Les sorties sont les nombres $\kappa(\GG,r)$ correspondants.

$$
\begin{tabular}{|c|cccc|}
\hline
 & $\mathcal{T}_1$ & $\mathcal{T}_2$ & $\mathcal{T}_3$ & $\mathcal{T}_4$\\
\hline
0 & 1 & 1  & 1 & 1 \\
1 &   & 2  & 2 & 2  \\
2 &   &    & 4 & 4   \\
3 &   &    &   & 7    \\
\hline
\end{tabular}
$$

En particulier, $\kappa(\GG,\Gamma,r)$ ne dépend pas de $r$, on le notera donc $\kappa(\GG,\Gamma)$.
\end{rema}

\subsection{Lemme de recollement}\label{3}

On cherche à présent à comprendre comment coller 2 polyèdres miroir, $L$ et $D$, le long d'une face triangulaire $T_L$ incluse dans $L$ et d'une face triangulaire $T_D$ incluse dans $D$ (les arêtes de $T_D$ et $T_L$ sont d'ordre 2). Notons $\G_L$ (resp. $\G_D$) le 3-circuit prismatique de $L$ (resp. $D$) qui fait le tour de la face triangulaire $T_L$ (resp. $T_D$).

Pour recoller $L$ et $D$ le long de $T_L$ et $T_D$, une condition nécessaire est que les quantités $R_{\G_L}(L)$ et $R_{\G_D}(D)$ soient égales. Nous allons voir avec le lemme \ref{colle} que c'est une aussi une condition suffisante. Le lemme \ref{colle} montre aussi que ce recollement n'est pas unique, on peut le paramétrer par un nombre réel. Pour paramétrer ce recollement, on va utiliser l'invariant $R$ d'un triplet de trois faces (pas nécessairement adjacentes) dont l'une appartient seulement à $L$, l'autre à $D$ et la dernière aux deux.

\begin{nota} Comme il y a plusieurs polyèdres en jeu dans le lemme qui
suit, on notera $R_{\Gamma}(P)$ la quantité que l'on notait
habituellement $R_{\Gamma}$ associée à un 3-circuit orienté
$\Gamma$ de $P$.
\end{nota}

\begin{lemm}\label{colle}
Soit $\GG$ un graphe étiqueté tel que $m(\GG)=0$ et $d(\GG) \geqslant
0$. On suppose que $\GG$ possède un 3-circuit orienté prismatique
essentiel $\Gamma$. On note $\GG^d_{\Gamma}$ et $\GG^g_{\Gamma}$ les
coupes droite et gauche de $\GG$ par rapport à $\Gamma$.

On suppose que $\GG^d_{\Gamma}$ est un bloc fondamental étiqueté.

Soient $s$ une face traversée par $\Gamma$, $l$ une face de $\GG$ à
gauche de $\Gamma$, et enfin $d$ une face de
$\GG$ à droite de $\Gamma$. On suppose que si $l$ (resp. $d$) et $s$ sont adjacentes alors l'arête commune n'est pas d'ordre 2. On suppose aussi que ni la face $l$, ni la face $d$ ne font partie de $\G$.

Alors l'application $\phi : P \in X_{\GG} \mapsto$
$(P^g_{\Gamma},R_{(l,s,d)}) \in X_{\GG^g_{\Gamma}} \times \R$ est
une fibration dont les fibres $\phi^{-1}(L,r)$ s'identifient à
$E_{L,r}=\{ D \in X_{\GG^d_{\Gamma}} \, | \, R_{\Gamma}(D) =
R_{\Gamma}(L) \}=X_{\GG^d_{\Gamma}}^{R_{\Gamma}(L)}$.

\end{lemm}

\begin{rema}
On rappelle que l'on peut définir la quantité $R_{(l,s,d)}$ pour n'importe quel triplet ordonnés de faces (notations \ref{defirgamma}). Il est essentiel de remarquer que les faces $l,s,d$ ne forment pas un 3-circuit.
\end{rema}

\begin{rema}
Pour mesurer le paramètre de recollement on a besoin de trois faces $l,s,d$. Pour que la quantité $R_{(l,s,d)}$ soit définie, on impose que si $l$ (resp. $d$) et $s$ sont adjacentes alors l'arête commune n'est pas d'ordre 2. Pour mesurer effectivement le recollement, on impose qu'aucune des faces $l$ et $d$ ne fait partie de $\G$. Comme $m(\G)=0$ et que $\G$ n'est pas un prisme exceptionnel, on peut toujours trouver un tel triplet.
\end{rema}

\begin{rema}
L'ensemble $E_{L,r}$ est difféomorphe à
$\kappa(\GG^d_{\Gamma},\Gamma)$ copies de
$X^r_{\GG^d_{\Gamma}} \simeq \R^{d(\GG^d_{\Gamma})-1}$, d'après la proposition \ref{coro3}. On verra qu'il ne dépend pas de $r$.
\end{rema}

\begin{proof}

On remarque pour commencer que $\GG^d$ ne peut être un prisme
exceptionnel car $\G$ est essentiel. Il faut aussi noter que $\GG^d$ est un bloc fondamental, pas un écimaèdre quelconque, donc on connaît la topologie de $X_{\GG}$. Soit $F$ le sous-fibré du fibré $X_{\GG^g} \times \R \times X_{\GG^d}$ de base
$X_{\GG^g} \times \R$ dont la fibre au-dessus du point $(P^g,r)$ est l'espace $\{ P^d \in X_{\GG^d}
\, | \, R_{\Gamma}(P^g) = R_{\Gamma}(P^d) \}$ (la proposition
\ref{coro3} et la remarque \ref{nombrefix} montrent que les fibres sont difféomorphes deux à deux). Le
lemme \ref{coupris}  fournit une application naturelle $\psi:P
\in X_{\mathcal{G}} \mapsto (P^g , R_{ ( l,s,d ) } , P^d) \in F
\subset X_{\GG^g} \times \R \times X_{\GG^d}$ et on a bien
$R_{\Gamma}(P^g) = R_{\Gamma}(P^d)$, donc $\psi$ est bien définie. On note $\pi$ la projection $\pi:F \rightarrow X_{\GG^g} \times \R$ donnée par l'oubli de la troisième coordonnée, on a $\phi = \pi \circ \psi$. Nous allons montrer que $\psi$ est un difféomorphisme, ce qui montrera la proposition.

On peut construire l'application réciproque de $\psi$ comme suit.

Soit $(L,r,D) \in F$, nous allons construire un polyèdre $P$ tel que $\psi(P) = (L,r,D)$. On commence par choisir une normalisation pour les faces de $L$ et $D$. On peut numéroter les faces traversées par $\Gamma$ de 1 à 3. Si on note $(e_i)_{i=1...4}$ une base de $\R^4$, alors on peut supposer que les faces traversées par
$\Gamma$ sont définies par les formes linéaires $-e_i^*$ pour
$i=1...3$, que la face $s$ est aussi la face 1, et que le plan
$\Pi_{\Gamma}$ est défini par la forme linéaire $-e_4^*$.
Les mêmes arguments que ceux de la démonstration de la proposition \ref{tropde2} montrent qu'on peut
aussi supposer que la face $l$ de $L$ est définie par
$-e_1^*-e_2^*-e_3^*-e_4^*$, que la face $d$ de $D$ est définie par
$-\lambda_1 e_1^*-\lambda_2 e_2^*- \lambda_3 e_3^*+ \lambda_4 e_4^*$, avec $\lambda_i>0$ pour $i=1,...,4$.

Les réflexions par rapport aux faces 1 (resp. 2, resp. 3) des polyèdres miroirs $L$ et $D$ sont les
mêmes, car $R_{\Gamma}(L) = R_{\Gamma}(D)$. On appellera face 4 de $L$ et $D$, la face définie par le plan $\Pi_{\G}$.
Il faut remarquer que $\alpha_4$ définit la face 4 du polyèdre $D$, alors que c'est $-\alpha_4$ qui définit la face 4 du polyèdre $L$.

On peut alors construire le polyèdre miroir $P$ dont les faces
et les réflexions sont les faces et les réflexions de $L$ ou $D$
moins la face de $L$ entourée par $\Gamma$ et la face de $D$
entourée par $\Gamma$. Mais on peut en construire d'autres. En
effet, pour tout $\lambda > 0$, considérons l'élément
$g_{\lambda}$ de $\sss_4(\R)$ dont la matrice dans la base
$(e_1,e_2,e_3,e_4)$ est:

$$
g_{\lambda} = \left( \begin{array}{c c c c} \lambda & & & \\
 & \lambda & & \\
 & & \lambda & \\
 & & & \lambda^{-3}
\end{array} \right).
$$

L'hyperplan propre de $\sigma_i$ est le noyau de $\alpha_i = -e_i^*$, pour $i=1,...,3$. Le vecteur propre $v_i$ de $\sigma_i$ appartient au noyau de la forme linéaire $-e_4^*$, pour $i=1,...,3$. L'élément $g_{\lambda}$ commute donc avec les réflexions $\sigma_1,\sigma_2,\sigma_3$. De plus, tout élément $g \in \sss_4(\R)$ qui commute avec $\sigma_1$, $\sigma_2$, $\sigma_3$ et qui fixe point par point dans $\mathbb{P}^+(\R^4)$ la face 4 de $L$ et $D$, est un $g_{\lambda}$ avec $\lambda >0$.


On peut construire un polyèdre miroir $P_{\lambda} \in X_{\mathcal{G}}$ dont les
réflexions à gauche de $\Gamma$ sont les réflexions de $L$ et les
réflexions à droite de $\Gamma$ sont les conjuguées des réflexions de $D$ par $g_{\lambda}$. Comme $g_{\lambda}$ commute
avec $\sigma_1,\sigma_2,\sigma_3$, cette définition n'est pas
ambiguë. Et il est clair que les polyèdres $P_{\lambda}$ sont
convexes. De plus, si on a un polyèdre $P'$ qui vérifie $P'^g_{\G} = L$ et $P'^d_{\G}=D$, alors il existe un unique $\lambda > 0$ tel que $P'$ et $P_{\lambda}$ soient équivalents. Nous allons relier les quantités $r=R_{l,s,d}$ et $\lambda$.

Résumons les coordonnées des différentes quantités, on met le symbole $\lambda$ en exposant pour signaler que ces quantités sont attachées au polyèdre $P_{\lambda}$ (p.ex: $v_i^{\lambda}$ et $\alpha_i^{\lambda}$). Pour faciliter les calculs on note $u_l$ (resp. $w_l$) la projection de $v_l$ sur l'hyperplan engendré par $e_1,e_2,e_3$ (resp. la droite engendrée par $e_4$). On note $u_d$ (resp. $w_d$) la projection de $v_d$ sur l'hyperplan engendré par $e_1,e_2,e_3$ (resp. la droite engendrée par $e_4$). Enfin, il faut se rappeler que la face $s$ et la face 1 sont identiques.


$$
\begin{array}{cc}
\left\{
\begin{array}{cc}
\alpha_l^{\lambda}= & (\alpha_1+  \alpha_2 + \alpha_3)+ \alpha_4\\
\alpha_1^{\lambda}= & \alpha_1\\
\alpha_d^{\lambda}= & \lambda^{-1}(\lambda_1 \alpha_1+ \lambda_2 \alpha_2 +\lambda_3 \alpha_3)- \lambda^{3} \lambda_4 \alpha_4
\end{array}
\right.

&

\left\{
\begin{array}{cccccccc}
v_l^{\lambda}= &  u_l + w_l\\
v_1^{\lambda}= & v_1\\
v_d^{\lambda}= & \lambda u_d + \lambda^{-3} w_d\\
\end{array}
\right.
\end{array}
$$

Il ne  reste plus qu'à montrer qu'il existe un unique réel
$\lambda > 0$ tel que $R_{(l,s,d)}(P_{\lambda})=r$, et cela résulte
du calcul suivant:
$$
\exp\Bigg(\frac{R_{l,s,d}(P_{\lambda})}{2}\Bigg)=\frac{\alpha^{\lambda}_l(v^{\lambda}_1) \alpha^{\lambda}_1(v^{\lambda}_d)
\alpha^{\lambda}_d(v^{\lambda}_l)}{\alpha^{\lambda}_l(v^{\lambda}_d) \alpha^{\lambda}_d(v^{\lambda}_1) \alpha^{\lambda}_1(v^{\lambda}_l)} =
K \frac{a+\lambda^4 b}{c+\lambda^{-4} d},
$$
Pour cela on fait le calcul suivant:
$$
\left\{
\begin{array}{ccc ccc}
\alpha^{\lambda}_l(v^{\lambda}_1) = &  \alpha_l(v_1) & < 0\\

\alpha^{\lambda}_1(v^{\lambda}_d) =  &\lambda\alpha_1(v_d) & < 0\\

\alpha^{\lambda}_d(v^{\lambda}_l) = &  \lambda^{-1}(\lambda_1 \alpha_1+ \lambda_2 \alpha_2 +\lambda_3 \alpha_3)(u_l)- \lambda^{3} \lambda_4 \alpha_4(w_l)\\

\alpha^{\lambda}_l(v^{\lambda}_d) = & \lambda (\alpha_1+  \alpha_2 + \alpha_3)(u_d)+ \lambda^{-3}\alpha_4(w_d)\\

\alpha^{\lambda}_d(v^{\lambda}_1) =  & \lambda^{-1} \alpha_d(v_1) &< 0\\
\alpha^{\lambda}_1(v^{\lambda}_l) =  & \alpha_1(v_l) &< 0\\

\end{array}
\right.
$$

Avant de définir $K,a,b,c,d$, justifions ces calculs. Pour les égalités, il s'agit simplement de calculer, en se rappelant que $s=1$ et $\alpha_4(v_1)=0$ car la face $4$ est la face définissant le plan $\Pi_{\G}$. On a aussi par définition $\alpha_4(u_l)  = \alpha_4(u_d) = 0$ et $\alpha_i(w_l) = \alpha_i(w_d)=0$, pour $i=1,...,3$. Pour les inégalités, elles résultent simplement du lemme \ref{lemme1}. \`{A} présent, on pose:

$$
\left\{
\begin{array}{clll}
K & = & \frac{\alpha_l(v_1)\alpha_1(v_d)}{\alpha_d(v_1)\alpha_1(v_l)}  & > 0\\
a & = & -(\lambda_1 \alpha_1+ \lambda_2 \alpha_2 +\lambda_3 \alpha_3)(u_l) & >0\\
b & = & \lambda_4 \alpha_4(w_l)                   & >0\\
c & = & -(\alpha_1+  \alpha_2 + \alpha_3)(u_d) & >0\\
d & = & -\alpha_4(w_d)                       & >0.\\
\end{array}
\right.
$$

Il faut vérifier que les quantités $K,a,b,c,d$ sont strictement positives. Pour $K$, c'est une conséquence des inégalité précédentes. Ensuite, il faut se rappeler que $\alpha_4$ définit la face 4 du polyèdre $P^d$, alors que c'est $-\alpha_4$ qui définit la face 4 du polyèdre $P^g$. Par conséquent, $\alpha_4(w_d) = \alpha_4(v_d) < 0$ et $\alpha_4(w_l)  = \alpha_4(v_l) >0$. Il reste à montrer que $a$ et $c$ sont strictement positifs. Mais le lemme \ref{lemme1} montre que $\alpha_i(u_l) = \alpha_i(v_l) < 0$ et $\alpha_i(u_d) = \alpha_i(v_d) < 0$ pour $i=1,...,3$.

Par conséquent, la fonction $\lambda > 0 \mapsto K \frac{a+\lambda^4 b}{c+\lambda^{-4} d}$ est continue, strictement croissante et tend vers 0 en 0, et vers $+\infty$ en $+\infty$. Il existe un unique réel $\lambda > 0$ tel que $R_{(l,s,d)}(P_{\lambda})=r$. De plus, il est clair que $\lambda$ dépend continûment de $r$.

On définit ainsi $\zeta:(L,r,D) \in F \mapsto P_{\lambda} \in
X_{\GG}$, où $P_{\lambda}$ est l'unique polyèdre projectif miroir
tel que $P_{\lambda}^d = L$, $P_{\lambda}^g = D$ et
$R_{l,s,d}(P_{\lambda})=r$. La fonction $\zeta$ est continue et il est clair que les
applications $\zeta$ et $\psi$ sont inverses l'une de l'autre. Ce
qui termine la démonstration.

\end{proof}

\begin{rema}
Ce lemme  montre que $X_{\GG}$ est une variété et qu'elle est
de dimension $d(\GG)$. Il s'agit d'une récurrence sur le nombre de
3-circuits prismatiques essentiels de $\GG$. Le cas initial est donné par la proposition \ref{blocfix}, et sinon on choisit un 3-circuit orienté prismatique essentiel $\Gamma$ tel que $\GG^d_{\Gamma}$ soit un bloc
fondamental. Le lemme \ref{colle}, la proposition \ref{coro3} et l'hypothèse de récurrence
montrent que $X_{\GG}$ est une variété de dimension
$d(\GG^g_{\Gamma}) + (d(\GG^d_{\Gamma}) -1 ) + 1 = d(\GG)$.
\end{rema}

\subsection{Explicitation du difféomorphisme}\label{last}

On peut construire explicitement le difféomorphisme entre $X_{\GG}$
et les copies de $\R^{d(\GG)}$. Pour cela, introduisons les
notations suivantes:

\begin{enumerate}
\item On commence par se donner un système puits-source
$\mathcal{S}$ de 3-circuits de $\GG$.

\item On note $V$ l'ensemble des sommets de valence 4 de la forêt
$\F_{\GG}$. Si $w \in V$, alors on note $(\Gamma_i^w)_{i=1...4}$
une suite des quatre 3-circuits de $\mathcal{S}$ qui correspondent
aux arêtes incidentes de $w$.

\item On note $C_3$ la partie de $\mathcal{S}$ formée des
3-circuits sans angle droit de $\GG$, et on pose

$$
\mathcal{E}= \left\{
\begin{array}{c|lc}
&
\forall \, \Gamma \in C_3,\, &  x_{\Gamma} \in I_{\Gamma}^{\GG}\\
 x \in \R^{C_3}  & et &\\

& \forall \, w \in V, \, & \sum_{i=1...4} x_{\Gamma_i^w} = 0
\end{array}
\right\}.
$$

Cet espace $\mathcal{E}$ est une réunion d'ouverts convexes disjoints indexée
naturellement par les orientations partielles admissibles de $\F_{\GG}$.

\item On note $C'_3$ la partie de $\mathcal{S}$ formée des
3-circuits prismatiques essentiels de $\GG$. Pour tout élément
$\Gamma$ de $C'_3$, on choisit une face $s_{\Gamma}$ de $\Gamma$,
une face $s_{\Gamma}^d$ à droite de $\Gamma$, et enfin une face $s_{\Gamma}^g$ à gauche de
$\Gamma$. On suppose que si $s_{\Gamma}^g$ (resp. $s_{\Gamma}^d$) et $s_{\Gamma}$ sont adjacentes alors l'arête commune n'est pas d'ordre 2, et qu'aucune des faces $s_{\Gamma}^g$ et $s_{\Gamma}^d$ ne fait partie de $\G$.

\item On note $C_4$ l'ensemble des sommets de $\A_{\GG}$ dont le
bloc fondamental correspondant vient d'un tétraèdre $T$ qui
possède exactement 2 arêtes d'ordre 2 qui ne sont pas sur la même
face de $T$. \`{A} tout élément $w \in C_4$ est associé un unique
4-circuit orienté (à orientation près) noté $\Delta_w$.

\item On note $\mathcal{M}_{\GG}$ l'espace $\mathcal{E} \times
\R^{C_4} \times \R^{C'_3}$.
\end{enumerate}

\begin{rema}\label{finfin}
Soit $\GG$ un écimaèdre étiqueté. Si $w$ est un sommet de
$\A_{\GG}$, et si on note $\mathcal{B}_w$ le bloc fondamental
associé à $w$ alors on a une projection naturelle
$p_w:\mathcal{M}_{\GG}\rightarrow \mathcal{M}_{\mathcal{B}_w}$.Le lemme \ref{colle} montre que cette application est surjective.
\end{rema}

On a un énoncé plus précis que le théorème \ref{theo}.

\begin{theo}\label{final}
Soit $\GG$ un écimaèdre étiqueté qui n'est pas un prisme
exceptionnel et tel que $m(\GG)=0$ et $d(\GG) \geqslant 0$, avec les
notations introduites ci-dessus. L'application
$$
\begin{array}{ccccccccccc}
\phi_{\GG}: & X_{\GG} &  \longrightarrow & \mathcal{E} & \times & \R^{C_4} &
\times
& \R^{C'_3}\\

\\

 & P & \longmapsto & \bigg( (R_{\Gamma})_{\Gamma \in C_3} & , & (R_{\Delta_w})_{w \in
C_4} & , & (R_{( s_{\Gamma}^g, s_{\Gamma}, s_{\Gamma}^d)})_{\Gamma \in C'_3} \bigg)
\end{array}
$$

\noindent est un difféomorphisme sur $\mathcal{M}_{\GG}$. En particulier, les
composantes connexes de $X_{\GG}$ sont en bijection avec les
orientations partielles admissibles de $\F_{\GG}$, et $X_{\GG}$
s'identifie à une réunion d'ouverts convexes d'un espace vectoriel
de dimension $d(\GG)$.
\end{theo}

\begin{proof}
Ce théorème a déjà été démontré lorsque $\GG$ est un bloc
fondamental ou un tétraèdre: il s'agit des propositions
\ref{tetra} et \ref{blocfix}. Nous allons procéder par étapes.

Commençons par montrer que $\phi_{\GG}$ est injective. Pour cela,
il suffit de procéder par récurrence sur le nombre de 3-circuits
prismatiques essentiels de $\GG$. Si $\GG$ ne possède pas de tel
3-circuit, alors $\GG$ est un bloc fondamental ou un tétraèdre et
le théorème est démontré dans ce cas. Si $\GG$ possède un 3-circuit
prismatique essentiel alors $\GG$ possède un 3-circuit orienté
prismatique essentiel  $\Gamma$ tel que $\GG^d_{\Gamma}$ soit un
bloc fondamental. Le lemme de collage (lemme \ref{colle}),
l'hypothèse de récurrence et la proposition \ref{blocfix}
montrent que tout polyèdre $P \in X_{\GG}$ est entièrement déterminé
par $\phi_{\GG}(P)$. Donc $\phi_{\GG}$ est injective.

Montrons à présent que l'image de $\phi_{\GG}$ contient
$\mathcal{M}_{\GG}$. Soit $m \in \mathcal{M}_{\GG}$. Pour tout
sommet $w$ de $\A_{\GG}$, on note $\mathcal{B}_w$ le bloc
fondamental étiqueté de $\GG$ correspondant au sommet $w$. Le lemme
\ref{colle}  montre qu'il existe un $P \in X_{\GG}$ tel que
$\phi_{\GG}(P) = m$ si et seulement si pour tout sommet $w$ de
$\A_{\GG}$ il existe un bloc fondamental $B_w \in
X_{\mathcal{B}_w}$ tel que $\phi_{\mathcal{B}_w}(B_w)= p_w(m)$, où
$p_w$ désigne la projection naturelle $p_w:\mathcal{M}_{\GG}
\rightarrow \mathcal{M}_{\mathcal{B}_w}$ (remarque \ref{finfin}). La proposition
\ref{blocfix} et le fait que toute orientation partielle admissible
de $\F_{\GG}$ induit une orientation partielle admissible de
$\F_{\mathcal{B}_w}$  montrent l'existence de tels blocs
fondamentaux. La surjectivité de $p_w$ permet de conclure que
$\mathcal{M}_{\GG}$ est inclus dans l'image de $\phi_{\GG}$.

Il faut montrer à présent que l'image de $\phi_{\GG}$ est incluse
dans $\mathcal{M}_{\GG}$. Soit $P \in X_{\GG}$, montrons que $P$
définit une orientation partielle admissible.
La proposition \ref{propoublie} montre que pour tout 3-circuit sans angle droit $\Gamma$, on a $R_{\Gamma}(P) \in I_{\G}$. Par conséquent, si le 3-circuit $\G$ est prismatique sans angle droit et affine ou sphérique alors en particulier $R_{\Gamma}(P) \neq 0$. On a vu au paragraphe \ref{orieninduite} qu'un tel polyèdre définit une orientation partielle de $\F_{\GG}$. Il ne nous reste plus qu'à montrer que cette orientation est admissible. La difficulté est la suivante, a priori le polyèdre $P$ peut posséder un 3-circuit $\G$ prismatique et hyperbolique qui vérifie $R_{\G}(P)=0$. La présence d'un tel 3-circuit empêche $P$ de définir une orientation partielle admissible à l'aide des signes des invariants $R$. Pour remédier à ce problème, il suffit de remarquer que l'espace des modules de tout bloc fondamental $\mathcal{B}$ est une variété connexe $X_{\mathcal{B}}$ et que l'ensemble des bloc fondamentaux miroirs $B$ qui possède un 3-circuit $\G$ qui vérifie $R_{\G}(B)=0$ est une réunion de sous-variété (de codimension 1) de cette variété $X_{\mathcal{B}}$. Par conséquent, pour tout polyèdre $P$, on peut trouver une famille continue $(P_t)_{t \in [0,1]} \in X_{\GG}^{[0,1]}$ telle que $P_0=P$ et pour tout $t \neq 0$ et tout 3-circuit $\G$ de $\GG$, on a $R_{\G}(P_t) \neq 0$. En effet, il est clair que l'on peut faire cela pour tout bloc fondamental et comme $\A_{\GG}$ est un arbre on peut s'assurer que le long de ces déformations les conditions de recollement (égalités des invariants $R$, lemme \ref{colle}) soient assurées.


Il faut à présent montrer la continuité de l'application réciproque. On a déjà vu avec la proposition \ref{tetra} et la proposition \ref{blocfix} que les réflexions par rapport aux faces des blocs fondamentaux dépendaient continûment des quantités $R_{\G}$ pour $\G \in C_3$. Il ne reste plus qu'à vérifier que les réflexions dépendent continûment du paramètre de recollement. On a vu durant la démonstration du lemme \ref{colle} que les réflexions dépendaient continûment de ce paramètre. Ce qui conclut la démonstration.
\end{proof}

\begin{rema}
Il est à présent clair que $\kappa(\GG)$ est un entier pair ou égal
à 1, qui vérifie $1 \leqslant \kappa(\GG) \leqslant 2^{n(\GG)}$. En effet, si on possède une orientation partielle admissible alors l'orientation obtenue en renversant toutes les flèches est aussi admissible.
\end{rema}

\section{Exemples}

\par{
On donne ici quelques exemples de calcul du nombre de composantes
connexes de $X_{\GG}$. Les figures sont décomposées en trois
parties: le graphe combinatoire associé à $\GG$ et l'arbre associé,
que l'on a représenté deux fois.
}
\\
\par{
Pour rendre les figures lisibles, nos exemples vérifient tous les
hypothèses suivantes: les arêtes extrémales sont sphériques avec
angle droit et les arêtes non extrémales sont sans angle droit.
}
\\
\par{
Enfin, on a inscrit sur les arêtes la nature géométrique de
l'arête, en utilisant les abréviations suivantes: a = affine ou sphérique et h
= hyperbolique.
}
\\
\par{
On désignera par $\kappa_1$ (resp. $\kappa_2$) le nombre de
composantes connexes de l'espace des modules associé à n'importe
quel graphe étiqueté qui vient du graphe étiqueté dessiné et dont
les 3-circuits vérifient les hypothèses données par les étiquettes des
arêtes de l'arbre dessiné à gauche (resp. droite).
}
\\
\par{
Rappelons qu'avec les notations que l'on vient de mettre en place, se donner une orientation partielle c'est orienter les arêtes marquées a. Se donner une orientation partielle admissible c'est se donner une orientation partielle, telle qu'il existe une orientation des arêtes marquées h telle que l'arbre ainsi orienté ne possède aucun puits et aucune source.
}
\\

\begin{figure}[!h]
\centerline{\psfig{figure=tetra4c.eps, width=3cm}} \vspace{3em}
\centerline{\psfig{figure=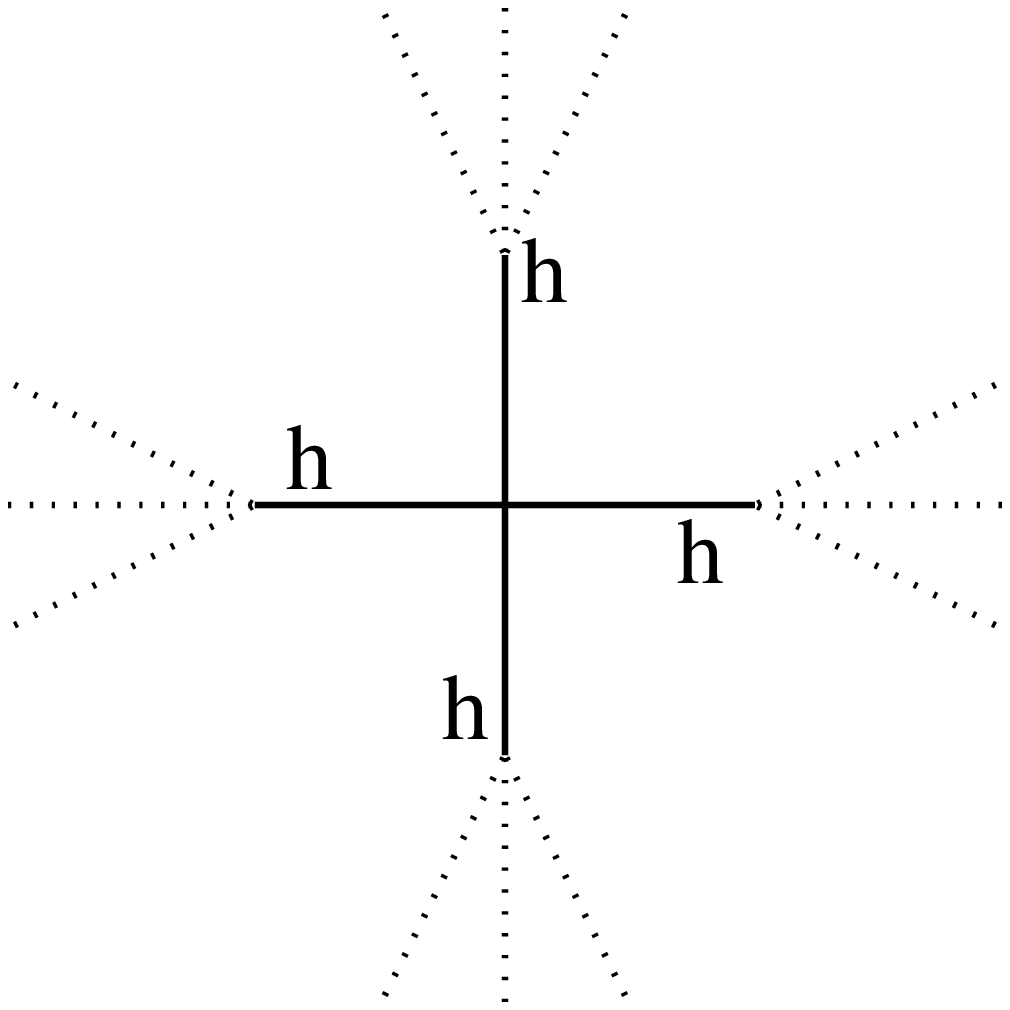, width=3cm}
\hspace{5em} \psfig{figure=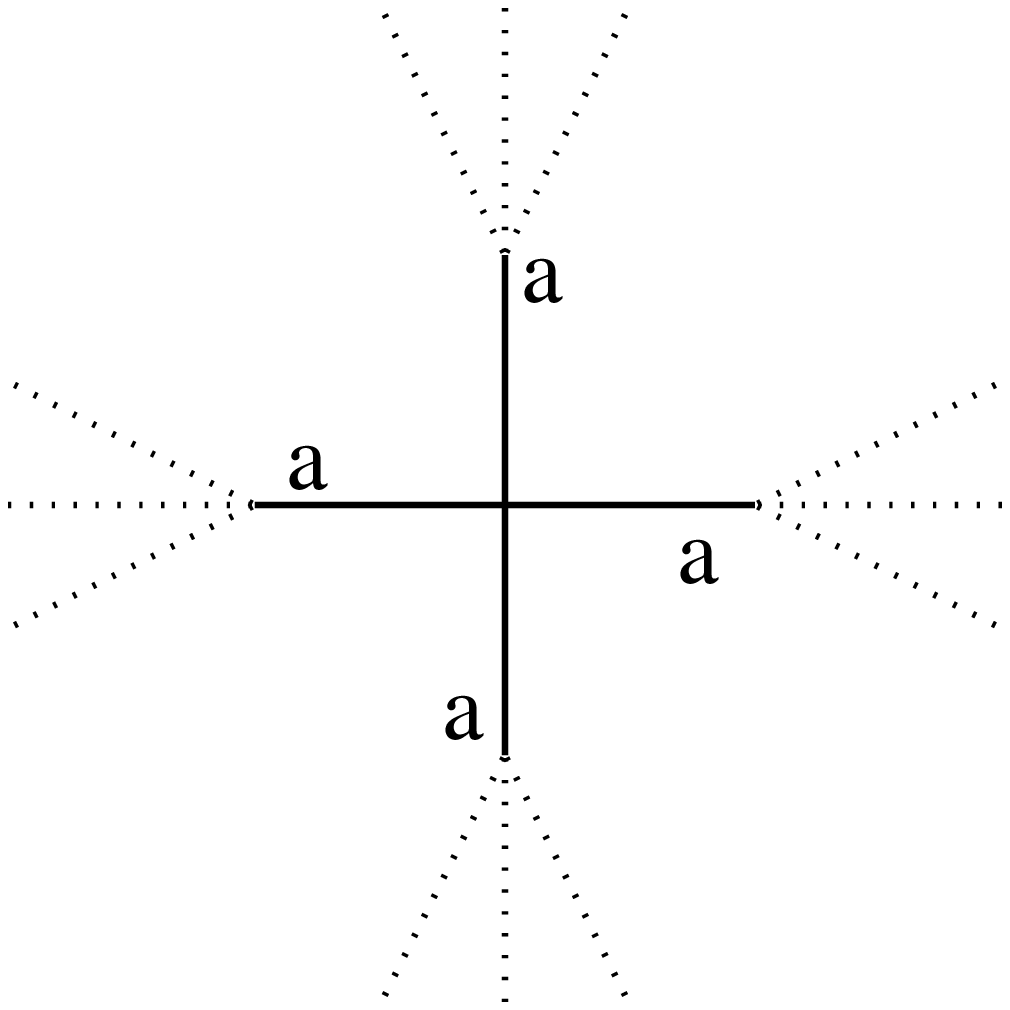, width=3cm}}
\caption{$\kappa_1=1$ et $\kappa_2=14$} \label{last1}
\end{figure}

\par{
1) L'arbre $\A_1$ de la figure \ref{last1} ne possède aucune arête marquée a, par conséquent, il n'y aucune arête à orienter d'où $\kappa_1=1$.
}
\\
\par{
2) L'arbre $\A_2$ de la figure \ref{last1} possède 4 arêtes marquées a, par conséquent, il y a $2^4 = 16$ orientations partielles. Si on oriente toutes ces arêtes vers l'intérieur (resp. l'extérieur) alors on obtient un puits (resp. une source). Il est facile de voir que ce sont les seules orientations partielles non admissibles. On a donc $\kappa_2=16-2=14$.
}
\\

\begin{figure}[!h]
\centerline{\psfig{figure=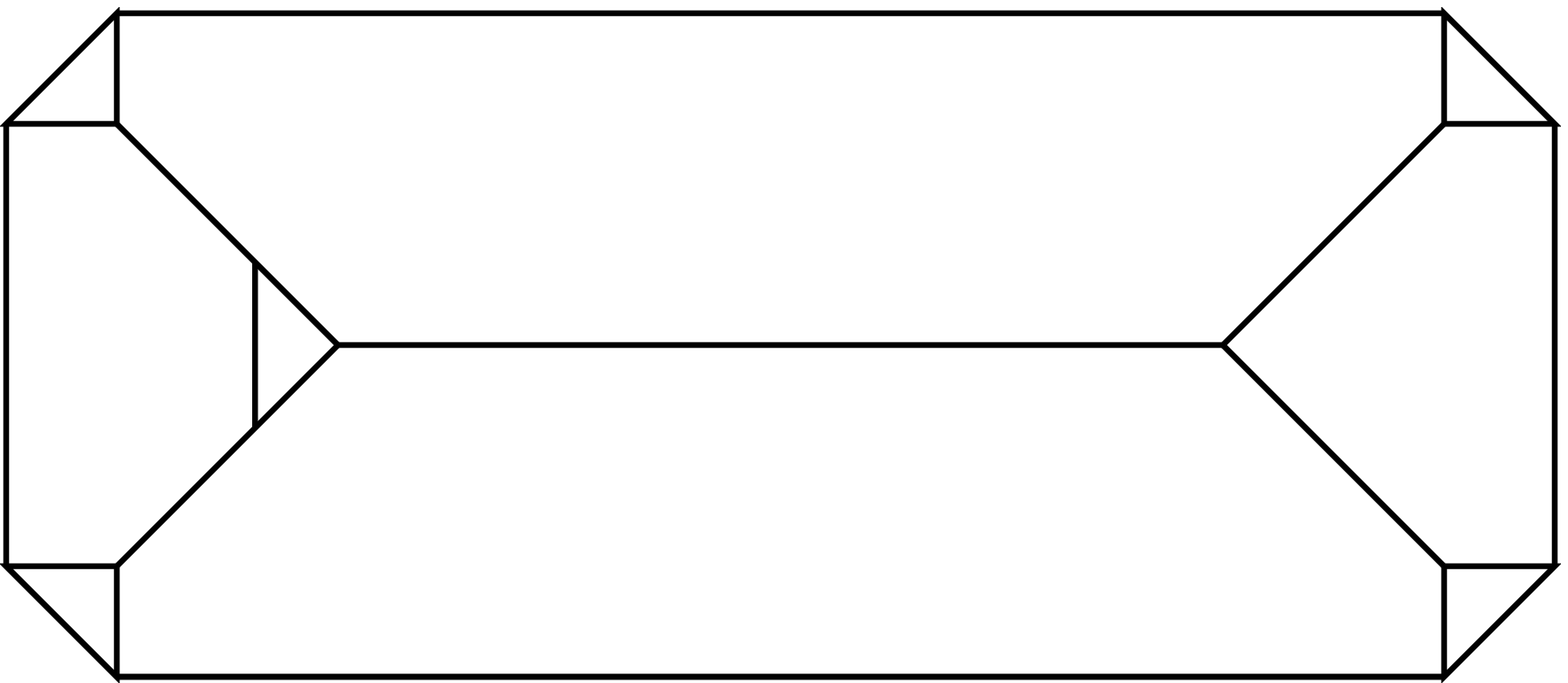, width=4cm}} \vspace{3em}
\centerline{\psfig{figure=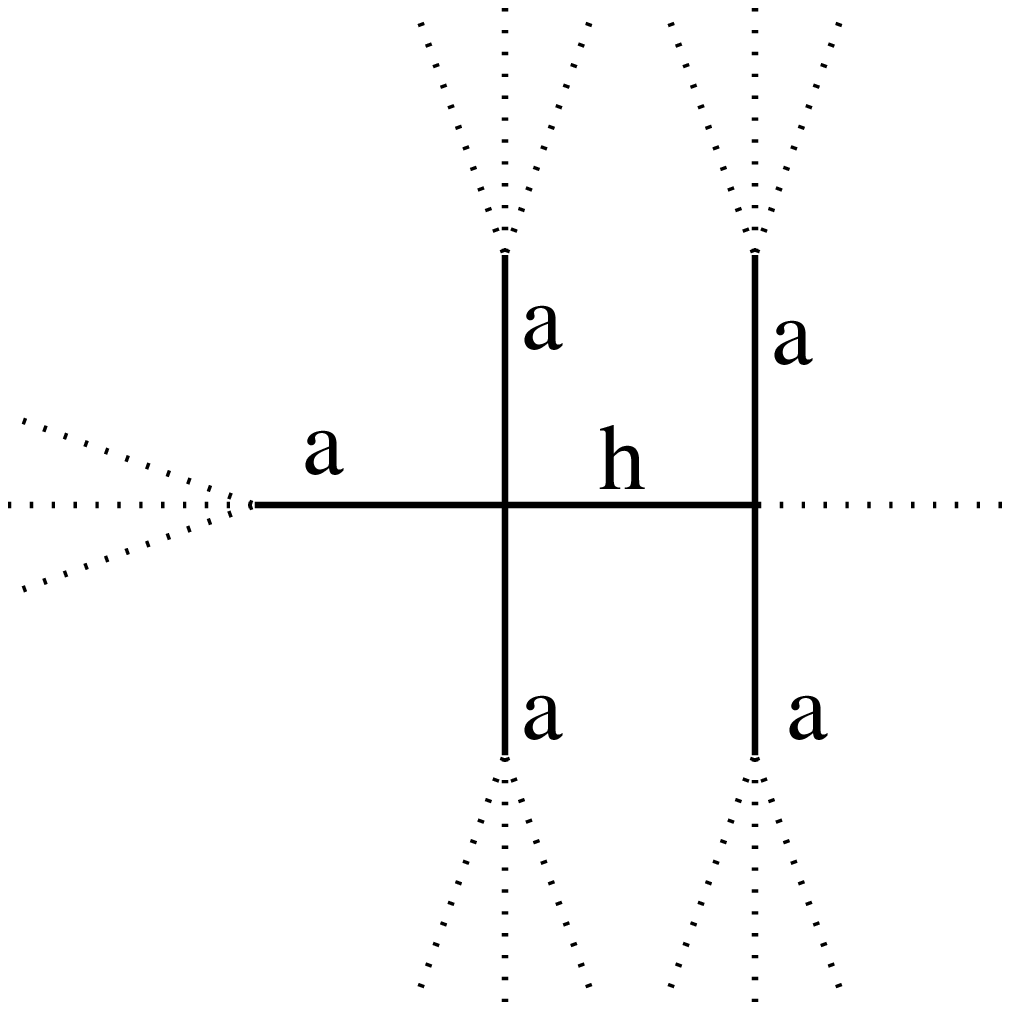, width=3cm} \hspace{5em}
\psfig{figure=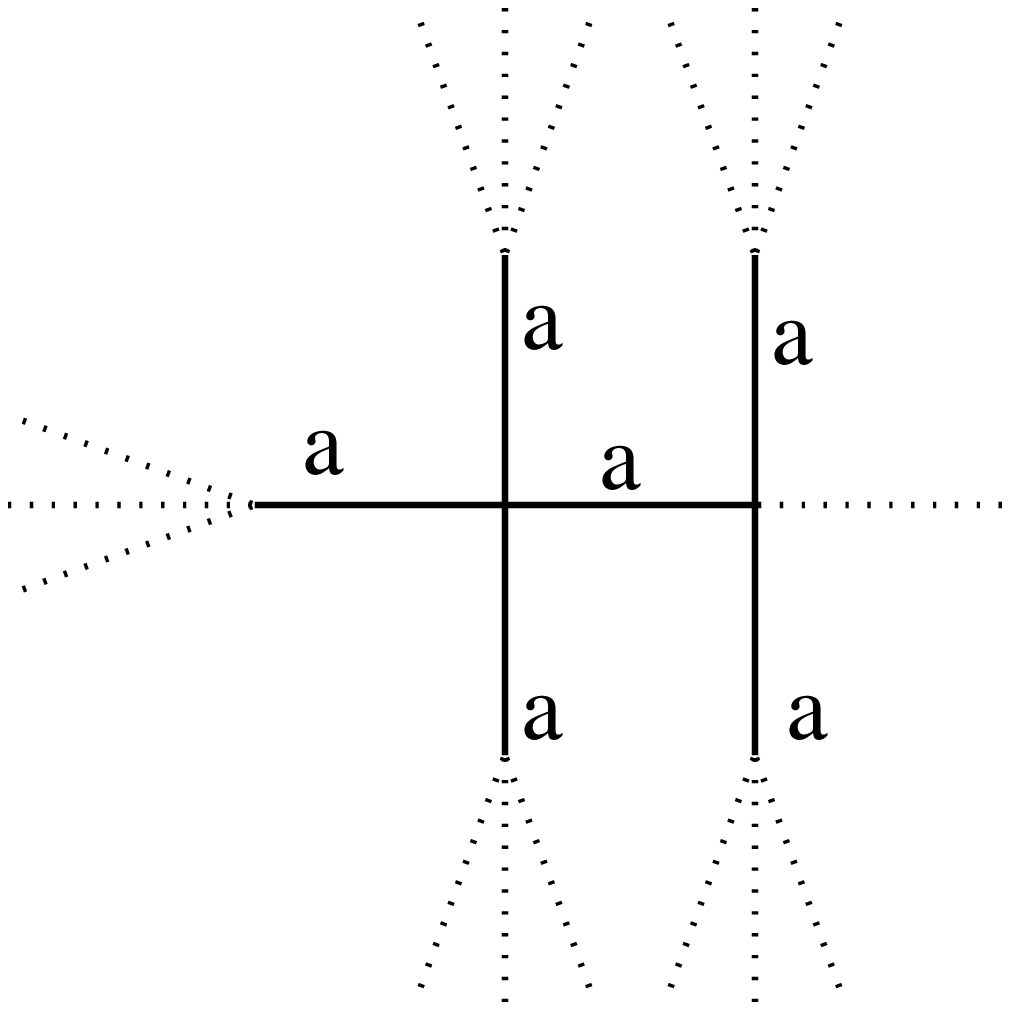, width=3cm}} \caption{$\kappa_1 =
2^5 = 32$ et $\kappa_2=14 \cdot 2^2 = 56$} \label{last2}
\end{figure}

\par{
1) L'arbre $\A_1$ de la figure \ref{last2} possède 5 arêtes marquées a, par conséquent, il y a $2^5 = 32$ orientations partielles. Il est facile de voir que toutes ces orientations partielles sont admissibles. On a donc $\kappa_1=2^5=32$.
}
\\
\par{
2) L'arbre $\A_2$ de la figure \ref{last2} possède 6 arêtes marquées a, par conséquent, il y a $2^6 = 64$ orientations partielles. Mais elles ne sont pas toutes admissibles. Se donner une orientation partielle admissible de $\A_2$, c'est se donner une orientation partielle admissible du sous-arbre formé des 4 arêtes marquées a qui ont une sommet en commun, puis se donner n'importe quelle orientation pour les 2 arêtes restantes. On a donc $\kappa_2 = 14 \cdot 2^2 = 56$.
}
\\

\begin{figure}[!h]
\centerline{\psfig{figure=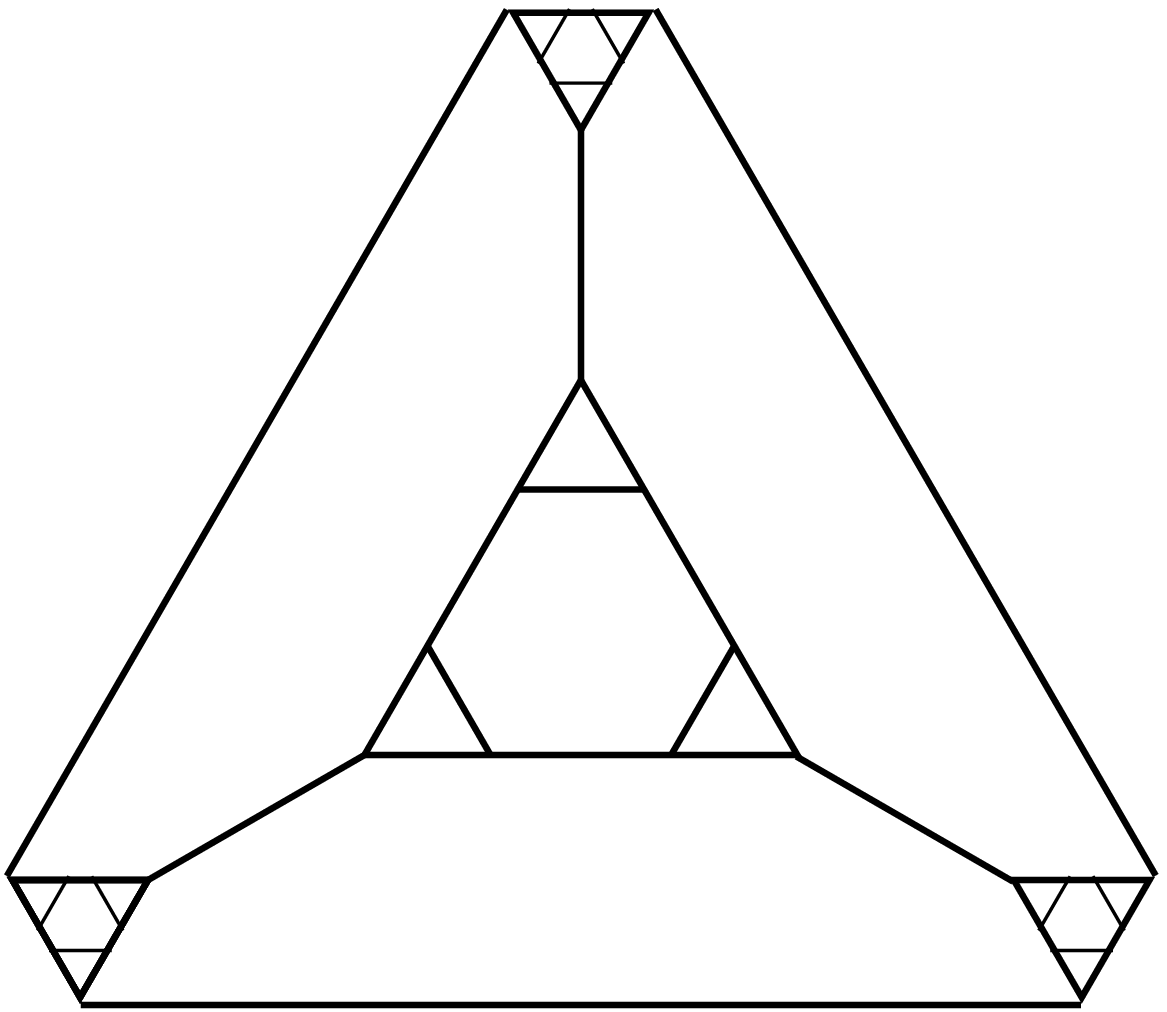, width=4cm}} \vspace{3em}
\centerline{\psfig{figure=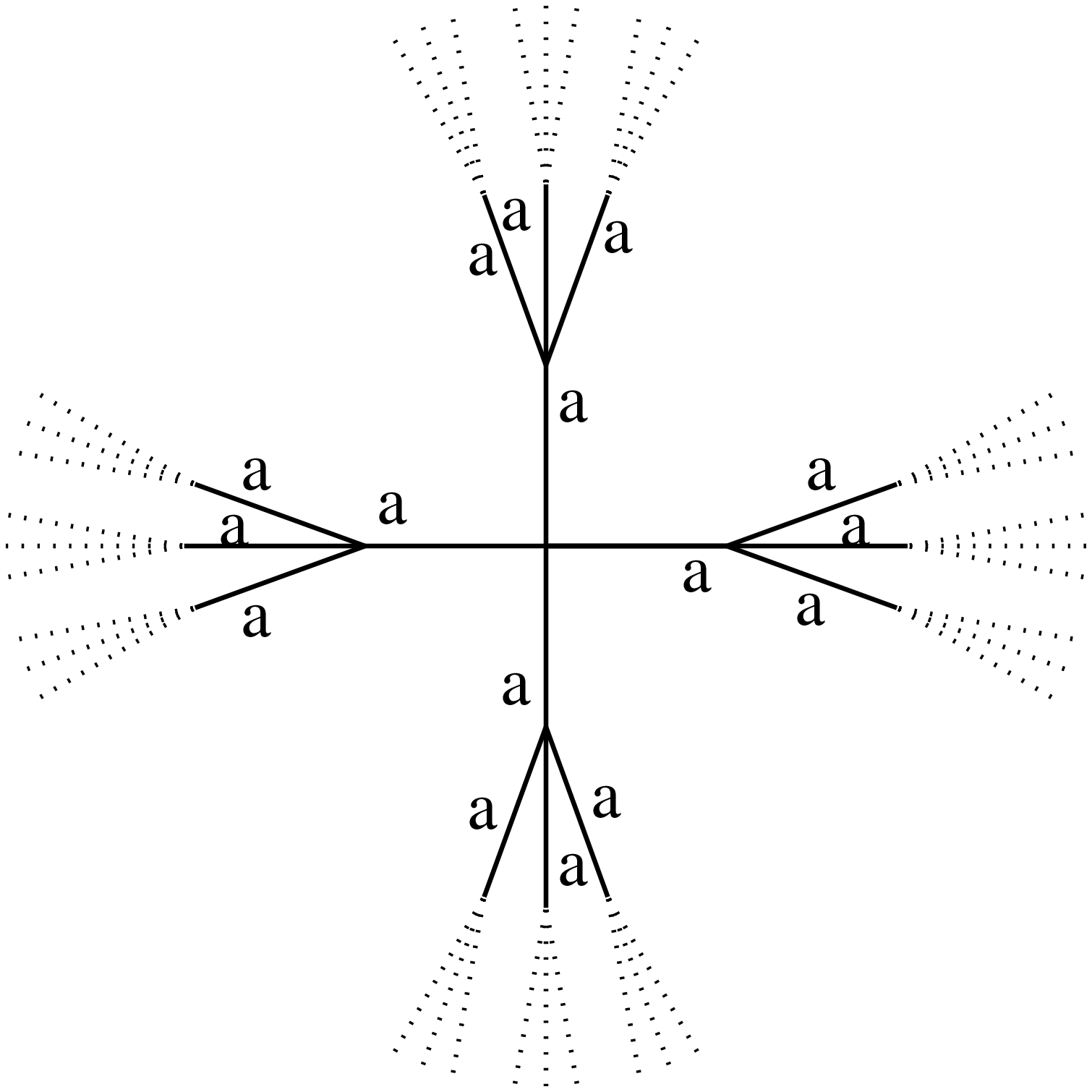, width=4cm} \hspace{5em}
\psfig{figure=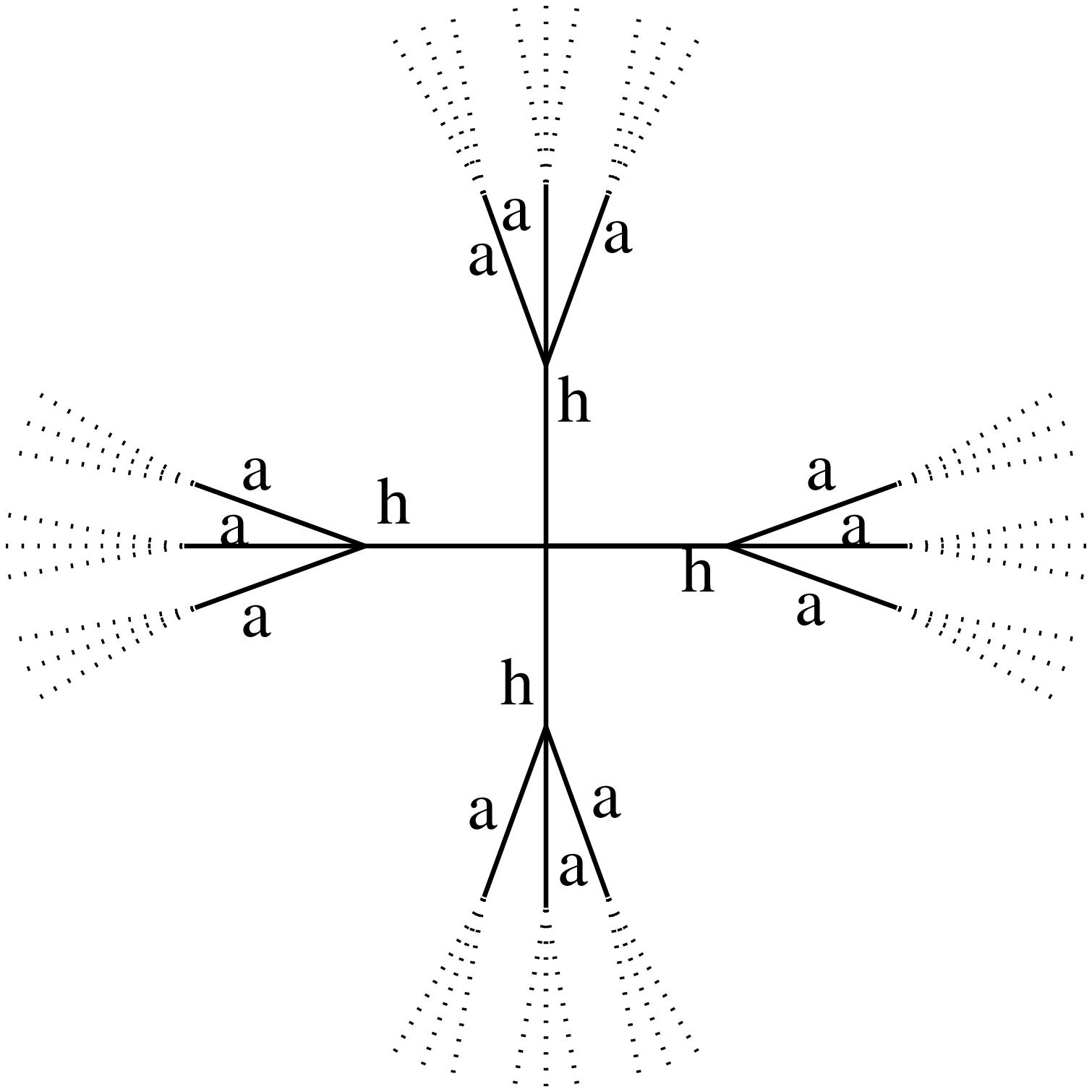, width=4cm}} \caption{$\kappa_1=14
\cdot 7^4 = 33\,614$ et $\kappa_2 = 2^{12} -2 = 4\,094$} \label{last3}
\end{figure}

\par{
1) L'arbre $\A_1$ de la figure \ref{last3} possède 16 arêtes marquées a, par conséquent, il y a $2^{16}=65 \, 536$ orientations partielles. Mais elles ne sont pas toutes admissibles. Se donner une orientation partielle admissible de $\A_1$, c'est se donner une orientation partielle admissible du sous-arbre formé des 4 arêtes marquées a qui ont un sommet en commun, au centre  de $\A_1$, puis se donner l'une des $7$ orientations partielles admissibles pour chacun des 4 groupes de 3 arêtes restantes. On a donc $\kappa_1=14
\cdot 7^4 = 33\,614$.
}
\\
\par{
2) L'arbre $\A_2$ de la figure \ref{last3} possède 12 arêtes marquées a, par conséquent, il y a $2^{12}=4 \, 096$ orientations partielles. Mais elles ne sont pas toutes admissibles. En effet, si l'on oriente toutes les arêtes marquées a vers l'intérieur (resp. l'extérieur), alors on doit orienter toutes les arêtes marqués $h$ vers l'extérieur (resp. l'intérieur). Il est facile de voir que ce sont les deux seules orientations partielles non admissibles.
}

\newpage

\begin{thebibliography}{Ben04b}

\bibitem[And70]{And}
E.M. Andreev.
\newblock On convex polyhedra in Lobacevskii spaces.
\newblock {\em Math. USSR Sbornik}, 10:p. 413--440, 1970.

\bibitem[Ben03]{Beno8}
Yves Benoist.
\newblock Convexes divisibles ii.
\newblock {\em Duke Math. Journ.}, 120:p97--120, 2003.

\bibitem[Ben04a]{Beno9}
Yves Benoist.
\newblock Convexes divisibles i.
\newblock {\em in Algebraic groups and arithmetic Tata Inst. Fund. Res. Stud.
  Math.}, 17:p.339--374, 2004.

\bibitem[Ben04b]{Beno7}
Yves Benoist.
\newblock Five lectures in semisimple lie groups, lecture one.
\newblock {\em Un cours à l'école d'été de Grenoble}, 2004.

\bibitem[Ben05]{Beno6}
Yves Benoist.
\newblock Convexes divisibles iii.
\newblock {\em Annales Scientifiques de l'ENS}, 38:p. 793--832, 2005.

\bibitem[Ben06]{Beno2}
Yves Benoist.
\newblock Convexes divisibles iv.
\newblock {\em Invent. Math.}, 164:p.249--278, 2006.

\bibitem[Bou]{Bou}
Nicolas Bourbaki.
\newblock {\em Groupes et algèbres de Lie Chapitres 4, 5 et 6.}

\bibitem[Cho06]{Choi2}
Suhyoung Choi.
\newblock The deformation spaces of projective structures on 3-dimensional
  Coxeter orbifolds.
\newblock {\em Geometriae Dedicata}, 119 (Vol 1):p 69--90, 2006.

\bibitem[Gol90]{Gold1}
William Goldman.
\newblock Convex real projective structures on compact surfaces.
\newblock {\em J. Differential Geom.}, 31:791--845, 1990.

\bibitem[RG96]{JG}
Jürgen Richter-Gebert.
\newblock Realization spaces of polytopes.
\newblock {\em Number 1643 in Lecture notes in mathematics}, 1996.

\bibitem[RHD07]{RHD}
Roland K.~W. Roeder, John~H. Hubbard, and William~D. Dunbar.
\newblock Andreev's theorem on hyperbolic polyhedra.
\newblock {\em Les Annales de l'Institut Fourier}, 57 (3):p. 825--882, 2007.

\bibitem[Vin71]{Vin3}
Ernest~Borisovich Vinberg.
\newblock Discrete linear groups generated by reflections.
\newblock {\em Math, USSR Izvestija}, 5:p. 1083--1119, 1971.

\end{thebibliography}

\end{document}